\documentclass[a4paper,12pt, oneside, reqno]{amsart}
\usepackage{mathrsfs}
\usepackage{amsfonts}
\usepackage{amsmath,amstext,amssymb,amsthm,bm}
\usepackage{mathtools}
\usepackage[polish, english]{babel}
\usepackage{amsxtra}
\usepackage{color}
\usepackage{dsfont}
\usepackage[margin=2cm, centering]{geometry}
\usepackage[bbgreekl]{mathbbol}
\usepackage{soul}
\usepackage{graphics}
\usepackage{accents}
\usepackage[colorlinks,citecolor=blue,urlcolor=blue,bookmarks=true]{hyperref}

\usepackage{makecell}

\usepackage{float}

\usepackage[leqno]{amsmath}

\makeatletter
\newcommand{\leqnomode}{\tagsleft@true}
\newcommand{\reqnomode}{\tagsleft@false}
\makeatother

\usepackage{multicol,calc}

\newenvironment{nindex}[1][51pt]%
    {\begin{list}{}%
        {%
            \setlength{\labelwidth}{#1}%
            \setlength{\leftmargin}{\labelwidth+\labelsep}%
            \setlength{\itemsep}{2.5pt}%
            \setlength{\parsep}{0pt}%
            \setlength{\rightmargin}{0pt}%
        }%
    }%
    {\end{list}}

\DeclareSymbolFontAlphabet{\mathbb}{AMSb}
\DeclareSymbolFontAlphabet{\mathbbl}{bbold}

\allowdisplaybreaks

\renewcommand{\leq}{\leqslant}
\renewcommand{\geq}{\geqslant}

\theoremstyle{plain}
\newtheorem{theorem}{Theorem}[section] 
\newtheorem{lemma}[theorem]{Lemma}
\newtheorem{proposition}[theorem]{Proposition}
\newtheorem{corollary}[theorem]{Corollary}
\newtheorem{fact}[theorem]{Fact}
 \theoremstyle{definition}
\newtheorem{definition}{Definition}

\newtheorem{example}{Example} 
\newtheorem*{remark*}{Remark} 
\newtheorem{remark}[theorem]{Remark}

\def \NN{{\mathbb N}}
\def \RR{{\mathbb R}}
\def \Rd{\RR^d}
\def \ind{{\bf 1}}
\def \pL{{\mathcal L}}

\newcommand{\rr}{\Upsilon}
\newcommand{\err}[2]{\rho_{#1}^{#2}} 
\newcommand{\errx}[2]{\bar{H}_{#1}^{#2}}
\newcommand{\erry}[2]{\bar{\bar{H}}_{#1}^{#2}}
\def \exdrf{b}

\def \efdrf{b}

\def \lmcc{c_{\!\scriptscriptstyle J}}
\def \ka{c_\kappa}

\def \exda{c_{\exdrf}}

\def \lah{\alpha_h}
\def \uah{\beta_h}
\def \lch{C_h}
\def \khe{\epsilon_{\kappa}}
\def \exdhe{\epsilon_b}
\def \indcs{\sigma}
\def \efcs{\sigma}

\def \indTr{\mathfrak{t_0}}
\def \indnj{N_0}
\def \efnj{N_0}
\newcommand{\indhei}[1]{\epsilon_{#1}}
\newcommand{\indcsi}[1]{s_{#1}}
\def \param{\varpi}
\def \boost{\eta}

\def \aR{\rm R0}
\def \dM{\theta}

\def \Sa{(\mathbf{A\!^\ast})}
\def \Sb{(\mathbf{A})}

\def \pf{\noindent{\bf Proof.} }
\def \qed{{\hfill $\Box$ \bigskip}}

\definecolor{ks}{rgb}{0.7,0.1,0.2}

\title{Non-symmetric L{\'e}vy-type operators}

\author[J. Minecki]{Jakub Minecki}
\address{Jakub Minecki\\
	Wydzia{\l} Matematyki,
	Politechnika Wroc{\l}awska\\
	Wyb. Wyspia\'{n}skiego 27\\
	50-370 Wroc{\l}aw\\
	Poland
}
\email{jakub.minecki@pwr.edu.pl}

\author[K. Szczypkowski]{Karol Szczypkowski}
\address{
	Karol Szczypkowski\\
	Wydzia{\l} Matematyki,
	Politechnika Wroc{\l}awska\\
	Wyb. Wyspia\'{n}skiego 27\\
	50-370 Wroc{\l}aw\\
	Poland
}
\email{karol.szczypkowski@pwr.edu.pl}

\date{}

\begin{document}

\begin{abstract}
We present a general  approach to the parametrix construction.
We apply it to prove the uniqueness and existence of a weak fundamental solution  for the equation $\partial_t u =\pL u$ with non-symmetric non-local operators
$$
\mathcal{L}f(x):= b(x)\cdot \nabla f(x)+ \int_{\mathbb{R}^d}( f(x+z)-f(x)- \ind_{|z|<1} \left<z,\nabla f(x)\right>)\kappa(x,z)J(z)\, dz\,,
$$
under certain assumptions on $b$,  $\kappa$ and $J$.
The result allows  more general coefficients even for
 $J(z)=|z|^{-d-1}$. 
\end{abstract}

\maketitle

\noindent {\bf AMS 2010 Mathematics Subject Classification}: Primary 60J35, 47G20; Secondary 
47D06, 47A55.

\noindent {\bf Keywords and phrases:}  heat kernel estimates,  
 L\'evy-type operator, non-symmetric operator, non-local operator, non-symmetric Markov process, Feller semigroup, Levi's parametrix method.

\section{Introduction}\label{sec:int}

An operator $\mathscr{L}$ is said to be
{\it L{\'e}vy-type} if 
it acts on
every smooth compactly supported function~$f$
according to the following formula
\begin{align*}
\mathscr{L}f(x) &= c(x)f(x)+
\exdrf(x)\cdot \nabla f(x)+\sum_{i,j=1}^d a_{ij}(x) \frac{\partial^2 f}{\partial x_i \partial x_j}(x) \\
&\qquad +\int_{\Rd}\Big(f(x+z)-f(x)- \ind_{|z|<1} \left<z,\nabla f(x)\right>\! \Big) N(x,dz)\,.
\end{align*}
Here $c(x)$, $b(x)$, $a_{ij}(x)$ and $N(x,dz)$
are called  coefficients, and need to satisfy certain natural conditions.
The coefficients model the  infinitesimal behaviour 
of a particle
at the point $x$.
For instance, the vector $b(x)$ defines the drift, while
$N(x,B)$ is the intensity of jumps from $x$ to the set $x+B\subset \Rd$.

The case of constant coefficients, i.e., when
$c(x)\equiv c$, $b(x)\equiv b$, $a_{ij}(x)\equiv a_{ij}$ and $N(x,dz)\equiv N(dz)$,
leads to convolution semigroups of operators or L{\'e}vy processes, which are prevalent
in probability, PDEs, physics, finance and statistics \cite{MR1833689}, \cite{WOS:000266711700004}.
In this case there is a well-established
one-to-one correspondence between
the  process,
the semigroup
and the operator 
$\mathscr{L}$
\cite{MR3185174}.

It is extremely important to 
understand operators with $x$-dependent coefficients. Due to the
Courr{\`e}ge-Waldenfels theorem,
the infinitesimal generator of a Feller semigroup with a sufficiently rich domain is a
L{\'e}vy-type operator;
see \cite[Theorem~2.21]{MR3156646}, \cite[Theorem~4.5.21]{MR1873235}.
However, it is  a highly non-trivial task to construct an operator semigroup for a given L{\'e}vy-type operator
with rough (non-constant) coefficients.
To resolve the problem many authors follow a scheme
known as {\it the parametrix method},
whose primary role is to provide a candidate for the integral kernel of the semigroup in question
(we list references later on in this section).
Each usage of the method
brings about different technical difficulties
to overcome that
depend on the
class of coefficients under consideration,
but most applications 
follow a characteristic pattern, e.g.,
the decomposition
of the candidate kernel
into
{\it a zero order approximation} and {\it a remainder}.

In the paper we provide a general
functional analytic framework 
for the parametrix 
method, namely, 
for the
construction
of the family of operators $\{P_t\colon t\in(0,1]\}$
for a given choice of
{\it the zero order approximation}
operator
$P_t^0$
and {\it the error term} operator $Q_t^0$.
The aforementioned remainder corresponds then to  an operator
expressed
by
$P_t^0$
and
multiple compositions of
$Q_t^0$, see
\eqref{def:op_P}--\eqref{def:op_Qn}.
In Section~\ref{sec:gen}
we single out natural hypotheses
on $P_t^0$ and
$Q_t^0$
that lead to the construction and basic properties of  $P_t$. We furthermore point out key hypotheses
on {\it the approximate solution} 
$P_{t,\varepsilon}$
that validate the semigroup property, non-negativity, etc., of $P_t$. 
We also 
discuss in a general context integral kernels associated with the constructed operators.
Thus, the success of
the parametrix method boils down to
 verifying
the proposed hypotheses.
It should be noted that they may only hold if $P_t^0$ is well chosen.
In
the literature
there is enough evidence  showing that the flexibility of the choice of the zero order approximation
is crucial.
Our results in
Section~\ref{sec:gen}
indeed
provide
such flexibility.
We emphasise that Section~\ref{sec:gen}
(together with Section~\ref{sec:a-proofs}) is general, self-contained and independent of the rest of the paper.

As an application of the framework developed in Section~\ref{sec:gen},
we consider the classical choice
of~$P_t^0$ by
{\it freezing coefficients at the endpoint}.
In doing so we focus on the case
 when
the measure $N(x,dz)$ is non-symmetric in $dz$.
In particular, a non-zero {\it internal drift} 
$$
\zeta_r^x =\int_{\Rd} z \left(\ind_{|z|<r}-\ind_{|z|<1}\right) N(x,dz)\,,\qquad r>0,
$$
may be induced by non-symmetric jumps.
For example, for
$c\equiv 0$, $b\equiv 0$, $a_{ij} \equiv 0$, we get
\begin{align*}
\mathscr{L}f(x)= 
\int_{\Rd}\Big(f(x+z)-f(x)- \ind_{|z|<r} \left<z,\nabla f(x)\right>\! \Big) N(x,dz)\\
+\left( \int_{\Rd} z \left(\ind_{|z|<r}-\ind_{|z|<1}\right) N(x,dz)\right) \cdot \nabla f(x)\,.
\end{align*}
The above can be interpreted
as a decomposition 
of the operator
into
{\it the leading non-local part}
and
{\it the internal drift part},
and we wish to control the latter.
As observed in examples,
the internal drift
$\zeta_r^x$
is much more difficult to handle than the {\it external drift} $b(x)$. 
To be more specific,
in Section~\ref{sec:appl}  we study the operator
\begin{align}\label{def:operator}
\pL f(x)=\exdrf(x)\cdot \nabla f(x)+ \int_{\Rd}\Big(f(x+z)-f(x)- \ind_{|z|<1} \left<z,\nabla f(x)\right>\! \Big) \,\kappa(x,z)J(z)dz\,.
\end{align}
Under the  assumptions
on
$\exdrf,J,\kappa$
specified in Section~\ref{sec:set}
we prove 
our main results collected in Section~\ref{sec:main}. This includes
the uniqueness and existence of
a weak fundamental solution for the equation $\partial_t u=\pL u$.
We also analyse the semigroup  associated with the solution
and discuss properties of its generator, which we identify as the closure of the operator 
$\pL$ in
\eqref{def:operator} acting on the space of smooth compactly supported functions. Pointwise estimates of the fundamental solution are established under additional conditions.

Our main focus is  the case when the mapping 
$z\mapsto \kappa(x,z)J(z)$
is non-symmetric.
A surprising fact is that despite extensive study of non-local operators
and rapidly growing literature
on this subject,
the following fundamental example has not yet been covered. 
Let $\exdrf=0$ and $J(z)=|z|^{-d-1}$ in \eqref{def:operator}, i.e.,
\begin{align}\label{def:1-stable}
\pL f(x)= \int_{\Rd}\Big(f(x+z)-f(x)- \ind_{|z|<1} \left<z,\nabla f(x)\right>\! \Big) \,\frac{\kappa(x,z)}{|z|^{d+1}} dz\,,
\end{align}
and suppose only that $\ka^{-1} \leq \kappa(x,z)\leq \ka$ and 
$|\kappa(x,z)-\kappa(y,z)|\leq \ka |x-y|^{\khe}$
for some 
$\ka>0$, 
$\khe \in (0,1]$ and all $x,y,z\in\Rd$.
We note that the 
results
available
in 
\cite{KKS}, \cite{MR4680974}, \cite{PJ}, \cite{MR3806688}
require further  assumptions on
the coefficient
$\kappa$ to treat~\eqref{def:1-stable}, which exclude some natural 
examples. Our results remove those restrictions
and we construct and estimate the semigroup $(P_t)_{t>0}$, see Example~\ref{ex:1}.
Of course we also
 cover other interesting operators, see Section~\ref{sec:ex}.
We emphasise that the general functional analytic framework developed in Section~\ref{sec:gen} should apply to other zero order approximations $P_t^0$ \cite{MR3907007}, but even for the choice of $P_t^0$ resulting from freezing coefficients at the endpoint we obtain new results for the non-symmetric case.

We now comment on the historical development of the parametrix method.
It was proposed by
Levi~\cite{zbMATH02644101}, Hadamard \cite{zbMATH02629781}, Gevery \cite{zbMATH02629782}
for differential operators and later extended by 
 Feller \cite{zbMATH03022319}  and
Dressel \cite{MR0003340}.
Here, we refer the reader to the classical monographs
by Friedman \cite{MR0181836} and Eidelman \cite{MR0252806}.
Non-local operators were later treated by
 Drin'~\cite{MR0492880}, Drin' and Eidelman \cite{MR616459}, Kochubei \cite{MR972089} and Kolokoltsov \cite{MR1744782}, see also the monograph by Eidelman, Ivasyshin, and Kochubei \cite{MR2093219}.
More recently, the method was used for non-local operators by 
Xie and Zhang \cite{MR3294616},
Knopova and Kulik \cite{MR3652202}, 
Bogdan, Knopova, and Sztonyk~\cite{MR4077549},
K{\"u}hn \cite{MR3912204},
Knopova, Kulik, and Schilling \cite{KKS},
 Chen and Zhang \cite{MR3500272, MR3806688}, Kim, Song, and Vondra\v{c}ek \cite{MR3817130}, Grzywny and Szczypkowski \cite{MR3996792}, and Szczypkowski~\cite{MR4680974}.
It was also employed in the analysis of
SDEs by Kohatsu-Higa and Li \cite{MR3544166},
Knopova and Kulik \cite{MR3765882}, Kulik \cite{MR3907007},
Kulczycki and Ryznar \cite{MR4167204}, 
Kulczycki, Ryznar, and Sztonyk \cite{MR4261305}, Kulczycki, Kulik, and Ryznar
\cite{KKR-2020},
Menozzi and Zhang \cite{MR4394312}.
We further refer the reader to
more recent monographs by
Knopova, Kochubei, and Kulik \cite{MR3965398},
Chen and Zhang \cite{MR3824209}, 
and
K\"{u}hn \cite{MR3701414}.

A different Hilbert-space approach to the parametrix method was developed by Jacob \cite{MR1254818, MR1917230}, Hoh \cite{MR1659620}, B{\"o}ttcher \cite{MR2163294, MR2456894}, Tsutsumi \cite{MR0367492, MR0499861} and Kumano-go  \cite{MR666870}.

For a more complete list of references,
we refer to Knopova, Kulik and Schilling \cite{KKS}
and Bogdan, Knopova, and Sztonyk \cite{MR4077549}.

The  paper is organized as follows.
Section~\ref{sec:big} is divided into three parts.
In Section~\ref{sec:set} we specify our setting and assumptions.
Section~\ref{sec:main} starts with a surprisingly general uniqueness result, after which we formulate the main existence result.
In Section~\ref{sec:ex} we give informative examples of applications.
Section~\ref{sec:gen} (independent and self-contained) is devoted to the general presentation of the parametrix method via
a functional analytic approach. 
In Section~\ref{sec:appl} the method is applied
to investigate the operator \eqref{def:operator} under  assumptions on the coefficients of $\pL$ described in Section~\ref{sec:set},
and we give proofs of our main results on this operator.
For a more compact presentation, the
statements of
Section~\ref{sec:f_a_a}
are proved
in Section~\ref{sec:a-proofs}.
Further technical results, that are used in Section~\ref{sec:appl}, are shifted to
Sections~\ref{sec:sac} and~\ref{sec:a-frozen}.

\subsection*{General notation}

We denote by $d$ a fixed (natural) dimension, that is, $d\in \NN=\{1,2,\ldots\}$.
By $c(d,\ldots)$ we mean a
positive number that depends only on the listed parameters $d,\ldots$. 
We use ``$:=$" as definition. As usual, $a\land b:=\min\{a,b\}$ and $a\vee b := \max\{a,b\}$.
The Beta function is denoted by $B(a,b)$
and the Gamma function by $\Gamma(a)$.
We let $\NN_0=\NN\cup\{0\}$ and
by $\bbbeta\in \NN_0^d$ we mean a multi-index.

We write  $x\cdot y$ or $\left<x,y\right>$ for the standard scalar product of $x,y\in\RR^n$, $n\in\NN$.
For a function $f\colon \RR^n \to\RR$
the supremum norm is denoted by
$$\|f\|_{\infty}=\sup_{x\in \RR^n}|f(x)|\,.$$
We use the following function spaces.
$B_b(\RR^n)$ denotes the bounded Borel measurable functions,
$C(E)$ are the continuous functions on $E\subseteq \RR^n$.
Furthermore, $C_b(E)$,  $C_0(E)$, $C_c(E)$ are subsets of $C(E)$
of functions that are  bounded,
vanish at infinity,
and have compact support, respectively.
We write
$f\in C^k(\RR^n)$ if
the function and all its partial derivatives up to (and including if finite) order $k\in \NN\cup \{\infty\}$ are elements of $C(\RR^n)$;
we similarly understand $C_b^k(\RR^n)$,  $C_0^k(\RR^n)$, $C_c^k(\RR^n)$.
In particular,
$C^{\infty}(\RR^n)$ are smooth functions, while
$C_c^{\infty}(\RR^n)$ are smooth functions with compact support.

We say that a Borel set is of full measure if its complement is of Lebesgue measure zero (the underlying Euclidean space shall be clear from the context).
For a Borel set $\mathcal{S}\subseteq (0,1]\times\Rd\times\Rd$ and $t\in (0,1]$, $x\in\Rd$, we consider the section $\mathcal{S}^{t,x}:=\{y\in\Rd\colon (t,x,y)\in \mathcal{S}\}$. The latter is a Borel set, see \cite[Lemma~1.28]{MR4226142}.

We refer to $\RR\times\Rd$ as the {\it space-time}.
Each of the sets $(a,b)$, $[a,b)$, $(a,b]$, $[a,b]$, where $-\infty<a<b<\infty$, is called a {\it bounded interval}.

\section*{List of important notation}

For the reader's convenience we collect important constants, functions and conditions.

\begin{table}[H]
\caption{Notation in Section~\ref{sec:gen} and Section~\ref{sec:a-proofs}}
{\footnotesize
{\raggedright
\noindent
\begin{multicols}{3}
\raggedright
\parindent0pt
\begin{nindex}

\item[$\varepsilon_0$] (H3) p.~\pageref{H3}  

\item[$r_t$]  ($\aR$) p.~\pageref{R}

\item[$L$, $D(L)$] (L0) p.~\pageref{L0}

\item[$\mathcal{D}$] (L3) p.~\pageref{L3}

\item[$P^0_t$, $Q^0_t$] \eqref{op:S1} p.~\pageref{op:S1}

\item[$Q_t$] \eqref{def:op_Q} p.~\pageref{def:op_Q}

\item[$Q_t^n$] \eqref{def:op_Qn} p.~\pageref{def:op_Qn}

\item[$P_t$] \eqref{def:op_P} p.~\pageref{def:op_P}

\item[$P_{t,\varepsilon}$] \eqref{def:approx_sol} p.~\pageref{def:approx_sol}

\item[$p_0(t,x,y)$] \eqref{op:S3} p.~\pageref{op:S3}

\item[$q_0(t,x,y)$]  \eqref{op:S3} p.~\pageref{op:S3}

\item[$q_n(t,x,y)$] \eqref{def:q_n-gen} p.~\pageref{def:q_n-gen}

\item[$q(t,x,y)$] \eqref{def:q-gen} p.~\pageref{def:q-gen}

\item[$p(t,x,y)$] \eqref{def:p-gen} p.~\pageref{def:p-gen}

\item[($\aR$)] p.~\pageref{R}

\item[(H0)]  p.~\pageref{H0}

\item[(H1)-(H4)]  p.~\pageref{H1}

\item[(H5), (H6)] p.~\pageref{H5}

\item[(L0)-(L3)] p.~\pageref{L0}

\item[(CoJ1)-(CoJ5)]  p.~\pageref{CoJ1}

\end{nindex}
\end{multicols}
}
}
\end{table}

\vspace{-0,6cm}

\begin{table}[H]
\caption{Notation in Sections~\ref{sec:big} and~\ref{sec:appl}, and Sections~\ref{sec:sac} and~\ref{sec:a-frozen}}
{\footnotesize
{\raggedright
\noindent
\begin{multicols}{3}
\raggedright
\parindent0pt
\begin{nindex}

\item[$\varepsilon_0$] \eqref{def:ve_0-A}, \eqref{def:ve_0-A*} p.~\pageref{def:ve_0-A}

\item[$r_t$] \eqref{def:r_t} p.~\pageref{def:r_t}

\item[$\indTr$] \eqref{def:indTr} p.~\pageref{def:indTr}

\item[$\nu(|x|)$] \eqref{def:nu} p.~\pageref{def:nu}

\item[$\kappa(x,z)$] \eqref{set:k-bound}, \eqref{set:k-holder} p.~\pageref{set:k-bound}

\item[$J(x)$] \eqref{set:J} p.~\pageref{set:J}

\item[$\lah$] \eqref{set:h-scaling} p.~\pageref{set:h-scaling}

\item[$\khe$] \eqref{set:k-holder} p.~\pageref{set:k-holder}

\item[$\indcs$] \eqref{set:indrf-cancellation-scale} p.~\pageref{set:indrf-cancellation-scale}

\item[$\indhei{j}$, $\indcsi{j}$] \eqref{set:indrf-holder-s-B}, \eqref{set:indrf-holder-s-A} p.~\pageref{set:indrf-holder-s-B}

\item[$\indhei{0}$, $\indcsi{0}$] \eqref{def:eps_0_and_s_0} p.~\pageref{def:eps_0_and_s_0}

\item[$\indnj$] \eqref{set:indrf-holder-s-B}, \eqref{set:indrf-holder-s-A} p.~\pageref{set:indrf-holder-s-B}

\item[$\efdrf_r^x$] \eqref{def:efdrf} p.~\pageref{def:efdrf}

\item[$\pL$]  \eqref{def:operator} p.~\pageref{def:operator}

\item[$D(\pL)$] \eqref{def:DpL} p.~\pageref{def:DpL}

\item[$\mathcal{D}$] \eqref{def:pD} p.~\pageref{def:pD}

\item[$\pL^{\mathfrak{K}_w}$] \eqref{def:pL_fr} p.~\pageref{def:pL_fr}

\item[$p^{\mathfrak{K}_w}(t,x,y)$] \eqref{def:p_fr} p.~\pageref{def:p_fr}

\item[$p_0(t,x,y)$] \eqref{def:p0q0} p.~\pageref{def:p0q0}

\item[$q_0(t,x,y)$]  \eqref{def:p0q0} p.~\pageref{def:p0q0}

\item[$P^0_t$, $Q^0_t$] \eqref{def:P0Q0} p.~\pageref{def:P0Q0}

\item[$\rr_t(x)$] \eqref{def:bound_function} p.~\pageref{def:bound_function}

\item[$\errx{\gamma}{\beta}(t,x,y)$] \eqref{def:errx_erry} p.~\pageref{def:errx_erry}

\item[$\erry{\gamma}{\beta}(t,x,y)$] \eqref{def:errx_erry} p.~\pageref{def:errx_erry}

\item[$\delta_{r}^{\mathfrak{K}_w}(t,x,y;z)$] \eqref{def:delta} p.~\pageref{def:delta}

\item[$\Sb$, $\Sa$] p.~\pageref{Sb-Sa}

\item[$\param$] Section~\ref{param} p.~\pageref{param}

\end{nindex}
\end{multicols}
}
}
\end{table}

\vspace{-0,5cm}
\noindent
$Q_t$, $Q^n_t$, $P_t$, $P_{t,\varepsilon}$ and 
$q(t,x,y)$,
$q_n(t,x,y)$,
$p(t,x,y)$
are used in 
Sections~\ref{sec:big} and~\ref{sec:appl}
as defined in Section~\ref{sec:gen},
see the summary at the end of Section~\ref{ssec:a-hc} on page \pageref{P_t}.

\vspace{\baselineskip}

\section{Main applications of the parametrix method}\label{sec:big}

\subsection{Setting and assumptions}\label{sec:set}

In this section we introduce 
two sets of assumptions $\Sb$ and $\Sa$ under which the main results in Section~\ref{sec:main} are stated.
We begin by defining objects and formulating conditions.
All the functions considered in the paper are assumed to be Borel measurable on their natural domains.

Let $d\in\NN$ and
$\nu:[0,\infty)\to[0,\infty]$ be a non-increasing  function ($\nu \not\equiv 0$) satisfying
\begin{align}\label{def:nu}
\int_{\Rd}  (1\land |x|^2)\, \nu(|x|)dx<\infty\,.
\end{align}
For $r>0$, we let 
\begin{align*}
h(r):= \int_{\Rd} \left(1\land \frac{|x|^2}{r^2}\right) \nu(|x|)dx\,,\qquad \quad
K(r):=r^{-2} \int_{|x|<r}|x|^2 \,\nu(|x|)dx\,.
\end{align*}
We say that  {\it the weak  scaling condition} at the origin holds if there are $\lah\in (0,2]$ and $\lch \in [1,\infty)$ such that 
\begin{align}\label{set:h-scaling}
h(r)\leq \lch \,\lambda^{\lah} h(\lambda r)\,,\qquad \lambda,\, r\in (0, 1]\, .
\end{align}
We refer the reader to \cite[Lemma~2.3]{MR4140542}  for other equivalent descriptions of \eqref{set:h-scaling},
which involve relations between $h$ and $K$.
We consider $J: \Rd  \to [0, \infty]$
 such that for some constant
$\lmcc\in [1,\infty)$ and all $x\in\Rd$,
\begin{align}\label{set:J}
\lmcc^{-1} \nu(|x|) \leq J(x) \leq  \lmcc\, \nu(|x|)\,.
\end{align}
Next, suppose that $\kappa \colon \Rd\times\Rd \to (0,\infty)$ is such that for a constant $\ka\in[1,\infty)$ and all $x,z\in\Rd$,
\begin{align}\label{set:k-bound}
\ka^{-1} \leq \kappa(x,z)\leq \ka\,,
\end{align}
and for some $\khe \in (0,1]$ and all $x,y,z\in\Rd$,
\begin{align}\label{set:k-holder}
|\kappa(x,z)-\kappa(y,z)|\leq \ka \,|x-y|^{\khe}\,.
\end{align}
Let  $\exdrf\colon \Rd \to \Rd$. For $x\in\Rd$ and $r>0$, we define (the effective drift)
\begin{align}\label{def:efdrf}
\efdrf_r^x:=\exdrf(x)+\int_{\Rd} z \left(\ind_{|z|<r}-\ind_{|z|<1}\right) \kappa(x,z)J(z)dz\,.
\end{align}
We ponder the existence of
a constant  $\indcs\in (0,1]$
such that for all $x\in\Rd$,
\begin{align}\label{set:indrf-cancellation-scale}
|\efdrf_r^x| \leq \ka r^{\indcs} h(r)\,, \qquad  r\in (0,1]\,.
\end{align}
Furthermore, we consider $\indnj\in\NN$,
$\indhei{j}\in (0,1]$ and $\indcsi{j} \in (0,1]$ such that
\begin{align}
\label{set:indrf-holder-s-B}
|\efdrf_r^x-\efdrf_r^y|\leq \ka
\sum_{j=1}^{\indnj}  |x-y|^{\indhei{j}} \, r^{\indcsi{j}} h(r)\,, \qquad |x-y|\leq 1,\,  r\in (0,1]\,.
\end{align}
Alternatively to \eqref{set:indrf-holder-s-B}, we consider the assumption
\begin{align}
\tag{\ref*{set:indrf-holder-s-B}${}^\ast$}
\label{set:indrf-holder-s-A}
|\efdrf_r^x-\efdrf_r^y|\leq \ka
\sum_{j=1}^{\indnj}  (|x-y|^{\indhei{j}}\land 1) \, r^{\indcsi{j}} h(r)\,, \qquad  x,y\in\Rd,\,r\in (0,1]\,.
\end{align}

The most natural example here is the isotropic $\alpha$-stable L{\'e}vy measure $J(x)=\nu(|x|)=|x|^{-d-\alpha}$, $\alpha\in (0,2)$. In this case $h(r)=r^{-\alpha}h(1)$, and depending on the values of the parameters $\indcs, \indcsi{j}\in (0,1]$, the right hand side of \eqref{set:indrf-cancellation-scale},
\eqref{set:indrf-holder-s-B} or
\eqref{set:indrf-holder-s-A} may be unbounded for small $r$, see also Example~\ref{ex:1}.

We now specify our framework. In what follows, the {\it dimension} $d$ and {\it the profile function} $\nu$ are always as described at the beginning of this section. 
We will alternatively use two sets of assumptions:

\vspace{\baselineskip}
\label{Sb-Sa}
\begin{itemize}
\item[$\Sb\colon$]  \eqref{set:h-scaling}--\eqref{set:k-holder} hold;
\eqref{set:indrf-cancellation-scale}
and \eqref{set:indrf-holder-s-B} hold with
\begin{align*}
\underset{j=1,\ldots,\efnj}{\forall}\quad 
\lah\land (\efcs \indhei{j})+ \indcsi{j} -1>0
\qquad \mbox{and}\qquad \lah+\efcs-1>0\,.
\end{align*}

\vspace{\baselineskip}

\item[$\Sa\colon$]  \eqref{set:h-scaling}--\eqref{set:k-holder} hold;
\eqref{set:indrf-cancellation-scale}
and \eqref{set:indrf-holder-s-A} hold with
\begin{align*}
\underset{j=1,\ldots,\efnj}{\forall}\quad 
\lah\land (\efcs \indhei{j})+ \indcsi{j} -1>0\,.
\end{align*}
\end{itemize}

\vspace{\baselineskip}

We shall always specify which one of
$\Sb$ and $\Sa$ 
is supposed to
hold.
We conclude with a few  comments.

\begin{remark}\label{rem:exdrf}
If
\eqref{set:indrf-cancellation-scale} holds, then
for some $\exda>0$ and all $x\in\Rd$,
\begin{align}\label{set:exdrf-bound}
|\exdrf(x)|\leq  \exda \,.
\end{align}
If \eqref{set:indrf-holder-s-B} 
or \eqref{set:indrf-holder-s-A} holds, then
for some $\exda>0$, $\exdhe \in (0,1]$ and all $|x-y|\leq 1$,
\begin{align}\label{set:exdrf-holder}
|\exdrf(x)-\exdrf(y)|\leq  \exda |x-y|^{\exdhe}\,.
\end{align}
Conversely, 
assuming \eqref{set:h-scaling},
if \eqref{set:exdrf-bound} and \eqref{set:exdrf-holder} hold, then there is $\ka>0$ such that for all $x,y\in\Rd$, $r\in (0,1]$,
$$
|\exdrf(x)| \leq \ka r^{\lah \land 1}h(r)\,,  \qquad \qquad
|\exdrf(x)-\exdrf(x)| \leq \ka (|x-y|^{\exdhe} \land 1)\, r^{\lah \land 1}h(r)\,.
$$
\end{remark}

\vspace{\baselineskip}

For the sake of the following discussion, we consider another condition: There are $\uah\in (0,2]$ and $c_h\in (0,1]$ such that
\begin{equation}\label{eq:intro:wusc}
 h(r)\geq c_h\,\lambda^{\uah}\,h(\lambda r)\, ,\qquad \lambda, r\in(0,1]\, .\\
\end{equation}
We relate our assumptions with those in \cite{MR3996792} and \cite{MR4680974}, without  addressing specific results in these papers.

\begin{remark}\label{rem:other-as}
Assume that \eqref{set:h-scaling}--\eqref{set:k-holder} hold. In each case below, 
$\Sa$  is satisfied, because
\eqref{set:indrf-cancellation-scale}
and \eqref{set:indrf-holder-s-A} hold with indicated parameters.

\begin{enumerate}
\item[(i)] If $\exdrf=0$ and $\kappa(x,z)J(z)=\kappa(x,-z)J(-z)$, then $\indcs=\indnj=\indhei{1}=\indcsi{1}=1$. This is exactly the case {\rm (P3)} in \cite{MR3996792}.
 
\item[(ii)] If  \eqref{set:exdrf-bound}--\eqref{set:exdrf-holder} hold and $\lah>1$, then $\indcs=\indnj=1$ and $\indhei{1}=\khe \land \exdhe$, $\indcsi{1}=1$, see \cite[Fact~1.1]{MR4680974}.
In particular, it covers the case {\rm (P1)} in \cite{MR3996792}.

\item[(iii)] If \eqref{eq:intro:wusc} holds, $0<\lah\leq \uah <1$ and
$\exdrf(x)=\int_{|z|<1} z\kappa(x,z) J(z)dz$, then
$\indcs=\indnj=1$ and $\indhei{1}=\khe$, $\indcsi{1}=1$. This is exactly the case {\rm (P2)} in \cite{MR3996792}. Note that in this case the operator in \eqref{def:operator},
for a smooth compactly supported function $f$, is equal to the following (absolutely convergent) integral
\begin{align}\label{eq:op-2}
\pL f(x)=\int_{\Rd}(f(x+z)-f(x)) \,\kappa(x,z)J(z)dz\,.
\end{align}

\item[(iv)] If {\rm (Q0)} from \cite{MR4680974} holds, then  $\indcs=\indnj=\indhei{1}=\indcsi{1}=1$. Note that
in \cite{MR4680974} the results were obtained under stronger conditions {\rm (Q1)} and {\rm (Q2)}.

\item[(v)] If $\exdrf=0$ and for all $x\in\Rd$ and $r\in (0,1]$,
\begin{align*}
\int_{r\leq |z|<1} z \kappa(x,z)J(z)dz =0\,,
\end{align*}
then $\indcs=\indnj=\indhei{1}=\indcsi{1}=1$.
The above condition was used in \cite{PJ} to analyse the operator~\eqref{def:1-stable}.
The condition was weakened by proposing 
{\rm (Q1)} in \cite{MR4680974}.

\end{enumerate}
\end{remark}

\begin{remark}\label{rem:close-to-1}
(a) Clearly, if \eqref{set:h-scaling} holds with $\lah\geq 1$ and some constant $\lch>0$, then it holds for every $\lah<1$ with the same constant $\lch$.
(b) Assume that \eqref{set:h-scaling}--\eqref{set:k-holder} and \eqref{set:exdrf-bound}--\eqref{set:exdrf-holder} hold. Then
\begin{enumerate}
\item[(i)] if $\lah\in(0,1)$, then 
\eqref{set:indrf-cancellation-scale}
and \eqref{set:indrf-holder-s-A} hold with $\indcs=\lah$
 and $\indnj=1$, $\indhei{1}=\khe \land \exdhe$, $\indcsi{1}=\lah$.
\item[(ii)] if $\lah\in(0,1)$ can be chosen arbitrarily close to $1$, then $\Sa$  is satisfied.
\end{enumerate}
For justification, see
Remark~\ref{rem:exdrf} and
Lemma~\ref{lem:drf}.
\end{remark}

\vspace{\baselineskip}

\subsection{Main results}\label{sec:main}

In this section, we present our main results.
To formulate them we use the notion of a kernel.
Recall that proofs are deferred to Section~\ref{sec:appl}.

\begin{definition}
We call $\mu$ {\it a kernel} on the space-time
if for every $(s,x)\in \RR\times \Rd$ and every bounded interval $I\subset \RR$
the mapping
$$E \longmapsto \mu(s,x,E)\,,$$
is a (finite) signed measure on $\mathcal{B}(I\times\Rd)$ --  the $\sigma$-field of Borel subsets of $I\times\Rd$.
\end{definition}

In our first result, we address the question of uniqueness of the fundamental solution for $\partial_t u=\pL u$, where $\pL$ is given in~\eqref{def:operator}.

\begin{theorem}\label{thm:uniq}
Assume $\Sb$ or $\Sa$.
Let $(s,x)\in \RR\times\Rd$ be fixed.
Suppose that $\mu_1$ and $\mu_2$ are 
kernels on the space-time
such that
for every $\phi \in C_c^{\infty}(\RR\times\Rd)$
\begin{align}\label{eq:uniq}
\iint\limits_{(s,\infty)\times \Rd}  \mu_j(s,x,dudz)\Big[\partial_u\, \phi(u,z) + \pL_z \phi(u,z)  \Big] = - \phi(s,x)\,,
\qquad j=1,2\,.
\end{align}
Then for every
bounded interval $I\subset (s,\infty)$ and every
set $E\in\mathcal{B}(I \times \Rd)$
we have
$$
\mu_1(s,x,E)=\mu_2(s,x,E)\,.
$$
\end{theorem}

The second result concerns the  question
of existence.
\begin{theorem}\label{thm:exist}
Assume  $\Sb$ or $\Sa$.
There exists a kernel $\mu$  on the space-time such that
for all $(s,x)\in\RR\times\Rd$ and $\phi\in C_c^{\infty}(\RR\times\Rd)$,
\begin{align*}
\iint\limits_{(s,\infty)\times\Rd} \mu(s,x,dudz)
\Big[\partial_u\, \phi(u,z) + \pL_z \phi(u,z)  \Big] = - \phi(s,x)\,.
\end{align*}
Furthermore, there is a Borel function $p\colon (0,\infty) \times\Rd\times\Rd\to \RR$ such that 
for every 
bounded interval $I\subset (s,\infty)$ and every
set $E\in\mathcal{B}(I \times \Rd)$
we have
$$
\mu(s,x,E)=\iint\limits_{E} p(u-s,x,z)\,dzdu\,.
$$
\end{theorem}

\noindent
In Section~\ref{ssec:a-hc} (see the summary on page \pageref{P_t}) for all $t>0$ and $x\in\Rd$, we explicitly define $p(t,x,y)$ of Theorem~\ref{thm:exist}, and in the formulas below, $p(t,x,y)$ means that specific function.
Further properties of $p(t,x,y)$ can be found in
Section~\ref{ssec:ker}.
Let
$$P_t f(x)=\int_{\Rd} p(t,x,y)f(y)\,dy\,.$$

\begin{theorem}\label{thm:sem_prop}
Assume $\Sb$ or $\Sa$.
The family $(P_t)_{t>0}$
is a strongly continuous positive contraction semigroup on $(C_0(\Rd),\|\cdot\|_{\infty})$.
Let $(\mathcal{A},D(\mathcal{A}))$ be its infinitesimal generator.
Then
\begin{itemize}
\item[(i)]  $P_t 1 =1$ for all $t>0$,
\item[(ii)] $P_t\colon B_b(\Rd)\to C_b(\Rd)$ for all $t>0$,
\item[(iii)] $C_0^2(\Rd)\subseteq D(\mathcal{A})$
and $\mathcal{A}=\pL$ on $C_0^2(\Rd)$,
\item[(iv)] $(\mathcal{A},D(\mathcal{A}))$ is the closure of $(\pL,C_c^{\infty}(\Rd))$,
\item[(v)] $(P_t)_{t>0}$ is differentiable.
\end{itemize}
\end{theorem}

\noindent
In
Proposition~\ref{prop:diff_closure},
we calculate $\frac{d}{dt}P_tf=\mathcal{A}P_t f$ and we give several representations of the derivative, two of which are expressed by means of the operator $\pL$, see Theorem~\ref{thm:der_p-t},
Corollary~\ref{cor:der_p-tb} and formula \eqref{eq:der3}.

Furthermore (under stronger assumptions), we get pointwise estimates of $p(t,x,y)$.
We will need the following notation.
For $t>0$, $x\in\Rd$ we define {\it the bound function},
\begin{align}\label{def:bound_function}
\rr_t(x):=\left( [h^{-1}(1/t)]^{-d}\land \frac{tK(|x|)}{|x|^{d}} \right)
\end{align}
and auxiliary functions
\begin{align}\label{def:errx_erry}
\begin{aligned}
\errx{\gamma}{\beta}(t,x,y) &:= [h^{-1}(1/t)]^{\gamma}\left(|y-x|^{\beta}\land 1\right)t^{-1}\rr_{t}(y-x-tb_{r_t}^{x})\,, \\
\erry{\gamma}{\beta}(t,x,y) &:= [h^{-1}(1/t)]^{\gamma}\left(|y-x|^{\beta}\land 1\right)t^{-1}\rr_{t}(y-x-tb_{r_t}^{y})\,.
\end{aligned}
\end{align}
We also put
\begin{align}\label{def:eps_0_and_s_0}
\indhei{0}:=\khe\,, \qquad \qquad \indcsi{0}:=1\,.
\end{align}
\begin{theorem}\label{thm:pointwise}
Assume $\Sa$ and suppose that
\begin{align}\label{cond:pointwise}
\boost:=2 \min\left\{ \frac{\lah}{2}\land (\efcs \indhei{j})+\indcsi{j}-1 \colon\quad j=0,\ldots,\efnj\right\} > 0\,.
\end{align}
There exists $c>0$ such that for all $t\in (0,1/h(1)]$, $x,y\in\Rd$, we have
$$
p(t,x,y)\leq 
ct \left(\erry{0}{0}+\sum_{j=0}^{\efnj}\erry{\indcsi{j}-1}{\indhei{j}}\right)(t,x,y)\,.
$$
The constant $c$ can be chosen to depend only on
$d,\lmcc,\ka,\lah,\lch,h,\indnj,\boost, \varepsilon_0$, $\min\limits_{j=0,\ldots,\efnj} (\indcsi{j})$.
\end{theorem}

\noindent
The parameter $\varepsilon_0$ is defined in Section~\ref{sec:appl}, see \eqref{def:ve_0-A} and \eqref{def:ve_0-A*}.
Moreover, in Proposition~\ref{prop:remainder-bound}, we estimate the difference between $p(t,x,y)$ and a certain
zero order approximation $p_0(t,x,y)$ selected in Section~\ref{ssec:a-zoet}. In particular, the difference converges to zero as $t\to 0^+$ in $L^1(\Rd,dy)$.

\begin{corollary}\label{cor:1}
Assume $\Sa$ and suppose that
 $\indcsi{j}=1$ for all $j=1,\ldots, \efnj$. 
There exists $c>0$ such that for all $t\in (0,1/h(1)]$, $x,y\in\Rd$ we have
$$
p(t,x,y)\leq c \rr_t(y-x-t\efdrf^y_{r_t})\,.
$$
\end{corollary}
\noindent
The term $\efdrf^y_{r_t}$ in the above estimates may be replaced by
$\efdrf^x_{r_t}$, see Corollary~\ref{cor-shifts}.

Even though we do not provide lower bounds here (other than non-negativity in \eqref{ineq:p-non-neg}), the pointwise upper bound in Corollary~\ref{cor:1} may be considered as sharp. It is a topic of our future studies to obtain two-sided estimates under relaxed assumptions.

\vspace{\baselineskip}

\subsection{Examples}\label{sec:ex}

We focus on the non-local part of \eqref{def:operator}, therefore in many examples $\exdrf \equiv 0$.
We stress that we can treat the operator \eqref{def:1-stable}
without any further restrictions than those mentioned in Section~\ref{sec:int}, very much like in the first example.

\begin{example}\label{ex:1}
Let $d=1$,
$J(z)=\nu(|z|)=|z|^{-2}$. Then $h(r)= r^{-1} h(1)$. Assume that
$\exdrf \equiv 0$ and $\kappa(x,z)=a(x)k(z)$, where
$$
a(x)=1+\sqrt{x}\, \ind_{(0,1)} +\ind_{[1,\infty)} \,,\qquad
k(z)=\frac12\ind_{(-\infty,0)}+\frac32\ind_{[0,\infty)}\,.
$$
Then $\Sa$ is satisfied.
Indeed, 
see Remark~\ref{rem:close-to-1}, or note that
for every $\sigma,s\in(0,1)$ there is $c>0$ such that for all $r\in (0,1]$,
$$
\left|\efdrf_r^x \right|
= |a(x)| \log(1/r) \leq  c r^{\sigma} h(r)\,,
$$
$$
\left| \efdrf_r^x - \efdrf_r^y \right|
= |a(x)-a(y)|\log(1/r) \leq c (|x-y|^{1/2}\land 1)\, r^{s} h(r)\,,
$$
so we can choose them so that  $(\sigma/2)+ s -1>0$.
Needless to say, Theorem~\ref{thm:exist}
and \ref{thm:sem_prop}
apply with these coefficients of $\pL$, and Theorem~\ref{thm:pointwise}
yields
$$
p(t,x,y)\leq c (1+ t^{s-1}(|y-x|^{1/2}\land 1)) \left(t^{-1} \land \frac{t}{|y-x-a(y)t\log(\frac1{th(1)}) |^2}\right),
$$
for all $t\in (0,1/h(1)]$, $x,y\in\RR$.
\end{example}

Here is a simple general observation.
\begin{fact}\label{fact:kappa_prod}
Assume that
\eqref{set:h-scaling} holds with $\lah \geq 1$, \eqref{set:J} holds  and $J(z)=J(-z)$ for all $|z|<1$.
Suppose further that $\exdrf \equiv 0$ and $\kappa(x,z)=a(x)k(z)$
is such that for some $c$ we have
\begin{itemize}
\item[{\it (i)}] for all $x\in\Rd$,
$$ 0<c^{-1}\leq a(x) \leq c <\infty\,,$$
\item[{\it (ii)}] for some $\beta\in (0,1]$
 and all $x,y\in\Rd$,
$$
|a(x)-a(y)|\leq c |x-y|^{\beta}\,,$$
\item[{\it (iii)}]  for all $z\in\Rd$,
$$0<c^{-1}\leq k(z) \leq c<\infty\,,$$
\item[{\it (iv)}]
for some
 $\eta\in(0,1]$ and all  $|z|<1$,
$$
|k(z)-k(0)|\leq c |z|^{\eta}\,,
$$
\end{itemize}
Then  $\Sa$ is satisfied with 
$\indcs=\indnj=1$ and $\indhei{1}=\beta$, $\indcsi{1}=1$.
\end{fact}
\pf
Clearly, \eqref{set:k-bound} and \eqref{set:k-holder} hold.
Now, by  \cite[Lemma~5.1]{MR3996792} we have,
\begin{align*}
\left| \int_{r\leq |z|<1} z \,k(z)J(z)dz\right|
&=\frac12 \left| \int_{r\leq |z|<1} z \big(k(z)-k(-z)\big)J(z)dz\right| \leq c \int_{r\leq |z|<1} |z|^{1+\eta}\,J(z)dz\\
&= c \int_r^1 s^{d+\eta}  \, \nu(s)ds\leq c \int_r^1 s^{\eta} h(s)ds\leq c \int_r^1 s^{\eta}(r/s) h(r)ds
\leq  c\, rh(r)\,.
\end{align*}
Thus, \eqref{set:indrf-cancellation-scale}
and \eqref{set:indrf-holder-s-A} hold with the desired parameters.
\qed

We give another concrete example to 
demonstrate how cancellations can validate
 conditions \eqref{set:indrf-cancellation-scale}
and \eqref{set:indrf-holder-s-A} with
$\indcs=\indnj=1$ and $\indhei{1}=\beta$, $\indcsi{1}=1$.
Note that
Fact~\ref{fact:kappa_prod} 
 does not apply here, because neither $\lah \geq 1$  nor condition {\it (iv)} hold.

\begin{example}\label{ex:2}
Let $d=1$ and
$J(z)=\nu(|z|)=|z|^{-2}\varphi(|z|)$, where
\begin{align*}
\varphi(r)=
\begin{cases}
\dfrac{1}{\log(1/r)}, \qquad &r \leq 1/2\,, \\
\\
\dfrac1{\log(2)}, & r >1/2\,.
\end{cases}
\end{align*}  
Then $h(r)$ is comparable to $r^{-1}\varphi(r)$ for $r\in (0,1]$.
Assume that 
$\exdrf \equiv 0$ and
$\kappa(x,z)=a(x)k(z)$, where
$a(x)$
satisfies {\it (i)} and {\it (ii)} from Fact~\ref{fact:kappa_prod}, while
\begin{align*}
k(z)=1+\begin{cases}
1 &z\in (2^{-2n},2^{-2n+1}], \\
k_n \qquad &z\in  [-2^{-2n},-2^{-2n-1}), \\
0 & else,
\end{cases}
\quad \qquad n=1,2,\ldots\,,
\end{align*}
and
\begin{align*}
k_n=\frac{\log(1+1/(2n-1))}{\log(1+1/(2n))}\,.
\end{align*}
Clearly, $k(z)$ is bounded from below and above. For $r\in (0,1/2]$ we also have
\begin{align*}
\left| \int_{r\leq |z|<1} z \,k(z)J(z)dz\right|=
\left|\int_r^{1/2} \big(k(z)-k(-z)\big)\frac{1}{z \log(z)}\,  dz \right|
\leq \frac{c}{\log(1/r)}\,.
\end{align*}
Thus, $\Sa$ is satisfied with $\indcs=\indnj=1$ and $\indhei{1}=\beta$, $\indcsi{1}=1$.
\end{example}

The next example shows that 
considering a decomposition of the coefficient $\kappa(x,z)$ can provide us with improved parameters in \eqref{set:indrf-holder-s-B} or \eqref{set:indrf-holder-s-A}.

\begin{example}
Let $d=1$, $J(z)=\nu(|z|)=|z|^{-1-3/4}$. Then $h(r)= r^{-3/4} h(1)$.
We define
\begin{align*}
a_1(x)=a(x)\,, \quad k_1(z)=k(z)\,,\quad 
a_2(x)= [a(x)]^{1/3}\,, \quad 
k_2(z)=
\begin{cases}
1 & z\in(2^{-2n},2^{-2n+1}],\\
2^{1/4} \qquad & -z\in(2^{-2n-1},2^{-2n}],\\
0 & else.
\end{cases}
\end{align*}
where $a(x)$ and $k(z)$ are from Example~\ref{ex:1}, and $n$ runs through all positive integers.
Assume that  $\exdrf \equiv 0$ and let $\kappa(x,z)= a_1(x)k_1(z)+a_2(x)k_2(z)$.
Then $\Sa$ is satisfied with
$\sigma=3/4$, $\indnj=2$ and $\epsilon_1=1/2$, $s_1=3/4$, $\epsilon_2=1/6$, $s_2=1$.
In a different manner, suppose that there is $c>0$ such that for all $|x-y|\leq 1$, $r\in (0,1]$,
$$
\left| \int_{r\leq |z|<1} z(\kappa(x,z)-\kappa(y,z)) J(z)dz\right|\leq c
|x-y|^{\epsilon} r^s h(r)\,.
$$
It is not hard to see that, necessarily,
$\epsilon \leq 1/6$ and $s\leq 3/4$. In particular, $\epsilon+s-1<0$.
\end{example}

In all subsequent examples we let
\begin{align*}
\nu(r)=r^{-d-1}\varphi(r)\,,
\end{align*}
and we specify  $\varphi\colon [0,\infty)\to [0,\infty]$ in such a way that \eqref{def:nu} and \eqref{set:h-scaling} holds.
The function $J(x)$ is assumed to satisfy~\eqref{set:J}, while $\kappa(x,z)$
both conditions \eqref{set:k-bound} and \eqref{set:k-holder}.
We note that in each example below
the functions $h(r)$, $K(r)$ and $r^d \nu(r)$
are comparable for $r\in (0,1]$.
In particular, the conditions \eqref{set:indrf-cancellation-scale}
and \eqref{set:indrf-holder-s-A} read as follows: for $r\in (0,1]$,
\begin{align*}
\left|\efdrf_r^x \right| 
 &\leq c \,r^{\indcs-1} \varphi(r)\,,
\\
\left|\efdrf_r^x -\efdrf_r^y \right|
&\leq c \sum_{j=1}^{\indnj} (|x-y|^{\indhei{j}}\land 1)\, r^{\indcsi{j}-1}  \varphi(r)\,. 
\end{align*}
In Examples~\ref{ex:3}
-- \ref{ex:0_log}
 condition $\Sa$ is satisfied, since
 comments from Remark~\ref{rem:close-to-1} apply.

\begin{example}\label{ex:3}
Let $\varepsilon\in(0,1]$ and
$$
\varphi(r)=\frac{1}{[\log(2+1/r)]^{1+\varepsilon}}\,.
$$
Then 
\eqref{set:h-scaling} holds for every $\lah<1$, but not with $\lah=1$,
and \eqref{eq:intro:wusc} holds with $\uah=1$, but fails for any $\uah<1$.
Furthermore, we let
$\exdrf(x)=\int_{|z|<1} z\kappa(x,z) J(z)dz$, thus
we actually consider the operator~\eqref{eq:op-2}.
The conditions \eqref{set:exdrf-bound} and \eqref{set:exdrf-holder} are satisfied.
Here, $\lim_{r\to 0^+}\varphi(r)=0$. Note also that
$\int_{|z|<r}|z| \nu(|z|)dz$ is comparable to $[\log(2+1/r)]^{-\varepsilon}$.
\end{example}

In Examples~\ref{ex:1_log} -- \ref{ex:wide_range-2} we have 
$$\int_{|z|<1}|z|\,\nu(|z|)dz=\infty\,.$$
Of course, this was also the case in Example~\ref{ex:1}, where the unbounded term $\log(1/r)$ appeared as a result of the integral $\int_{r\leq |z|<1}|z|\nu(|z|)dz$. Example~\ref{ex:2} is a concrete realisation of Example~\ref{ex:1_log}.

\begin{example}\label{ex:1_log}
Let   $\exdrf \equiv 0$ and
$$\varphi(r)=\frac{1}{\log(2+1/r)}\,.$$ Then 
\eqref{set:h-scaling} holds for every $\lah<1$, but not with $\lah=1$,
and \eqref{eq:intro:wusc} holds with $\uah=1$, but fails for any $\uah<1$. The order of the corresponding operator is logarithmically smaller than 1. 
Here, $\lim_{r\to 0^+} \varphi(r)=0$.
Note that $\int_{r\leq |z|<1}|z|\nu(|z|)dz$ is
comparable to $\log[\log(2+1/r)]$ for small $r$.
\end{example}

\begin{example}\label{ex:0_log}
Let  $\exdrf \equiv 0$ and
$$\varphi(r)=\log(2+1/r)\,.$$
Then
\eqref{set:h-scaling}
holds with $\lah=1$, but fails for any $\lah>1$, and
\eqref{eq:intro:wusc} holds with  every $\uah>1$,
but not with $\uah=1$.
Roughly speaking, the order of the corresponding operator is logarithmically greater than 1.
Here, $\lim_{r\to 0^+} \varphi(r)=\infty$.
Note that $\int_{r\leq |z|<1}|z|\nu(|z|)dz$
comparable to \mbox{$[\log(2+1/r)]^{2}$} for small $r$.
\end{example}

In the last two examples we propose other interesting
functions $\varphi$,
which can be considered to put
our assumptions and results to a test.
The idea of constructing oscillating symbols (in our context related and comparable to the function $h$) can be traced back to \cite[Example~1.1.15]{MR1840499}.

\begin{example}\label{ex:wide_range-1}
Let 
\begin{align*}
\varphi(r):=
\begin{cases} 
c_k\,  r^{-1/4}\,,\qquad \qquad \quad &r \in \left[ ((2k+1)!)^{-1}, ((2k)!)^{-1}\right],\\
c_k \sqrt{(2k+1)!}\, r^{1/4}\,, &r \in \left[ ((2k+2)!)^{-1}, ((2k+1)!)^{-1}\right],
\end{cases}
\end{align*}
and $c_k=((2k)!!)^{-1/2}$.
We put $\varphi(r)=0$ if $r>1$.
Then 
\eqref{set:h-scaling} holds with $\lah=3/4$, but fails for any $\lah>3/4$, and \eqref{eq:intro:wusc} holds with $\uah=5/4$, but fails for any $\uah<5/4$. Roughly speaking, the order of the corresponding operator  ranges from $3/4$ to $5/4$. 
Interestingly, here $\liminf_{r \to 0^+} \varphi (r)=0$ and $\limsup_{r \to 0^+} \varphi (r)=\infty$.
\end{example}

\begin{example}\label{ex:wide_range-2}
Let 
\begin{align*}
\varphi(r):=
\begin{cases} 
c_k\,  r^{-1/4}\,,\qquad \qquad \quad &r \in \left[ ((3k+2)!)^{-1}, ((3k)!)^{-1}\right],\\
c_k \sqrt{(3k+2)!}\, r^{1/4}\,, &r \in \left[ ((3k+3)!)^{-1}, ((3k+2)!)^{-1}\right],
\end{cases}
\end{align*}
and $c_k=((3k)!!!)^{-1/2}$.
We put $\varphi(r)=0$ if $r>1$.
All the conclusions of Example~\ref{ex:wide_range-1}
are valid here, except that now $\lim_{r\to 0^+}\varphi(r)=\infty$.
\end{example}

One can extend Example~\ref{ex:wide_range-1}
and Example~\ref{ex:wide_range-2} by considering $\nu(r)=r^{-d-1} [\varphi(r)]^a$, $a\in [0,4)$, which gives an operator of order that ranges from $1-a/4$ to $1+a/4$, but we leave such discussion to the interested reader.

\vspace{\baselineskip}

\section{General approach}\label{sec:gen}

In this section we present
a general approach to the
parametrix construction motivated by \cite{KKS} and other papers on  this topic.
Our goal is to provide a transparent and flexible, yet rigorous framework. 
It will be used in Section~\ref{sec:appl} to obtain and analyse the (weak) fundamental solution for the equation $\partial_tu-\pL u=0$, where $\pL$ is given by \eqref{def:operator}.
We should note that to implement our strategy in a concrete situation, one needs to specify the zero order approximation $P_t^0$, 
which in turn determines the error term $Q_t^0$
(in Section~\ref{sec:appl} we have $Q_t^0=-(\partial_t-\pL)P_t^0$, which is typical).
A proper choice is such that allows
proving certain natural hypothesis.
In this section we focus on pointing  out these hypotheses and deriving its logical consequences.
Throughout this section we assume that:
\vspace{\baselineskip}
\begin{enumerate}
\item[($\aR$)]\label{R}
\quad $t \mapsto r_t$
is a non-negative, 
non-decreasing function such that
for all $t\in (0,1]$
$$\forall_{0<\lambda\leq 1} \,\, r_{\lambda t} \leq \sqrt{\lambda}\, r_t\,.$$
\end{enumerate}

\vspace{\baselineskip}
\noindent
We stress again that this section (Section~\ref{sec:gen}) and Section~\ref{sec:a-proofs} are general, self-contained and independent of the rest of the paper.

\subsection{Construction - functional analytic approach}\label{sec:f_a_a}

We 
use the space
$B_b(\Rd)$ equipped 
with the supremum norm $\|f\|_{\infty}$.
For $t\in (0,1]$ we consider two {\it linear} operators
\leqnomode
\begin{align}\tag{S1} \label{op:S1}
\begin{aligned}
P_t^0 \colon B_b(\Rd) \to B_b(\Rd)\,,\\
Q_t^0 \colon B_b(\Rd) \to B_b(\Rd)\,.
\end{aligned}
\end{align}
\reqnomode
According to the ideas presented in \cite[Section~5]{KKS}, we 
investigate 
objects of the form
\begin{align}\label{def:op_P}
P_t f:= P_t^0 f+\int_0^t P_{t-s}^0\, Q_s f\, ds\,,
\end{align}
where
\begin{align}\label{def:op_Q}
Q_t f := Q_t^0 f +  \sum_{n=1}^{\infty}Q_t^nf \,,
\end{align}
and
\begin{align}\label{def:op_Qn}
Q_t^nf:= \idotsint\limits_{0<s_1<\ldots<s_n<t} Q_{t-s_n}^0 \ldots Q_{s_1}^0 f \,\,ds_1 \ldots ds_n\,.
\end{align}
We will make sense of \eqref{def:op_P} provided \eqref{op:S1} holds together with the following hypotheses\\
\begin{enumerate}
\item[(H1)]\label{H1} There is $C_1>0$ such that for all $t\in (0,1]$, $f \in B_b(\Rd)$,
$$\| P_t^0 f\|_{\infty}\leq C_1 \|f\|_{\infty}\,.$$
\item[(H2)] For all $t\in (0,1]$, $f \in B_b(\Rd)$,
$$\lim_{s \to t} \|P_s^0 f - P_t^0 f\|_{\infty} =0\,.$$
\end{enumerate}
\begin{enumerate}
\item[(H3)]\label{H3} There are $C_3,\varepsilon_0 >0$ such that for all $t\in (0,1]$, $f \in B_b(\Rd)$,
$$\| Q_t^0 f\|_{\infty}\leq C_3\, t^{-1} r_t^{\varepsilon_0} \|f\|_{\infty}\,.$$ 
\item[(H4)] For all $t\in (0,1]$, $f \in B_b(\Rd)$,
$$\lim_{s \to t} \|Q_s^0 f - Q_t^0 f\|_{\infty} =0\,.$$
\end{enumerate}

\noindent
To this end we will use the notion of the Bochner integral, which in our setting boils down to the continuity and absolute integrability of the integrands in \eqref{def:op_P}--\eqref{def:op_Qn}, see, e.g., \cite{MR3617205}.

Here is a consequence of ($\aR$),
cf. \cite[Lemma~5.15]{MR3996792}.
\begin{lemma}\label{lem:time_conv}
For all $t \in (0,1]$, $\varepsilon>0$ and $k\in \NN$,
$$
\int_0^t (t-s)^{-1} r_{t-s}^{\varepsilon}\, s^{-1} r_s^{k \varepsilon} \,ds \leq  B(\varepsilon/2, (k\varepsilon)/2)\, t^{-1} r_t^{(k+1)\varepsilon}\,,
$$
and
$$
\int_0^t s^{-1} r_s^{\varepsilon}\,ds\leq B(1,\varepsilon/2)\, r_t^{\varepsilon} \,,
$$
where $B(a,b)$ is the Beta function.
\end{lemma}

The first inequality in Lemma~\ref{lem:time_conv} allows us  to propagate the bound of $Q_t^0f$ in {\rm (H3)} to the summands of the series in \eqref{def:op_Q} as follows. 
To lighten the presentation, we postpone the proofs to Section~\ref{sec:a-proofs}.

\begin{fact}\label{fact:Qn_well_Bb}
Assume \eqref{op:S1}. Suppose {\rm (H3)} and {\rm (H4)} hold. Then the integral in
\eqref{def:op_Qn}
is well defined as a Bochner integral in
$(B_b(\Rd),\|\cdot\|_{\infty})$. Further,
for all $t\in (0,1]$, $n\in\NN$ and $f\in B_b(\Rd)$,
\begin{enumerate}
\item $Q_t^n \colon B_b(\Rd)\to B_b(\Rd)$,
\item $\|Q_t^n f\|_{\infty}\leq C_3^{n+1} \prod_{k=1}^n B\!\left(\frac{\varepsilon_0}{2},\frac{k\varepsilon_0}{2}\right) t^{-1} r_t^{(n+1)\varepsilon_0}\, \|f\|_{\infty}$,
\item $\lim_{s \to t} \|Q_s^n f - Q_t^n f\|_{\infty} =0$.
\end{enumerate}
\end{fact}

\begin{fact}\label{fact:Q_well_Bb}
Assume \eqref{op:S1}. Suppose {\rm (H3)} and {\rm (H4)} hold. Then
 the series in \eqref{def:op_Q} converges absolutely in
$(B_b(\Rd),\|\cdot\|_{\infty})$.
Further,
there is $c>0$ such that for all $t\in (0,1]$ and $f\in B_b(\Rd)$,
\begin{enumerate}
\item $Q_t\colon B_b(\Rd)\to B_b(\Rd)$,
\item $\|Q_t f\|_{\infty}\leq c t^{-1} r_t^{\varepsilon_0} \|f\|_{\infty}$,
\item $\lim_{s \to t} \|Q_s f - Q_t f\|_{\infty} =0$.
\end{enumerate}
\end{fact}

Having the operators $Q_t^n$ and $Q_t$ well defined, it is convenient to note that for every $f\in B_b(\Rd)$ they
satisfy the equations
\begin{align}\label{eq:Qn}
Q_t^n f = \int_0^t Q_{t-s}^0\, Q_{s}^{n-1}f \,ds\,,\qquad n=1,2,\ldots\,,
\end{align}
and
\begin{align}\label{eq:Q}
Q_t f = Q_t^0f + \int_0^t Q_{t-s}^0\, Q_{s}f \,ds\,. 
\end{align}
Actually, we verify \eqref{eq:Qn} and \eqref{eq:Q} in the proofs of Fact~\ref{fact:Qn_well_Bb} and~\ref{fact:Q_well_Bb}.

\begin{fact}\label{fact:P_well_Bb}
Assume \eqref{op:S1}. Suppose {\rm (H1)} -- {\rm (H4)} hold.
The integral in  \eqref{def:op_P} is well defined as a Bochner integral in  $(B_b(\Rd),\|\cdot\|_{\infty})$.
Further,
there is $c>0$ such that for all $t\in (0,1]$ and $f\in B_b(\Rd)$,
\begin{enumerate}
\item $P_t\colon B_b(\Rd)\to B_b(\Rd)$,
\item $\|P_t f\|_{\infty}\leq c \|f\|_{\infty}$,
\item $\lim_{s \to t} \|P_s f - P_t f\|_{\infty} =0$.
\end{enumerate}
\end{fact}

\begin{remark}\label{rem:const}
The constant $c$ in Fact~\ref{fact:Q_well_Bb}
depends only on $C_3,\varepsilon_0,r_1$.
The constant $c$ in Fact~\ref{fact:P_well_Bb}
depends only on $C_1,C_3,\varepsilon_0,r_1$.
See \eqref{ineq:Q_norm}  and \eqref{ineq:P_norm}, respectively.
\end{remark}

In the next result we make use of a condition stronger than \eqref{op:S1},
which will allow us to assess what is known in the probability theory as {\it the strong Feller property}. Namely,
for $t\in (0,1]$ we consider the two linear operators to satisfy
\leqnomode
\begin{align}\tag{S1\!\textsuperscript{*}} \label{op:S1*}
\begin{aligned}
P_t^0 \colon B_b(\Rd) \to C_b(\Rd)\,,\\
Q_t^0 \colon B_b(\Rd) \to C_b(\Rd)\,.
\end{aligned}
\end{align}
\reqnomode
\begin{fact}\label{fact:strong-F}
Assume \eqref{op:S1*}.  
Suppose {\rm (H1)} -- {\rm (H4)} hold.
Then, for all $t\in (0,1]$, $n\in \NN$,
$$P_t,\, Q_t,\, Q_t^n\colon 
B_b(\Rd)\to C_b(\Rd)\,.$$
\end{fact}

The above delivers a convenient strategy for the construction 
and proving basic features
of the family $\{P_t: t\in (0,1]\}$. 
To obtain more properties, we
shall require that $P_t^0$ and $Q_t^0$
have yet another mapping property.
We suppose that for $t\in (0,1]$,
\leqnomode
\begin{align}\tag{S2} \label{op:S2}
\begin{aligned}
P_t^0 \colon C_0(\Rd) \to C_0(\Rd)\,,\\
Q_t^0 \colon C_0(\Rd) \to C_0(\Rd)\,.
\end{aligned}
\end{align}
\reqnomode

We prove what is known in  probability theory as {\it the Feller property}.

\begin{fact}\label{fact:PQ_well_C0}
Assume \eqref{op:S1} and \eqref{op:S2}. Suppose {\rm (H1)} -- {\rm (H4)} hold.
Then
the integrals in \eqref{def:op_P} and \eqref{def:op_Qn}
are well defined as Bochner integrals 
and the series in \eqref{def:op_Q} converges absolutely in
$(C_0(\Rd),\|\cdot\|_{\infty})$.
In particular,
\begin{enumerate}
\item $Q_t\colon C_0(\Rd)\to C_0(\Rd)$,
\item $P_t\colon C_0(\Rd)\to C_0(\Rd)$.
\end{enumerate}
\end{fact}

One of the central hypotheses here, crucial for the choice of $P_t^0$ in applications, is
the continuity at zero in the norm, as follows.
\\
\begin{enumerate}
\item[(H0)]\label{H0} For all $f\in C_0(\Rd)$,
$$
\lim_{t\to 0^+} \|P_t^0f  - f\|_{\infty}=0\,.
$$
\end{enumerate}

\vspace{\baselineskip}
\noindent
At this point it immediately translates into the same property for $P_t$.
\begin{fact}\label{fact:P-s-cont}
Assume \eqref{op:S1} and \eqref{op:S2}. Suppose {\rm (H0)} -- {\rm (H4)} hold. Then for all $f\in C_0(\Rd)$,
$$
\lim_{t\to 0^+} \|P_t f -f\|_{\infty}=0\,.
$$
\end{fact}

\vspace{\baselineskip}
The importance of \eqref{op:S1}, {\eqref{op:S1*}, \eqref{op:S2} and {\rm (H0)} -- {\rm (H4)} should be clear from the above presentation.
In particular, they give that
$\{P_t\colon t\in [0,1]\}$
is a continuous family of bounded operators 
on $C_0(\Rd)$.
Next, we want to discuss
a positivity preserving sub-Markovian semigroup corresponding in a sense to some operator $L$.
Our aim is to verify the following conjectures\\
\begin{enumerate}
\item[(CoJ1)]\label{CoJ1} For all $t\in (0,1]$ and $f\in C_0(\Rd)$,
$$
 f \geq 0 \quad\implies\quad  P_t f \geq 0\,.
$$
\item[(CoJ2)] For all $t\in (0,1]$ and $f\in C_0(\Rd)$,
$$
 f \leq 1 \quad\implies\quad  P_t f \leq 1\,.
$$
\item[(CoJ3)] 
For all $t,s>0$ (such that $t+s\leq 1$) and $f\in C_0(\Rd)$,
$$
P_t P_s f = P_{t+s}f\,.
$$
\item[(CoJ4)] For all $t\in (0,1]$ and $f\in \mathcal{D} \cap C_0(\Rd)$,
$$
P_tf(x) -f(x)= \int_0^t P_s\, L f(x) \,ds\,,
\qquad x\in\Rd
$$
(the set $\mathcal{D}$ and the operator $L$ have to be specified).
\end{enumerate}

\vspace{\baselineskip}
\noindent
In the end we shall also consider conservativeness\\
\begin{enumerate}
\item[(CoJ5)] For all $t\in (0,1]$,
$$
P_t 1 = 1\,.
$$
\end{enumerate}

\vspace{\baselineskip}
We focus on the verification of conjectures {\rm (CoJ1)} -- {\rm (CoJ4)} in the next section.
The conjecture {\rm (CoJ5)} is treated only in 
Section~\ref{ssec:ker}, where applications to the operator \eqref{def:operator} are considered.
We note that the above mentioned verification requires continuity properties similar to Fact~\ref{fact:P_well_Bb}~(3).
We collect such properties below.
For $\tau \in [0,1]$~, we let
$$
\Omega_{[\tau,1]} :=\{(t,\varepsilon)\in\RR^2\colon\,\,  t\in [\tau, 1] ,\,\,\varepsilon\geq 0 \quad\mbox{and} \quad t+\varepsilon\leq 1\}\,.
$$
\begin{lemma}\label{lem:continuity}
Assume \eqref{op:S1} and \eqref{op:S2}. Suppose {\rm (H0)} -- {\rm (H4)} hold.
Then for every $f\in C_0(\Rd)$
the following mappings are uniformly continuous into $(C_0(\Rd),\|\cdot\|_{\infty})$,
\begin{enumerate}
\item 
$
t \longmapsto P_t^0 f
$ on $[0,1]$,
\item 
$
(t,\varepsilon)\longmapsto \int_0^t P_{t-s+\varepsilon}^0\,Q_s f \,ds
$ on $\Omega_{[0,1]}$,
\item 
$
t \longmapsto Q_t^0 f
$ on $[\tau, 1]$ for every $\tau\in (0,1]$,
\item 
$
(t,\varepsilon)\longmapsto \int_0^t Q_{t-s+\varepsilon}^0\,Q_s f \,ds
$ on $\Omega_{[\tau,1]}$ for every $\tau\in (0,1]$,
\item 
$
(t,\varepsilon) \longmapsto P_{\varepsilon}^0\, Q_t f
$ on $[\tau,1]\times [0,1]$ for every $\tau\in (0,1]$.
\end{enumerate}
\end{lemma}

\begin{remark}\label{rem:1toT}
Note that 
if ($\aR$), \eqref{op:S1}, \eqref{op:S2} and {\rm (H0)} -- {\rm (H4)}
are true for $t\in (0,T]$, then
all results of this section
hold for $t\in (0,T]$
in place of $t\in (0,1]$.
In particular, in such a case, in Lemma~\ref{lem:continuity}
we have $[0,T]$ and $[\tau, T]$  instead of
$[0,1]$ and $[\tau,1]$, respectively.
Furthermore, {\rm (CoJ1)} -- {\rm (CoJ5)}
shall also be considered on $[0,T]$.
The same applies to
the results of Sections~\ref{sec:coj-gen} and~\ref{sec:realization}.
\end{remark}

\vspace{\baselineskip}

\subsection{Conjectures via approximate solution}\label{sec:coj-gen}

We assume in this subsection that ($\aR$), \eqref{op:S1}, \eqref{op:S2} and {\rm (H0)} -- {\rm (H4)} are satisfied, and we continue to discuss the general approach.
The aim is to assert that the conjectures {\rm (CoJ1)} -- {\rm (CoJ4)} hold true for  $P_tf$ 
given by \eqref{def:op_P}.
Let $\mathfrak{R}(\Rd)$ be the set of real-valued functions on $\Rd$. We ponder on the following properties
\\
\begin{enumerate}
\item[(L0)]\label{L0} $L$ is a linear operator  
defined on a linear subspace 
$D(L)$ of $\mathfrak{R}(\Rd)$ such that $Lf\in \mathfrak{R}(\Rd)$ for $f\in D(L)$.\\
\item[(L1)] $L$ satisfies the {\it positive maximum principle}, that is, \\
\begin{center}
\it
if $f\in D(L)$ is such that $f(x_0)=\sup_{x\in\Rd} f(x) \geq 0$, then  $Lf(x_0)\leq 0\,.$
\end{center}
{\,}
\item[(L2)] Non-zero constant function on $\Rd$ belongs to $D(L)$.\\
\item[(L3)]\label{L3} $\mathcal{D} \subseteq D(L)$ is such that $L f \in C_0(\Rd)$ for every $f\in \mathcal{D}$.
\end{enumerate}

\vspace{\baselineskip}
A classical method of proving
{\rm (CoJ1)} -- {\rm (CoJ4)}
is to associate $P_t f$ with an operator $L$ satisfying {\rm (L0)} -- {\rm (L3)}
by showing that $u(t,x):=P_tf(x)$ solves the equation
$\partial_t u - L u =0$, in other words, 
that it is harmonic for $\partial_t-L$.
In general, a typical problem  is that
we do not know whether $P_t f \in D(L)$ for $t>0$, which is usually a question of sufficient regularity.
In a series of papers 
\cite{MR3765882},
\cite{MR4003136},
\cite{KKS}
the problem was resolved for certain operators by 
introducing the notion of approximate harmonicity
and by specifying
the so called {\it approximate fundamental solution}. 
Namely, 
for $t,\varepsilon\geq 0$, $0<t+\varepsilon\leq 1$
and $f\in C_0(\Rd)$ we let (in the sense of Bochner)
\begin{align}\label{def:approx_sol}
P_{t,\varepsilon}f := P_{t+\varepsilon}^0 f + \int_0^t P_{t-s+\varepsilon}^0\, Q_s f\, ds\,.
\end{align}
We argue here
that such an approach is successful in general, provided
$P_{t,\varepsilon}$ is adequately regular
and  a proper relation holds between $P_t^0$ and $Q_t^0$.
This decisive relation is an equality
that we point out in the hypothesis {\rm (H5)} below.

\vspace{\baselineskip}

\begin{enumerate}
\item[(H5)]\label{H5} For all $t,\varepsilon>0$, $t+\varepsilon\leq 1$, $x\in\Rd$
and $f\in C_0(\Rd)$,
$$
(\partial_t-L_x) P_{t,\varepsilon} f(x)
= P_{\varepsilon}^0\, Q_t f(x)
-Q_{t+\varepsilon}^0 f(x) - \int_0^t Q_{t-s+\varepsilon}^0\, Q_s f(x)\,ds\,.
$$
\vspace{\baselineskip}
\item[(H6)] For all $t,\varepsilon>0$, $t+\varepsilon\leq 1$, $x\in\Rd$
and $f\in C_0(\Rd)$,
\begin{align*}
\int_0^{1-\varepsilon} |\partial_s P_{s,\varepsilon} f (x)|ds <\infty\,, \quad\quad
\int_0^{1-\varepsilon} | L_x P_{s,\varepsilon} f (x)|ds<\infty\,,
\end{align*}
\begin{align*}
L_x \int_0^t P_{s,\varepsilon} f (x)\,ds =
\int_0^t L_x P_{s,\varepsilon} f (x)\,ds \,.
\end{align*}
\end{enumerate}
\vspace{\baselineskip}
\begin{remark}\label{rem:H5_H6}
An inherent part of {\rm (H5)}
is that $t\mapsto P_{t,\varepsilon}f(x)$ is differentiable, and $P_{t,\varepsilon}f\in D(L)$.
Similarly, in {\rm (H6)}
we require that $t\mapsto P_{t,\varepsilon}f(x)$ is differentiable, both $P_{t,\varepsilon}f$ and $\int_0^t P_{s,\varepsilon} f\,ds$ belong to $D(L)$, and the integrands are measurable.
\end{remark}

\begin{remark}\label{rem:H5-ver}
In practice, to verify the hypothesis {\rm (H5)}, one initiates the construction in Section~\ref{sec:f_a_a} so that
$Q_t^0f = -(\partial_t-L)P_t^0f$,
and
 performs direct computations to evaluate $\partial_t P_{t,\varepsilon}f(x)$.
\end{remark}

Here are some general properties of
$P_{t,\varepsilon}f$ that stem from Lemma~\ref{lem:continuity}.

\begin{corollary}\label{cor:approx_sol-gen_prop}
For all $\varsigma\in(0,1)$ and $f\in C_0(\Rd)$ we have
\begin{enumerate}
\item $$\lim_{\varepsilon\to 0^+}\sup_{t\in [0,1-\varsigma]}\|P_{t,\varepsilon}f -P_tf\|_{\infty}=0\,,$$
\item 
$$
\underset{\varepsilon \in (0,\varsigma]}{\forall} \quad \lim_{|x|\to \infty} \sup_{t\in [0,1-\varsigma]} |P_{t,\varepsilon}f(x)|=0\,,
$$
\item 
$$
\lim_{\substack{(t,\varepsilon)\to (0,0)\\ t,\varepsilon \geq 0}} \|P_{t,\varepsilon}f-f\|_{\infty} = 0\,,
$$
\item
$$
\underset{\varepsilon \in (0,\varsigma]}{\forall} \quad  (t,x)\mapsto P_{t,\varepsilon}f(x) \mbox{ is an element of }\, C_0([0,1-\varsigma]\times\Rd)\,.
$$
\end{enumerate}
\end{corollary}
\pf
The reason to introduce $\varsigma$ is only to create room for $\varepsilon>0$ such that $t+\varepsilon\leq 1$.
The first statement follows from Lemma~\ref{lem:continuity},
parts (1) and (2).
The same induces that
$P_{t,\varepsilon}f$ is continuous in $t\in[0,1-\varsigma]$, hence the set $\{P_{t,\varepsilon}f\colon t\in [0,1-\varsigma]\}$ is a  compact subset of $(C_0(\Rd),\|\cdot\|_{\infty})$ and therefore the second statement holds true
as a direct consequence of sequential compactness (see \cite[Chapter~9, Theorem~16]{zbMATH05172308}).
The third one results from
$\|P_{t,\varepsilon}f-f\|_{\infty} \leq \|P_{t,\varepsilon}f-P_{t}f\|+\|P_{t}f-f\|_{\infty}$,
the first statement  and
Fact~\ref{fact:P-s-cont}.
Finally, since
$P_{t,\varepsilon}f(x)-P_{t_0,\varepsilon}f(x_0)
= P_{t,\varepsilon}f(x)-P_{t_0,\varepsilon}f(x)+P_{t_0,\varepsilon}f(x)-P_{t_0,\varepsilon}f(x_0)$, the continuity in the last statement holds by Lemma~\ref{lem:continuity},
parts (1) and (2), and because $P_{t_0,\varepsilon}f\in C_0(\Rd)$. Convergence to zero is assured by the second statement.
\qed

In what follows we prove each conjecture
{\rm (Coj1)} -- {\rm (Coj4)}
under a proper selection of conditions
{\rm (L0)} -- {\rm (L3)} and hypotheses
{\rm (H5)} -- {\rm (H6)}.
For the sake of this short discussion we use the notion of approximate harmonicity
\cite[Definition~5.4]{MR4003136}.
Roughly speaking, 
we show that
$P_tf$, $1-P_tf$ and $P_{t+s}f-P_tP_sf$ 
for $f\in C_0(\Rd)$,
as well as
$P_tf -f -\int_0^t P_sL f ds$
for $f\in \mathcal{D}\cap C_0(\Rd)$,
are approximate harmonic for $\partial_t-L$.
The corresponding approximating families are
$P_{t,\varepsilon} f$, $1-P_{t,\varepsilon}f$, $P_{t+s,\varepsilon}f-P_{t,\varepsilon}P_sf$,
and
$P_{t,\varepsilon} f -f -\int_0^t P_{s,\varepsilon} L f ds$, respectively.
This will lead to the conjectures.
The methodology of the proof is similar to that of \cite[Proposition~5.5 and Corollary~5.6]{MR4003136}. Nevertheless, since we consider an abstract operator $L$ and we do not always require (L2) in Propositions~\ref{prop:coj123} and~\ref{prop:coj4}, as well as for completeness of the paper, we provide short proofs.

\begin{fact}\label{fact:approx_sol-zero-1}
Suppose {\rm (L0)} and {\rm (H5)} hold.
For all $\tau, \varsigma>0$, $\tau+ \varsigma\leq 1$
and $f\in C_0(\Rd)$ we have
$$
\lim_{\varepsilon \to 0^+}\sup_{t\in [\tau,1-\varsigma ]} \|(\partial_t-L) P_{t,\varepsilon}f \|_{\infty} = 0\,.
$$
\end{fact}
\pf
Due to Lemma~\ref{lem:continuity} parts 
(3), (4) and (5), the expression in
the hypothesis {\rm (H5)}
converges 
in the supremum norm
uniformly with respect to $t\in  [\tau,1-\varsigma ]$ as $\varepsilon\to 0^+$ to the limit
$Q_t f - Q_t^0f + \int_0^t Q_{t-s}^0\, Q_{s}f \,ds$,
which by \eqref{eq:Q} equals zero.
This completes the proof.
\qed

We recall that our assumptions are commented in more detail in
Remark~\ref{rem:H5_H6}.

\begin{proposition}\label{prop:coj123}
Suppose {\rm (L0)}, {\rm (L1)} and {\rm (H5)} hold.
Then the conjectures {\rm (CoJ1)} and {\rm (CoJ3)} are true.
If additionally {\rm (L2)} holds, then {\rm (CoJ2)} also holds. 
\end{proposition}
\pf
Let $f\in C_0(\Rd)$ be non-negative.
Suppose there are 
$\varsigma, \theta >0$,
$t'\in [0,1-\varsigma]$ and $x'\in\Rd$ such that $e^{-\theta t'} P_{t'}f(x')<0$.
Define $u_{\varepsilon}(t,x)= e^{-\theta t} P_{t,\varepsilon}f(x)$.
By part (1) of Corollary~\ref{cor:approx_sol-gen_prop},
there are $\eta>0$ and $0<\bar{\varepsilon}_0<\varsigma$
such that
$$
u_\varepsilon(t',x') <-\eta\,,
$$
for all $0<\varepsilon<\bar{\varepsilon}_0$.
By part (4) of Corollary~\ref{cor:approx_sol-gen_prop},
$u_{\varepsilon}\in C_0([0,1-\varsigma]\times \Rd)$
and it attains  minimum at a point denoted by $(t_\varepsilon,x_\varepsilon)$. Since $f\geq 0$, by part (3) of  
Corollary~\ref{cor:approx_sol-gen_prop},
there are $\bar{t}_1>0$ and $0<\bar{\varepsilon}_1<\bar{\varepsilon}_0$
such that
$$
u_{\varepsilon}(t,x)\geq -\eta/2\,,
$$
for all $t\in [0,\bar{t}_1]$, $\varepsilon\in(0,\bar{\varepsilon}_1]$ and $x\in\Rd$. This gives $t_\varepsilon\in (\bar{t}_1,1-\varsigma]$ for all $
\varepsilon\in(0,\bar{\varepsilon}_1]$.
Therefore, by Remark~\ref{rem:H5_H6}  and {\rm (L1)}, we get
$$
\partial_t u_{\varepsilon}(t_\varepsilon,x_\varepsilon)\leq 0\,,
\qquad \quad
L_x u_\varepsilon(t_\varepsilon,x_\varepsilon)\geq 0\,,
$$
and so for all $\varepsilon\in (0,\bar{\varepsilon}_1]$ we have
$$
(\partial_t-L_x) u_\varepsilon(t_\varepsilon,x_\varepsilon) \leq 0\,.
$$
However, Fact~\ref{fact:approx_sol-zero-1} gives
\begin{align*}
(\partial_t-L_x) u_\varepsilon(t_\varepsilon,x_\varepsilon)
&= e^{-\theta t_\varepsilon} (\partial_t-L_x) P_{t_\varepsilon,\varepsilon}f(x_\varepsilon)
-\theta u_{\varepsilon}(t_\varepsilon,x_\varepsilon)\\
&\geq  e^{-\theta t_\varepsilon} (\partial_t-L_x) P_{t_\varepsilon,\varepsilon}f(x_\varepsilon)+\theta\eta 
\qquad\qquad 
\xrightarrow{\varepsilon\to 0^+} \quad \theta \eta  >0\,,
\end{align*}
which is a contradiction, and shows that $P_tf \geq 0$ for $t\in [0,1]$, see Fact~\ref{fact:P_well_Bb}.

The proof of {\rm (CoJ3)} is similar: we consider $u_\varepsilon(t,x)=e^{-\theta t}(P_{t+s,\varepsilon}f(x) - P_{t,\varepsilon} P_s f(x))$ for $t\in [0,1-s-\varsigma]$. Note that $P_sf \in C_0(\Rd)$. We also use
Lemma~\ref{lem:continuity}, parts (1) and (2),
to assert that $u_\varepsilon(t,x)\geq -\eta/2$.
Finally, we put $-f$ in place of $f$.

The proof of  {\rm (CoJ2)} is also similar: we consider
$u_\varepsilon(t,x)=e^{-\theta t}(1-P_{t,\varepsilon}f(x))$
for $t\in [0,1-\varsigma]$.
If for some point $u_\varepsilon(t',x')<-\eta$, by part (4) of Corollary~\ref{cor:approx_sol-gen_prop}, it attains its minimum.
We write $u_\varepsilon(t,x)=e^{-\theta t}(1-f(x))+e^{-\theta t}(f(x)-P_{t,\varepsilon}f(x))$ to show $u_\varepsilon(t,x)\geq -\eta/2$.
To obtain a contradiction from
$$
(\partial_t-L_x) u_\varepsilon(t_\varepsilon,x_\varepsilon)
=e^{-\theta t}(-L1) -e^{-\theta t_\varepsilon} (\partial_t-L_x) P_{t_\varepsilon,\varepsilon}f(x_\varepsilon)
-\theta u_{\varepsilon}(t_\varepsilon,x_\varepsilon)\,,
$$
we additionally use that
$L1 \leq 0$, which follows from {\rm (L1)} and {\rm (L2)}.
\qed

\begin{fact}\label{fact:approx_sol-zero-2}
Suppose {\rm (L0)}, {\rm (H5)} and {\rm (H6)} hold.
For all $\varsigma \in (0,1)$ and $f\in C_0(\Rd)$ we have
$$
\lim_{\varepsilon\to 0^+}\int_0^{1-\varsigma} \|(\partial_t-L) P_{t,\varepsilon}f \|_{\infty} \,dt=0\,.
$$
\end{fact}
\pf
By Fact~\ref{fact:approx_sol-zero-1}
it suffices to show that  
$\int_0^{\tau}\|(\partial_t-L) P_{t,\varepsilon}f \|_{\infty} \,dt$
can be made arbitrarily small uniformly in $0<\varepsilon \leq \bar{\varepsilon}_1$, by the choice of $\tau>0$ and $\bar{\varepsilon}_1>0$. To guarantee that, we apply the supremum norm to each element in the expression of
{\rm (H5)}.
We get
$\|P_\varepsilon^0 Q_t \|_{\infty}\leq C_1 c t^{-1}r_t^{\varepsilon_0}\|f\|_{\infty}$
by {\rm (H1)} and Fact~\ref{fact:Q_well_Bb};
$\|Q_{t+\varepsilon}^0f\|_{\infty}\leq C_3 (t+\varepsilon)^{-1}r_{t+\varepsilon}^{\varepsilon_0}\|f\|_{\infty}$
by {\rm (H3)};
and
$\int_0^t \| Q_{t-s+\varepsilon}^0 Q_s f\|_{\infty} ds\leq \int_0^t C_3 (t-s+\varepsilon)^{-1}r_{t-s+\varepsilon}^{\varepsilon_0} c s^{-1}r_s^{\varepsilon_0}\,ds$
by {\rm (H3)} and Fact~\ref{fact:Q_well_Bb}.
Clearly, the result follows after integrating against $dt$ over $(0,\tau)$ and using Lemma~\ref{lem:time_conv} (a change of the order of integration simplifies calculations for the third term): we just have to make $r_{\tau+\varepsilon}^{\varepsilon_0}$ small, see~($\aR$).
\qed

Again, we refer to Remark~\ref{rem:H5_H6} for important comments on our assumptions.

\begin{proposition}\label{prop:coj4}
Suppose {\rm (L0)}, {\rm (L1)}, {\rm (L3)} and {\rm (H5)}, {\rm (H6)} hold.
Then {\rm (CoJ4)} is true.
\end{proposition}
\pf
Let $f\in \mathcal{D}\cap C_0(\Rd)$.
Note that 
$f\in D(L)$ and
$L f\in C_0(\Rd)$ by {\rm (L3)}.
Suppose that there are $\varsigma, \theta>0$, $t'\in[0,1-\varsigma]$ and $x'\in \Rd$ such that 
$e^{-\theta t'}(P_{t'}f(x')-f(x')-\int_0^{t'} P_s L f(x')\,ds) <0$. 
We consider
$$
u_{\varepsilon}(t,x)= e^{-\theta t}\left(P_{t,\varepsilon}
f(x)-f(x)-\int_0^t P_{s,\varepsilon} L f(x)\,ds \right) .
$$
By part (1) of Corollary~\ref{cor:approx_sol-gen_prop},
there are $\eta>0$ and $0<\bar{\varepsilon}_0<\varsigma$
such that
$u_\varepsilon(t',x') <-\eta$
for all $0<\varepsilon<\bar{\varepsilon}_0$.
By part (4) of Corollary~\ref{cor:approx_sol-gen_prop}, $u_{\varepsilon}\in C_0([0,1-\varsigma]\times \Rd)$, hence
 $u_{\varepsilon}$ attains its minimum at a point denoted by $(t_\varepsilon,x_\varepsilon)$.
By $\int_0^t P_{s,\varepsilon} L f(x)\,ds=\int_0^t \left[ P_{s,\varepsilon} L f(x) - L f(x)\right]ds+ t L f(x)$, and
Corollary~\ref{cor:approx_sol-gen_prop} part (3),
there are $\bar{t}_1>0$ and $0<\bar{\varepsilon}_1<\bar{\varepsilon}_0$
such that 
$u_{\varepsilon}(t,x)\geq -\eta/2$
for all $t\in [0,\bar{t}_1]$, $\varepsilon\in(0,\bar{\varepsilon}_1]$ and $x\in\Rd$.
Thus $t_\varepsilon\in (\bar{t}_1,1-\varsigma]$ for all $
\varepsilon\in(0,\bar{\varepsilon}_1]$.
Therefore, by {\rm (L1)}, for all $\varepsilon\in (0,\bar{\varepsilon}_1]$ we have
$$
(\partial_t-L_x) u_\varepsilon(t_\varepsilon,x_\varepsilon) \leq 0\,.
$$
Note that we can apply 
$\partial_t$ 
and $L_x$
to $u_{\varepsilon}(t,x)$
due to our assumptions, see
Remark~\ref{rem:H5_H6},
and (for the differentiability of the integral in $t$) Corollary~\ref{cor:approx_sol-gen_prop} part (4).
However,
$$
(\partial_t-L_x) u_\varepsilon(t,x)=
e^{-\theta t}(\partial_t-L_x) (\ldots)_{(t,x)}
-\theta u_\varepsilon(t,x)
$$
where, by {\rm (H6)},
\begin{align*}
(\partial_t-L_x) (\ldots)_{(t,x)}
&:=
(\partial_t-L_x) P_{t,\varepsilon}
f(x) + L f(x) - P_{t,\varepsilon} L f(x)
+\int_0^t L_x P_{s,\varepsilon} L f(x)\,ds\\
&= (\partial_t-L_x) P_{t,\varepsilon}
f(x)
+ L f(x)
-P_{0,\varepsilon} L f(x)
-\int_0^t (\partial_s-L_x) P_{s,\varepsilon} L f(x)\,ds\,,
\end{align*}
and so by Facts~\ref{fact:approx_sol-zero-1} and~\ref{fact:approx_sol-zero-2},
$$
(\partial_t-L_x) u_\varepsilon(t_\varepsilon,x_\varepsilon)\geq
e^{-\theta t_\varepsilon}(\partial_t-L_x) (\ldots)_{(t_\varepsilon,x_\varepsilon)}
+\theta \eta \quad
\xrightarrow{\varepsilon\to 0^+} \quad \theta\eta >0\,,
$$
which is a contradiction.
Finally, it suffices to consider $-f$ in place of $f$.
\qed

\subsection{Realization via integral kernels}\label{sec:realization}

We assume in the whole subsection that ($\aR$), \eqref{op:S1} and {\rm (H1)} -- {\rm (H4)} are satisfied.
Suppose there are Borel functions
$p_0, q_0\colon (0,1] \times \Rd\times \Rd \to \RR$ such that for all $t\in (0,1]$, $x\in\Rd$ and $f\in B_b(\Rd)$,
\leqnomode
\begin{align}\tag{S3} \label{op:S3}
\begin{aligned}
P_t^0 f (x) &= \int_{\Rd} p_0(t,x,y) f(y)\,dy\,,\\
Q_t^0 f (x) &= \int_{\Rd} q_0(t,x,y) f(y)\,dy\,,
\end{aligned}
\end{align}
\reqnomode
the integrals being absolutely convergent.
The aim of this section is to make sure that 
under \eqref{op:S3}
the operator $P_t$ is an integral operator on $B_b(\Rd)$ with a Borel measurable kernel. Clearly, we have

\begin{fact} Assume \eqref{op:S3}.
The hypothesis {\rm (H1)} is tantamount to: For all $t\in (0,1]$, $x\in\Rd$,
\begin{align}
\int_{\Rd} |p_0(t,x,y)|\,dy &\leq  C_1\,. \label{ineq:H1-gen} 
\end{align}
The hypothesis {\rm (H3)} is tantamount to: For all $t\in (0,1]$, $x\in\Rd$,
\begin{align}
\int_{\Rd} |q_0(t,x,y)|\,dy \leq \, & C_3 \, t^{-1} r_t^{\varepsilon_0}\,. \label{ineq:H3-gen}
\end{align}
\end{fact}

Before treating the operator
\eqref{def:op_P}
we concentrate on \eqref{def:op_Q} and \eqref{def:op_Qn}.
We often use that
the evaluation at a point is a continuous functional on $(B_b(\Rd),\|\cdot\|_{\infty})$ and therefore commutes with Bochner integral, see \cite{MR3617205}.
We shall now assert that there exists a Borel function 
$q_n \colon (0,1] \times \Rd\times \Rd \to \RR$,
which is a kernel of
the operator $Q_t^n$, see \eqref{def:op_Qn}.

\begin{lemma}
Assume \eqref{op:S3}. 
There exist Borel sets $\mathcal{S}_n\subseteq (0,1]\times\Rd\times\Rd$
such that 
$\mathcal{S}_n$
as well as
every section $\mathcal{S}^{t,x}_n\subseteq \Rd$
are of full measure, and for
all $n\in \NN$, $t\in (0,1]$, $x\in\Rd$ and $f\in B_b(\Rd)$,
\begin{align}\label{eq:Q_n_via_q_n}
Q_t^n f(x) = \int_{\Rd} q_n(t,x,y) f(y)\,dy\,,
\end{align}
where Borel functions $q_n$ are given 
inductively for all $t\in (0,1]$, $x,y\in\Rd$
by the following absolutely convergent double integral
\begin{align}\label{def:q_n-gen}
q_n(t,x,y):=\int_0^t \int_{\Rd} \ind_{\mathcal{S}_n}(t,x,y)\, q_0(t-s,x,z)q_{n-1}(s,z,y)\,dzds\,.
\end{align}
For all $n\in\NN$, $t\in (0,1]$ and $x\in\Rd$,
\begin{align}\label{ineq:q_n-gen-b}
\int_{\Rd} |q_n(t,x,y)|\,dy\leq 
C_3^{n+1} \prod_{k=1}^n B\!\left(\frac{\varepsilon_0}{2},\frac{k\varepsilon_0}{2}\right) t^{-1} r_t^{(n+1)\varepsilon_0}\,.
\end{align}
\end{lemma}
\pf
For $n \in \NN$, $t, s_1,\ldots, s_n\in (0,1]$ and $x,y, z_1,\ldots,z_n \in \Rd$ we define Borel function
$F(t,x,s_n,z_n,\ldots,s_1,z_1,y)$ by
$$\ind_{0<s_1<\ldots<s_n<t}\,| q_0(t-s_n,x,z_n) q_0(s_n-s_{n-1},z_n,z_{n-1})\ldots q_0(s_2-s_1,z_2,z_1) q_0(s_1,z_1,y)|\,,$$
and Borel set (see \cite[Lemma~1.28]{MR4226142})
$$
\mathcal{S}_n :=\{(t,x,y)\colon \int_{(\RR)^n} \int_{(\Rd)^n} F(t,x,s_n,z_n,\ldots,s_1,z_1,y)\, dz_1\ldots dz_n \,ds_1\ldots ds_n<\infty \}\,.
$$
Let $n=1$. Clearly, 
$\mathcal{S}_1^{t,x}=\{y\colon \int_{\RR}\int_{\Rd} F(t,x,s,z,y)\,dzds<\infty \}$. 
By \eqref{ineq:H3-gen}
and Lemma~\ref{lem:time_conv} we get
\begin{align}\label{ineq:Fub-aux}
\int_{\Rd} \int_{\RR}\int_{\Rd} F(t,x,s,z,y)\,dzds\,dy < \infty\,, \tag{$\star$}
\end{align}
which confirms that $\mathcal{S}_1^{t,x}$ is of full measure. Similarly,
for every compact  $D\subset \Rd$,
$$
\int_0^1 \int_{D}\int_{\Rd} \int_{\RR}\int_{\Rd} |F(t,x,s,z,y)|\,dzds \,dy\,\,dxdt < \infty\,,
$$
hence $\mathcal{S}_1$ is of full measure.
Since 
\eqref{ineq:Fub-aux}
validates the change the order of integration, we have
\begin{align*}
Q_t^1 f(x)&=\int_{0<s<t} Q_{t-s}^0 Q_s^0 f(x)\,ds
=\int_0^t \int_{\Rd} q_0(t-s,x,z) \int_{\Rd} q_0(s,z,y)f(y)\,dy\,dzds\\
&= \int_{\Rd} \int_0^t \int_{\Rd} q_0(t-s,x,z)q_0(s,z,y)\,dzds \,f(y)\,dy = \int_{\Rd} q_1(t,x,y) f(y)\,dy \,.
\end{align*}
Finally, \eqref{ineq:q_n-gen-b} is equivalent to
part (2) of Fact~\ref{fact:Qn_well_Bb}.
The case of $n \geq 2$ goes by similar lines.
\qed

We make a similar observation about a kernel of the operator $Q_t$ defined in \eqref{def:op_Q}.

\begin{lemma}
Assume \eqref{op:S3}. 
There exists
a Borel set $\mathcal{S}_q\subseteq (0,1]\times\Rd\times\Rd$ 
such that $\mathcal{S}_q$
as well as every section $\mathcal{S}_q^{t,x}\subseteq \Rd$ are of full measure,
and for all $t\in (0,1]$, $x\in\Rd$ and $f\in B_b(\Rd)$,
\begin{align}\label{eq:Q_via_q}
Q_tf(x) = \int_{\Rd}q(t,x,y) f(y)\,dy\,,
\end{align}
where the following series converges absolutely for all $t\in (0,1]$, $x,y\in\Rd$,
\begin{align}\label{def:q-gen}
q(t,x,y):=\sum_{n=0}^\infty  \ind_{\mathcal{S}_q}(t,x,y) \,q_n(t,x,y)\,.
\end{align}
Further, there exists $c>0$ such that  for all $t\in (0,1]$ and $x\in\Rd$, 
\begin{align}\label{ineq:q-gen-b}
\int_{\Rd} |q(t,x,y)|\,dy \leq 
c \, t^{-1} r_t^{\varepsilon_0}\,.
\end{align}
\end{lemma}
\pf
Let $F(t,x,y):=\sum_{n=0}^\infty |q_n(t,x,y)|$, where $q_n$ are Borel functions defined in \eqref{def:q_n-gen}.
We introduce a Borel set $\mathcal{S}_q:=\{(t,x,y)\colon  F(t,x,y) <\infty \}$.
Using \eqref{ineq:q_n-gen-b}
and Lemma~\ref{lem:time_conv}, 
see 
proof of Fact~\ref{fact:Q_well_Bb},
for every compact  $D\subset \Rd$,
$$
\int_{\Rd} F(t,x,y)\,dy<\infty\,,\qquad\quad
\int_0^1 \int_{D}\int_{\Rd} F(t,x,y)\,dy\,dxdt<\infty\,,
$$
hence $\mathcal{S}_q^{t,x}$ and $\mathcal{S}_q$ are of full measure.
By
Fact~\ref{fact:Q_well_Bb}
we can calculate the following limit pointwise
\begin{align*}
Q_tf(x)=\lim_{N\to \infty} \sum_{n=0}^N Q_t^n f(x)\,.
\end{align*}
Therefore, by
\eqref{eq:Q_n_via_q_n} and since $\mathcal{S}_q^{t,x}$ is of full measure, and by the dominated convergence theorem,
\begin{align*}
Q_tf(x)
&=\lim_{N\to\infty} \int_{\Rd} \sum_{n=0}^N q_n(t,x,y)f(y)\,dy
=\lim_{N\to\infty} \int_{\Rd} \ind_{\mathcal{S}_q} (t,x,y)\sum_{n=0}^N q_n(t,x,y)f(y)\,dy\\
&= \int_{\Rd} \lim_{N\to\infty} \ind_{\mathcal{S}_q}(t,x,y)\sum_{n=0}^N q_n(t,x,y)f(y)\,dy
=\int_{\Rd} q(t,x,y) f(y)\,dy\,.
\end{align*}
Finally,
\eqref{ineq:q-gen-b}
is equivalent to Fact~\ref{fact:Q_well_Bb} part (2).
\qed

We deal with the operator $P_t$ defined in \eqref{def:op_P}.

\begin{proposition}\label{prop:p-gen}
Assume \eqref{op:S3}.
There exists a Borel set $\mathcal{S}_p\subseteq (0,1]\times\Rd\times\Rd$ 
such that $\mathcal{S}_p$
as well as every section $\mathcal{S}_p^{t,x}\subseteq \Rd$ are of full measure, and
for  all $t\in (0,1]$, $x\in\Rd$ and $f\in B_b(\Rd)$,
\begin{align}\label{eq:P_via_p}
P_tf(x) = \int_{\Rd} p(t,x,y) f(y)\,dy\,,
\end{align}
where a Borel function $p$ is given for all $t\in (0,1]$, $x,y\in\Rd$ by the
following formula with absolutely convergent
double integral
\begin{align}\label{def:p-gen}
p(t,x,y):=p_0(t,x,y)+ \int_0^t \int_{\Rd}
\ind_{\mathcal{S}_p}(t,x,y)\, p_0(t-s,x,z)q(s,z,y)\,dzds\,.
\end{align}
Further, there exists $c>0$ such that for all $t\in (0,1]$ and $x\in\Rd$, 
\begin{align}\label{eq:p-gen}
\int_{\Rd} | p(t,x,y)| \,dy \leq c\,.
\end{align}

\end{proposition}
\pf
Let
$F(t,x,s,z,y):=\ind_{0<s<t} \,|p_0(t-s,x,z)q(s,z,y)|$, where $q$ is a Borel function given by \eqref{def:q-gen}.
We define a Borel set
$\mathcal{S}_p:=\{(t,x,y)\colon \int_{\RR}\int_{\Rd} F(t,x,s,z,y)\,dzds<\infty \}$.
Integrating over $\Rd$ against $dy$ 
or  over $(0,1]\times D \times \Rd$
against $dy\,dxdt$, for arbitrary compact $D\subset \Rd$,
and using
\eqref{ineq:q-gen-b},
\eqref{ineq:H1-gen} and
Lemma~\ref{lem:time_conv},
we conclude that
$\mathcal{S}_p^{t,x}$ and $\mathcal{S}_p$ are of full measure.
It also validates the change the order of integration
\begin{align*}
\int_0^t P_{t-s}^0\, Q_s f (x)\, ds
=\int_0^t \int_{\Rd} p_0(t-s,x,z)\int_{\Rd} q(s,z,y)f(y)\,dy\,dzds\\
=\int_{\Rd} \int_0^t\int_{\Rd} \ind_{\mathcal{S}_p}(t,x,y) p_0(t-s,x,z) q(s,z,y)\,dzds\,  f(y)\,dy\,.
\end{align*}
Now, \eqref{eq:P_via_p} follows from \eqref{def:op_P}, while 
\eqref{eq:p-gen} is equivalent to
Fact~\ref{fact:P_well_Bb} part (2).
\qed

\section{Proofs of the main results}\label{sec:appl}

In this section we finally treat the operator \eqref{def:operator}.
We assume that either $\Sb$ or $\Sa$ holds, see
Section~\ref{sec:set}.
In what follows,  $T>0$ is {arbitrary and fixed}.
Throughout this section, for $t>0$,
\begin{align}\label{def:r_t}
r_t:=h^{-1}(1/t)\,.
\end{align}
Note that $r_t$ in \eqref{def:r_t} satisfies ($\aR$) for all $t\in (0,T]$, see \cite[Lemma~5.1]{MR3996792} and
Remark~\ref{rem:1toT}.
We will also need a fixed number $\varepsilon_0$ such that
under $\Sb$,
\begin{align}\label{def:ve_0-A}
\boxed{
0<\varepsilon_0< \min\{\lah\land (\efcs \indhei{j})+ \indcsi{j} -1\colon\,\, j=0,\ldots,\efnj\}
\qquad \mbox{and}\qquad
0<\varepsilon_0\leq \lah+\efcs-1,}
\end{align}
while under $\Sa$,
\begin{align}\label{def:ve_0-A*}
\boxed{
0<\varepsilon_0< \min\{\lah\land (\efcs \indhei{j})+ \indcsi{j} -1\colon\,\, j=0,\ldots,\efnj\}.}
\end{align}
Its existence results from the assumptions, and will be crucial in what follows.
We also put
\begin{align}\label{def:indTr}
\indTr:=1/h(1)\,.
\end{align}

In the subsequent sections we prove that
the general approach from 
Section~\ref{sec:gen} applies and so all the conclusions of results in Section~\ref{sec:gen} hold.
While using the auxiliary results of Section~\ref{sec:sac}
or~\ref{sec:a-frozen},
the reader may conveniently 
notice that
all the assumptions therein are implied 
by $\Sb$ - the only exclusion here being  Lemma~\ref{lem:q_0-aux} - and separately by $\Sa$.


\subsection{Zero order approximation and error term}\label{ssec:a-zoet}
In this subsection we specify
$P_t^0$
and $Q_t^0$, both of which shall be integral operators.
We start by ,,freezing''  coefficients of the operator~\eqref{def:operator}.
Namely, for $w\in\Rd$ we define
the operator 
\begin{align}\label{def:pL_fr}
\pL^{\mathfrak{K}_w} f(x):=\exdrf(w)\cdot \nabla f(x)+ \int_{\Rd}\Big(f(x+z)-f(x)- \ind_{|z|<1} \left<z,\nabla f(x)\right>\! \Big) \,\kappa(w,z)J(z)dz\,,
\end{align}
which makes sense whenever $f\in C_0^2(\Rd)$. 
Already when considered on $C_c^{\infty}(\Rd)$, $\pL^{\mathfrak{K}_w}$ uniquely determines a L{\'e}vy process with a transition density $p^{\mathfrak{K}_w}(t,x)$,
see \cite[Section~6]{MR3996792}.
Recall  that the density is given via the Fourier inversion formula \cite[(97)]{MR3996792}. In 
Section~\ref{sec:a-frozen}
we provide more properties of that function, which are necessary for the discussion below.
For $t>0$, $x,y\in\Rd$
we write
\begin{align}\label{def:p_fr}
p^{\mathfrak{K}_w}(t,x,y):=p^{\mathfrak{K}_w}(t,y-x)\,,
\end{align}
and we let
\begin{align}\label{def:p0q0}
\begin{aligned}
p_0(t,x,y)&:=p^{\mathfrak{K}_y}(t,x,y)\,,\\
q_0(t,x,y)&:=-(\partial_t-\pL_x)\,p_0(t,x,y)\,.
\end{aligned}
\end{align}
We record that these functions are jointly continuous in all three variables. This is a consequence of Lemma~\ref{lem:a-continuity} and \eqref{eq:par_t_is_L}.
The specific choice of $q_0$ is typical for the parametrix method \cite[Section~1.7]{MR3701414}, \cite[Section~5]{KKS}.
In our considerations we  substantially use that for $r>0$,
\begin{align}
\pL_x^{\mathfrak{K}_{v}} p^{\mathfrak{K}_{w}}(t,x,y)
= 
\efdrf_{r}^{v}\cdot \nabla_x p^{\mathfrak{K}_{w}}(t,x,y)
+
\int_{\Rd}\delta_{r}^{\mathfrak{K}_{w}}(t,x,y;z)\kappa(v,z)J(z)dz\,, \label{eq:pL}
\end{align}
where
\begin{align}\label{def:delta}
\delta_{r}^{\mathfrak{K}_w} (t,x,y;z)&:=p^{\mathfrak{K}_w}(t,x+z,y)-p^{\mathfrak{K}_w}(t,x,y)-\ind_{|z|<r}\left< z,\nabla_x p^{\mathfrak{K}_w}(t,x,y)\right>.
\end{align}
Thus
\begin{align}\label{eq:delta-alt}
\begin{aligned}
q_0(t,x,y)&=-(\partial_t-\pL_x^{\mathfrak{K}_x}) p^{\mathfrak{K}_y}(t,x,y)
=(\pL_x^{\mathfrak{K}_x}-\pL_x^{\mathfrak{K}_y}) p^{\mathfrak{K}_y}(t,x,y) \\
&= 
(\efdrf_{r}^{x}-\efdrf_{r}^{y}) \cdot \nabla_x p^{\mathfrak{K}_y}(t,x,y)
+
\int_{\Rd} \delta_{r}^{\mathfrak{K}_y} (t,x,y;z) (\kappa(x,z)-\kappa(y,z))J(z)dz.
\end{aligned}
\end{align}
Now, for $t>0$, $x\in\Rd$ and $f\in B_b(\Rd)$ we set
\begin{align}\label{def:P0Q0}
\begin{aligned}
P_t^0 f (x) &:= \int_{\Rd} p_0(t,x,y) f(y)\,dy\,,\\
Q_t^0 f (x) &:= \int_{\Rd} q_0(t,x,y) f(y)\,dy\,.
\end{aligned}
\end{align}
These operators are well defined, as observed right at the beginning in  Lemma~\ref{lem:well_def}.

\subsection{Hypothesis - construction}\label{ssec:a-hc}

Our aim in this section is to 
show that
with the above choice of $P_t^0$ and $Q_t^0$, the conditions
\eqref{op:S1}, \eqref{op:S1*}, \eqref{op:S2}, and {\rm (H0)} -- {\rm (H4)}
are satisfied for $t\in (0,T]$.
\begin{lemma}\label{lem:well_def}
There exist  constants $C_1, C_3>0$ such that for all $t\in (0,T]$, $x\in\Rd$,
\begin{align}
\int_{\Rd} p_0(t,x,y)\,dy &\leq  C_1\,, \label{ineq:H1} \\
\int_{\Rd} |q_0(t,x,y)|\,dy \leq \, & C_3 \, t^{-1} r_t^{\varepsilon_0}\,. \label{ineq:H3}
\end{align}
\end{lemma}
\pf
The inequality \eqref{ineq:H1}
follows from Proposition~\ref{prop:gen_est}, Corollary~\ref{cor-shifts} and Lemma~\ref{lem:conv}(a).
We focus on \eqref{ineq:H3} and $t\in (0,\indTr]$. 
Under $\Sb$, Lemmas~\ref{lem:q_0-aux*1} 
and~\ref{lem:q_0-aux*2} guarantee that
\begin{align*}
\int_{\Rd} |q_0(t,x,y)|\,dy 
&=\int_{|x-y|<R} |q_0(t,x,y)|\,dy
+ \int_{|x-y|\geq R} |q_0(t,x,y)|\,dy\\
&\leq 
c 
\sum_{j=0}^{\efnj} \int_{\Rd} \erry{\indcsi{j}-1}{\indhei{j}}(t,x,y)\,dy
+ c t^{-1}r_t^{\lah+\efcs-1}\,,
\end{align*}
while under $\Sa$, by
Lemma~\ref{lem:q_0-aux}, we have
\begin{align*}
\int_{\Rd} |q_0(t,x,y)|\,dy\leq c 
\sum_{j=0}^{\efnj} \int_{\Rd} \erry{\indcsi{j}-1}{\indhei{j}}(t,x,y)\,dy\,.
\end{align*}
Now, due to the choice of $\varepsilon_0$ in \eqref{def:ve_0-A} and \eqref{def:ve_0-A*}, for each $j=0,\ldots,\efnj$ there exists
$0<\ell_j < \lah \land (\efcs \indhei{j})$
such that $\varepsilon_0=\ell_j+\indcsi{j}-1$.
Applying Corollary~\ref{cor-shifts} and Lemma~\ref{lem:conv}(a) with
$\beta_0=\ell_j$ we get
$$
\int_{\Rd} \erry{\indcsi{j}-1}{\indhei{j}}(t,x,y)\,dy
\leq \frac{c_1}{\lah-\ell_j}t^{-1} r_t^{\ell_j+\indcsi{j}-1}= \frac{c_1}{\lah-\ell_j} t^{-1} r_t^{\varepsilon_0}\,.
$$
Obviously, $r_t^{\lah+\efcs-1}\leq r_t^{\varepsilon_0}$ for $t\in (0,\indTr]$.
This ends the proof of \eqref{ineq:H3} if $0<T\leq \indTr$.
Now, let $T>\indTr$.
Applying
\eqref{eq:delta-alt} with $r=1$,
\begin{align*}
| q_0(t,x,y)|
=
&\Big|\,(\efdrf_{r_1}^x-\efdrf_{r_1}^y) \cdot \nabla_x  p^{\mathfrak{K}_y}(t,x,y)
+
\int_{\Rd} \delta_{1}^{\mathfrak{K}_y} (t,x,y;z)(\kappa(x,z)-\kappa(y,z))J(z)dz \,\Big|\,.
\end{align*}
Applying
\eqref{set:indrf-cancellation-scale} and
\eqref{ineq:L1_uni_time-1}
to the first term of the sum,
and
\eqref{set:J},
\eqref{set:k-bound},
\eqref{ineq:L1_uni_time-2},
\eqref{ineq:L1_uni_time-3}
to the second one, we get
that for $t\in [\tau,T]$,
\begin{align}\label{ineq:q_0-away}
|q_0(t,x,y)|\leq c \left(\rr_{\tau}(y-x)
+\int_{|z|\geq 1} \rr_{\tau}(y-x-z)\nu(|z|)dz\right).
\end{align}
We take $\tau = \indTr$ in \eqref{ineq:q_0-away}, integrate the inequality against $dy$ and  use \cite[Lemma~5.6]{MR3996792}, which yields $\int_{\Rd}|q_0(t,x,y)|\,dy\leq c$. It remains to notice that $1\leq (T r_\indTr^{-\varepsilon_0}) t^{-1}r_t^{\varepsilon_0}$.
\qed

Clearly, $P_t^0 f$ and $Q_t^0 f$ are Borel functions, see \cite{MR4226142}. We have the following conclusion.
\begin{corollary}\label{cor:S1H1H3}
The condition \eqref{op:S1} as well as hypotheses {\rm (H1)} and {\rm (H3)} are satisfied.
\end{corollary}

\vspace{0.25\baselineskip}

\begin{lemma}
For all $t>0$ we have
\begin{align*}
\lim_{s\to t}\sup_{x\in\Rd} \int_{\Rd} | p_0(t,x,y) - p_0(s,x,y)|\,dy&=0\,,\\
\lim_{s\to t}\sup_{x\in\Rd} \int_{\Rd} | q_0(t,x,y) - q_0(s,x,y)|\,dy&=0\,.
\end{align*}
\end{lemma}
\pf
Using \eqref{ineq:L1_uni_time-1} we get for $t,s\geq \tau$,
\begin{align*}
|p_0(t,x,y) - p_0(s,x,y)|
= | p^{\mathfrak{K}_y}(t,x,y)-p^{\mathfrak{K}_y}(s,x,y)|
\leq c |t-s|\, \rr_{\tau}(y-x)\,.
\end{align*}
It suffices to integrate the inequality, see \cite[Lemma~5.6]{MR3996792}.
Using
\eqref{eq:delta-alt} with $r=1$,
\begin{align*}
| q_0(t,x,y) - q_0(s,x,y)|
=
&\Big|\,(\efdrf_{r_1}^x-\efdrf_{r_1}^y) \cdot \Big[ \nabla_x  p^{\mathfrak{K}_y}(t,x,y)
-\nabla_x p^{\mathfrak{K}_y}(s,x,y)\Big]\\
&+
\int_{\Rd} \Big[ \delta_{1}^{\mathfrak{K}_y} (t,x,y;z)-\delta_{1}^{\mathfrak{K}_y} (s,x,y;z)\Big] (\kappa(x,z)-\kappa(y,z))J(z)dz \,\Big|\,.
\end{align*}
Further, using
\eqref{set:indrf-cancellation-scale} and
\eqref{ineq:L1_uni_time-1}
for the first term of the sum,
and
\eqref{set:J},
\eqref{set:k-bound},
\eqref{ineq:L1_uni_time-2},
\eqref{ineq:L1_uni_time-3}
for the second one, we have
for $t,s\geq \tau$,
$$
| q_0(t,x,y) - q_0(s,x,y)|
\leq c |t-s| \, \left(\rr_{\tau}(y-x)
+\int_{|z|\geq 1} \rr_{\tau}(y-x-z)\nu(|z|)dz\right).
$$
Integrating against $dy$ and using \cite[Lemma~5.6]{MR3996792} gives the result.
\qed

\begin{corollary}
The hypotheses {\rm (H2)} and {\rm (H4)} are satisfied.
\end{corollary}

\vspace{0.25\baselineskip}

\begin{lemma}
The condition \eqref{op:S2} is satisfied.
\end{lemma}
\pf
Note that
$P_t^0 f(x) = \int_{\Rd} p_0(t,x,y+x)f(y+x)dy$ and by \eqref{ineq:L1_uni_time-1},
$$
|p_0(t,x,y+x)f(y+x)|\leq c \rr_t(y) \|f\|_{\infty}\,.
$$
The estimate validates the use of the dominated convergence theorem
\cite[Lemma~5.6]{MR3996792}, hence the continuity of the integral follows from that of the integrand, see Lemma~\ref{lem:a-continuity}. Similarly, convergence to zero at infinity follows from that of $f$ and the boundedness of $p_0(t,x,y)$.
Likewise, $Q_t^0 f(x) = \int_{\Rd} q_0(t,x,y+x)f(y+x)dy$, and by \eqref{ineq:q_0-away},
$$|q_0(t,x,y+x)f(y+x)|\leq c\left( \rr_t(y)+\int_{|z|\geq 1}
\rr_t(y-z) \nu(|z|)dz
\right)\|f\|_{\infty}\,.$$
This ends the proof, see
Lemma~\ref{lem:a-continuity} and \eqref{eq:par_t_is_L}.
\qed

We show the strong continuity at zero.

\begin{lemma}
For every $f\in C_0(\Rd)$,
\begin{align*}
\lim_{t\to 0^+}\sup_{x\in\Rd} \left|\int_{\Rd} p_0(t,x,y) f(y)\,dy - f(x)\right| =0\,.
\end{align*}
\end{lemma}
\pf
We have
\begin{align*}
\int_{\Rd} p_0(t,x,y) f(y)\,dy - f(x)
= \int_{\Rd} p^{\mathfrak{K}_y}(t,x,y)\big[f(y)-f(x)\big]\,dy
+ \left[ \int_{\Rd}p^{\mathfrak{K}_y}(t,x,y)\,dy-1\right] f(x) \,.
\end{align*}
Since $f$ is bounded, due to
Proposition~\ref{prop:strong_at_zero}
the second term converges to zero in the desired fashion.
Since $f\in C_0(\Rd)$, there exists $\delta>0$ such that $|f(y)-f(x)|\leq \varepsilon$ whenever $|y-x|\leq \delta$.
Thus, by \eqref{ineq:H1}, we get
\begin{align*}
\int_{\Rd} p^{\mathfrak{K}_y}(t,x,y)| f(y)-f(x)|\,dy
\leq 
\varepsilon \,C_1
+2 \|f\|_{\infty}\int_{|x-y|>\delta} p^{\mathfrak{K}_y}(t,x,y)\,dy\,.
\end{align*}
Again, the second term here converges to zero by Lemma~\ref{lem:strong_at_zero-aux}.
This ends the proof, because the choice of $\varepsilon$ was arbitrary.
\qed

\begin{corollary}\label{cor:H0}
The hypothesis {\rm (H0)} is satisfied.
\end{corollary}

\begin{lemma}\label{lem:strong-F}
The condition \eqref{op:S1*} is satisfied.
\end{lemma}
\pf
Since \eqref{op:S1} holds by Corollary~\ref{cor:S1H1H3},
it suffices to show the continuity of $P_t^0f$ and $Q_t^0f$. Fix $t\in (0,T]$ and $R>0$. By \cite[Corollary~5.10]{MR3996792} for all $|x|\leq R$ and $y\in\Rd$,
\begin{align}\label{ineq:remove_x}
\rr_t(y-x)\leq c \rr_t(y),
\end{align}
for some $c>0$. Thus, by \eqref{ineq:L1_uni_time-1} we have $p_0(t,x,y) \leq c \rr_t(y)$,
while by \eqref{ineq:q_0-away} we get
$|q_0(t,x,y)|\leq c\rr_t(y)+ c\int_{|z|\geq 1} \rr_{t}(y-z)\nu(|z|)dz$.
Recall that $p_0$ and $q_0$ are continuous.
Now, the desired continuity follows by 
\cite[Corollary~5.6]{MR3996792} and the dominated convergence theorem.
\qed

In summary, we verified \eqref{op:S1}, \eqref{op:S1*},  \eqref{op:S2} and {\rm (H0)} -- {\rm (H4)} for $t\in (0,T]$, see also Remark~\ref{rem:1toT}. They are the foundations of the functional analytic approach, therefore all the conclusions of results in Section~\ref{sec:f_a_a} are at hand. 
In particular, we constructed
operators $P_t$,\label{P_t} $Q_t$, $Q_t^n$ as in
\eqref{def:op_P},
\eqref{def:op_Q},
\eqref{def:op_Qn}, respectively.
Furthermore,
since $P_t^0$ and $Q_t^0$ were defined in \eqref{def:P0Q0} as integral operators, i.e., the condition
\eqref{op:S3} holds,
all the conclusions of
Section~\ref{sec:realization}
are valid for $t\in (0,T]$.
Among other things, the operators $P_t$, $Q_t$, $Q_t^n$ are integral operators with integral kernels $p(t,x,y)$,  $q(t,x,y)$, $q_n(t,x,y)$,
see \eqref{eq:P_via_p}, \eqref{eq:Q_via_q}, \eqref{eq:Q_n_via_q_n}, respectively.
We also refer to Remark~\ref{rem:t-global} below.
We continue more detailed analysis in the next subsection.

\vspace{\baselineskip}

\subsection{Hypothesis - conjectures}\label{ssec:b-hc}

The aim in this section is to
specify the operator $(L,D(L))$ and the set $\mathcal{D}$ such that
{\rm (L0)} -- {\rm (L3)} are satisfied,
and to ensure that {\rm (H5)} and {\rm (H6)} hold true. 
Due to the results of Section~\ref{sec:coj-gen},
this will lead to {\rm (CoJ1)} -- {\rm (CoJ4)} for the objects we constructed.
To be more precise,
we consider $L f:=\pL f$ given by \eqref{def:operator} for $f\in D(L)=D(\pL)$, where
\begin{align}\label{def:DpL}
\begin{aligned}
D(\pL):=\big\{f\colon\,  \nabla f(x) \, &\mbox{\it exists, is finite and the integral in \eqref{def:operator}}\\
&\mbox{\it converges absolutety for every } x\in\Rd \big\}\,,
\end{aligned}
\end{align}
and we choose 
\begin{align}\label{def:pD}
\mathcal{D}:=C_0^2(\Rd)\,.
\end{align}

\vspace{\baselineskip}

\noindent
Note that {\rm (L0)}, {\rm (L1)},
{\rm (L2)}
and {\rm (L3)} are indeed satisfied.
We need a few direct computations before we 
can conclude that {\rm (H5)} and {\rm (H6)} are true on $(0,T]$, see Remark~\ref{rem:H5_H6} and~\ref{rem:H5-ver}.
From the technical point of view, the main ingredient of all the proofs is
Corollary~\ref{cor:L1_uni_time},
which provides various estimates of the integral kernel $p_0(t,x,y)$.

We point out that in our setting the choice of $D(\pL)$ and $\mathcal{D}$ is not the only one possible. Another good one, which appears in the literature, is  $C_0^2(\Rd)$ as $D(\pL)$, and $C_c^{\infty}(\Rd)$ as $\mathcal{D}$.
In general, the choice has to guarantee that $P_{t,\varepsilon}f\in D(\pL)$ for every $f\in C_0(\Rd)$, see Remark~\ref{rem:H5_H6}.

\begin{lemma}\label{lem:point-1}
For all $t\in (0,T]$, $x\in\Rd$
and $f\in C_0(\Rd)$,
\begin{align*}
\partial_t P_t^0 f (x) = \int_{\Rd}  \partial_t \,p_0(t,x,y) f(y)\,dy\,,
\qquad
\quad
\pL_x P_t^0 f (x) = \int_{\Rd} \pL_x   \,p_0(t,x,y) f(y)\,dy\,.
\end{align*}
In particular,
$$
Q_t^0 f (x) = -(\partial_t-\pL_x) P_t^0f(x)\,.
$$
\end{lemma}
\pf
The first equality follows from the
dominated convergence theorem
justified by
\eqref{ineq:L1_uni_time-1}
and \cite[Lemma~5.6]{MR3996792}.
The rest of the result will follow if we can differentiate under the integral in the first spatial variable (by dominated convergence theorem)
and change the order of integration (by Fubini's theorem),
which validate the steps below, where
$
F(x):=P_t^0f(x)
$,
\begin{align*}
\pL_x F(x)
&=\exdrf(x)\cdot \nabla_x 
F(x)+
\int_{\Rd}\Big(F(x+z)-F (x)-\ind_{|z|<1} \left<z,\nabla_x F(x)\right>\!\Big)\kappa(x,z)J(z)dz\\
&=\int_{\Rd}\exdrf(x)\cdot \nabla_x p^{\mathfrak{K}_y}(t,x,y) f(y)\,dy
+\int_{\Rd}\int_{\Rd}
\delta_1^{\mathfrak{K}_y}(t,x,y;z)f(y)\,dy\,
\kappa(x,z)J(z)dz\\
&=\int_{\Rd}\Bigg( \exdrf(x)\cdot \nabla_x p^{\mathfrak{K}_y}(t,x,y)
+\int_{\Rd}
\delta_1^{\mathfrak{K}_y}(t,x,y;z)\,
\kappa(x,z)J(z)dz\Bigg)
f(y)\,dy\\
&= \int_{\Rd} \pL_x p_0(t,x,y) f(y)\,dy\,,
\end{align*}
 see \eqref{def:operator} and \eqref{eq:pL}.
The differentiability in the spatial variable
is justified by
\cite[Lemma~5.6]{MR3996792}
and the fact that given $\tau>0$, for all $t\in [\tau,T]$, $x,y\in\Rd$ and $0<|\varkappa|<1$,
we have
\begin{align}\label{lem:point-1-aux}
\begin{aligned}
|\frac1{\varkappa}\left( p^{\mathfrak{K}_y}(t,x+\varkappa e_i,y)-p^{\mathfrak{K}_y}(t,x,y)\right)|
&\leq  \int_0^1 |\partial_{x_i}p^{\mathfrak{K}_y}(t,x+ \theta \varkappa  e_i,y) |\,d\theta \\
&\leq c \rr_{\tau}(y-x-\theta \varkappa e_i)
\leq c \rr_{\tau}(y-x)\,,
\end{aligned}
\end{align}
where the inequalities follow from
\eqref{ineq:L1_uni_time-1}
and \cite[Corollary~5.10]{MR3996792}, since $|\theta\varkappa e_i|<1= (r_{\tau}^{-1}) r_{\tau}$.
Fubini's theorem can then be used due to
\eqref{ineq:L1_uni_time-2}
and \eqref{ineq:L1_uni_time-3}
accompanied by
\eqref{set:J},
\eqref{set:k-bound}
and
\cite[Lemma~5.6]{MR3996792}.
In particular, we get $P_t^0f \in D(\pL)$.
\qed

We now obtain
results 
similar to Lemma~\ref{lem:point-1}, but
for the integral part of \eqref{def:approx_sol}.

\begin{lemma}\label{lem:point-2}
For all $t,\varepsilon>0$, $t+\varepsilon\leq T$,
$x\in\Rd$ and $f\in C_0(\Rd)$,
\begin{align*}
\partial_t \int_0^t P_{t-s+\varepsilon}^0\, Q_s f (x)\, ds
= P_{\varepsilon}^0\, Q_t f(x)
+ \int_0^t \partial_t P_{t-s+\varepsilon}^0\, Q_s f (x)\, ds\,.
\end{align*}
\end{lemma}
\pf
All the integrals are absolutely convergent Lebesgue integrals.
We divide the following expression by $t-t_1$,
\begin{align*}
&\int_0^t P_{t-s+\varepsilon}^0\, Q_s f (x)\, ds
- \int_0^{t_1} P_{t_1-s+\varepsilon}^0\, Q_s f (x)\, ds\\
&\quad=\int_{t_1}^t P_{t-s+\varepsilon}^0\, Q_s f (x)\, ds
+\int_0^{t_1} \left[P_{t-s+\varepsilon}^0\, Q_s f (x)-P_{t_1-s+\varepsilon}^0\, Q_s f (x)\right] ds\,.
\end{align*}
Then due to the continuity in
Lemma~\ref{lem:continuity} part (5)
the first term on the right hand side converges to $P_{\varepsilon}^0\, Q_{t_1} f(x)$ as $t\to t_1$.
We rewrite the second term as
$$
\int_0^{t_1} \int_{\Rd}
\big[p_0(t-s+\varepsilon,x,y) - p_0(t_1-s+\varepsilon,x,y)\big]\, Q_sf(y)\, dy ds\,.
$$
After dividing by $t-t_1$, then 
applying bounds in \eqref{ineq:L1_uni_time-1}
and part (2) of Fact~\ref{fact:Q_well_Bb},
and using integrability 
from \cite[Lemma~5.6]{MR3996792}
and Lemma~\ref{lem:time_conv},
we can
pass to the limit $t\to t_1$ under the integral to obtain
$
\int_0^{t_1} \int_{\Rd} \partial_{t_1} 
p_0(t_1-s+\varepsilon,x,y)\, Q_sf(y)\, dy ds\,.
$
The first equality in
Lemma~\ref{lem:point-1} finally yields the result.
\qed

\begin{lemma}\label{lem:point-3}
For all $t,\varepsilon>0$, $t+\varepsilon\leq T$,
$x\in\Rd$ and $f\in C_0(\Rd)$,
\begin{align*}
\pL_x \int_0^t P_{t-s+\varepsilon}^0\, Q_s f (x)\, ds
= \int_0^t \pL_x P_{t-s+\varepsilon}^0\, Q_s f (x)\, ds\,.
\end{align*}
\end{lemma}
\pf
Similarly to Lemma~\ref{lem:point-1},
the result will follow if we can differentiate under the integrals in the first spatial variable 
and change the order of integration,
which validate the steps below, where now
$
F(x):=\int_0^t P_{t-s+\varepsilon}^0\, Q_s f (x)\, ds
$,
\begin{align*}
\pL_x F(x)
&=\exdrf(x)\cdot \nabla_x 
F(x)+
\int_{\Rd}\Big(F(x+z)-F (x)-\ind_{|z|<1} \left<z,\nabla_x F(x)\right>\!\Big)\kappa(x,z)J(z)dz\\
&=
\int_0^t\int_{\Rd}\exdrf(x)\cdot \nabla_x p^{\mathfrak{K}_y}(t-s+\varepsilon,x,y) Q_sf(y)\,dyds\\
&\quad +\int_{\Rd}\int_0^t \int_{\Rd}
\delta_1^{\mathfrak{K}_y}(t-s+\varepsilon,x,y;z)Q_sf(y)\,dyds\,
\kappa(x,z)J(z)dz\\
&=\int_0^t\int_{\Rd} \Bigg(\exdrf(x)\cdot \nabla_x p^{\mathfrak{K}_y}(t-s+\varepsilon,x,y)\\
&\qquad\qquad\qquad+
\int_{\Rd}\delta_1^{\mathfrak{K}_y}(t-s+\varepsilon,x,y;z)\kappa(x,z)J(z)dz \Bigg) Q_s f(y)\,dyds\\
&=\int_0^t\int_{\Rd} \pL_x p_0(t-s+\varepsilon,x,y)Q_s f(y) \,dyds = \int_0^t \pL_x P_{t-s+\varepsilon}Q_sf(x)\,ds\,.
\end{align*}
In the last equality we have used
Lemma~\ref{lem:point-1} with $\tilde{f}(x)=Q_sf(x)$ with
Fact~\ref{fact:PQ_well_C0}.
The differentiability
in the spatial variable
is validated since by \eqref{lem:point-1-aux} and
Fact~\ref{fact:Q_well_Bb} we get
\begin{align*}
|\frac1{\varkappa}\left(p^{\mathfrak{K}_y}(t-s+\varepsilon,x+\varkappa e_i,y)-p^{\mathfrak{K}_y}(t-s+\varepsilon,x,y)\right)Q_sf(y)|
\leq c \rr_{\varepsilon}(y-x) s^{-1} r_s^{\varepsilon_0} \|f\|_{\infty}\,,
\end{align*}
where the right hand side is integrable over $(0,t)\times \Rd$ against 
$dyds$, see 
\cite[Lemma~5.6]{MR3996792} and
Lemma~\ref{lem:time_conv}.
Fubini's theorem is in force by
bounds provided by
\eqref{ineq:L1_uni_time-2},
\eqref{ineq:L1_uni_time-3},
Fact~\ref{fact:Q_well_Bb},
\eqref{set:J},
\eqref{set:k-bound},
and
integrability properties in
\cite[Lemma~5.6]{MR3996792},
Lemma~\ref{lem:time_conv},
and that of $\nu(|z|)dz$.
In particular, we get $ \int_0^t P_{t-s+\varepsilon}^0\, Q_s f\, ds \in D(\pL)$.
\qed

We are in a position to handle
$
(\partial_t-\pL) P_{t,\varepsilon}f
$.
Note that from \eqref{def:approx_sol},
\begin{align}\label{eq:approx_sol_point}
P_{t,\varepsilon}f(x) = P_{t+\varepsilon}^0 f(x) + \int_0^t P_{t-s+\varepsilon}^0\, Q_s f(x)\, ds\,,
\end{align}
since the evaluation at a point is a continuous functional on $C_0(\Rd)$ and therefore commutes with Bochner integral, see \cite{MR3617205}.
From Lemmas~\ref{lem:point-1}, \ref{lem:point-2} and~\ref{lem:point-3} we get an immediate corollary.

\begin{corollary}\label{cor:point-1-3}
For all $t,\varepsilon>0$, $t+\varepsilon\leq T$,
$x\in\Rd$ and $f\in C_0(\Rd)$,
the mapping
$t\mapsto P_{t,\varepsilon}f(x)$ is differentiable,
$P_{t,\varepsilon}f\in D(\pL)$ and
\begin{align*}
\partial_t  P_{t,\varepsilon}f(x) &= \partial_t P_{t+\varepsilon}^0 f (x) + P_{\varepsilon}^0\, Q_t f(x)
+ \int_0^t \partial_t P_{t-s+\varepsilon}^0\, Q_s f (x)\, ds\,,\\
\pL_x  P_{t,\varepsilon}f(x)&=\pL_x P_{t+\varepsilon}^0 f (x) + \int_0^t \pL_x P_{t-s+\varepsilon}^0\, Q_s f (x)\, ds\,.
\end{align*}
\end{corollary}

\begin{corollary}
The hypothesis {\rm (H5)} is satisfied.
\end{corollary}
\pf
By Corollary~\ref{cor:point-1-3}
we can
apply the operator $(\partial_t-\pL_x)$ to
$P_{t,\varepsilon}f(x)$ in compliance with
Remark~\ref{rem:H5_H6}, and
\begin{align}\label{eq:approx_sol_der_t_L}
(\partial_t-\pL_x) P_{t,\varepsilon} f(x)
&=
\partial_t P_{t+\varepsilon}^0 f (x)+
P_{\varepsilon}^0\, Q_t f(x)
+ \int_0^t \partial_t P_{t-s+\varepsilon}^0\, Q_s f (x)\, ds \nonumber \\
&\quad - \pL_x P_{t+\varepsilon}^0 f (x) - \int_0^t \pL_x P_{t-s+\varepsilon}^0\, Q_s f (x) \nonumber \\
&= P_{\varepsilon}^0\, Q_t f(x)
-Q_{t+\varepsilon}^0 f(x) - \int_0^t Q_{t-s+\varepsilon}^0\, Q_s f(x)\,ds\,.
\end{align}
\qed

\begin{lemma}\label{lem:point-4}
For all $t,\varepsilon>0$, $t+\varepsilon\leq T$,
$x\in\Rd$ and $f\in C_0(\Rd)$,
\begin{align*}
\pL_x \int_0^t P_{s,\varepsilon} f (x)\,ds =
\int_0^t \pL_x P_{s,\varepsilon} f (x)\,ds \,.
\end{align*}
\end{lemma}
\pf
The following two equalities will entail the result (see Lemma~\ref{lem:point-1} and~\ref{lem:point-3}),
\begin{align*}
&\pL_x \int_0^t P_{s+\varepsilon}^0 f (x)\,ds =
\int_0^t \int_{\Rd} \pL_x \, p_0 (s+\varepsilon,x,y) f (y)\,dyds\,,\\
&\pL_x \int_0^t \int_0^u P_{u-s+\varepsilon}^0\, Q_s f (x)\, dsdu
= \int_0^t \int_0^u \int_{\Rd} \pL_x p_0(u-s+\varepsilon,x,y)\, Q_s f (y)\,dy dsdu\,,
\end{align*}
Similarly to
Lemma~\ref{lem:point-1} and~\ref{lem:point-3}
they hold if we can differentiate under the integrals in the first spatial variable 
and change the order of integration.
The steps to verify that are very much like those in
Lemma~\ref{lem:point-1} and~\ref{lem:point-3} (the main ingredients being Corollary~\ref{cor:L1_uni_time}, \cite[Corollary~5.10, Lemma~5.6]{MR3996792}, Fact~\ref{fact:Q_well_Bb} and Lemma~\ref{lem:time_conv}), and
therefore we omit further details. In particular, 
the integrals are absolutely convergent and
the discussed objects belong to $D(\pL)$.
\qed

\begin{corollary}\label{cor:H6}
The hypothesis {\rm (H6)} is satisfied.
\end{corollary}
\pf
The integrability conditions required for {\rm (H6)} follow 
from the formulas in Corollary~\ref{cor:point-1-3}
and Lemma~\ref{lem:point-1},
and the estimates in \eqref{ineq:L1_uni_time-1} and Fact~\ref{fact:Q_well_Bb}
combined with integrability properties from
\cite[Lemma~5.6]{MR3996792}
and Lemma~\ref{lem:time_conv}.
The equality required for {\rm (H6)} is provided by Lemma~\ref{lem:point-4}.
The regularity indicated in Remark~\ref{rem:H5_H6}
is guaranteed by
Corollary~\ref{cor:point-1-3}
and Lemma~\ref{lem:point-4}.
\qed

To sum up, we showed that for $(L,D(L))=(\pL, D(\pL))$ and $\mathcal{D}=C_0^2(\Rd)$
the conditions
\mbox{{\rm (L0)} -- {\rm (L3)}} are satisfied,
and the hypothesis 
{\rm (H5)} and {\rm (H6)}
hold true on $(0,T]$. Therefore, 
all the conclusions of Section~\ref{sec:coj-gen} are available. In particular,
Proposition~\ref{prop:coj123}
and~\ref{prop:coj4}
assert that
the conjectures {\rm (CoJ1)} -- {\rm (CoJ4)}
are true on $(0,T]$.

\vspace{\baselineskip}

\subsection{Integral kernels}\label{ssec:ker}

As concluded at the end of the Section~\ref{ssec:a-hc},
the condition
\eqref{op:S3} holds and
all the conclusions of
Section~\ref{sec:realization}
are valid for $t\in (0,T]$.
Beyond that, we have 

\begin{lemma}\label{lem:p-non-neg}
There exists 
a Borel set $\mathcal{S}_{p_+} \subseteq (0,T]\times\Rd\times\Rd$ 
such that $\mathcal{S}_{p_+}$
as well as every section $\mathcal{S}_{p_+}^{t,x}\subseteq \Rd$ are of full measure,
and for
 all $t\in (0,T]$, $x,y\in\Rd$,
\begin{align}\label{ineq:p-non-neg}
\ind_{\mathcal{S}_{p_+}}(t,x,y)\, p(t,x,y) \geq 0\,.
\end{align}
\end{lemma}
\pf
We let $\mathcal{S}_{p_+}:=\{(t,x,y)\colon p(t,x,y)\geq 0\}$, which is a Borel set by Proposition~\ref{prop:p-gen}, and so is $\mathcal{S}_{p_+}^{t,x}=\{y\colon p(t,x,y)\geq 0\}$.
Using 
Lusin's theorem, see \cite[Corollary on p. 56]{MR924157},
\eqref{eq:P_via_p}, \eqref{eq:p-gen}
and the dominated convergence theorem,
we extend {\rm (CoJ1)}:
for every $f\in B_b(\Rd)$, if $f \geq 0$, then
$$
\int_{\Rd} p(t,x,y) f(y)\,dy \geq 0\,.
$$
Let $B_n:=\{(t,x,y)\colon p(t,x,y)\leq -1/n\}$.
Clearly, each $B_n^{t,x}$ is of measure zero. Therefore, since $\ind_{B_n^{t,x}}(y)=\ind_{B_n}(t,x,y)$, we get
$\int_0^T \int_{\Rd}\int_{\Rd} \ind_{B_n}(t,x,y) p(t,x,y) dy\,dxdt=0$.
This asserts that
$\mathcal{S}_{p_+}^{t,x}$ 
and
$\mathcal{S}_{p_+}$
are of full measure.
\qed

\begin{lemma}\label{lem:CoJ5}
The conjecture {\rm (CoJ5)} holds true.
\end{lemma}
\pf Let
$f_n(x)=\varphi(x/n)$,
where $\varphi\in C_c^{\infty}(\Rd)$ is
such that $0\leq \varphi(z) \leq 1$,
$\varphi(z)=1$ if $|z|\leq 1$
and
$\varphi(z)=0$ if $|z|\geq 2$.
Clearly, $f_n\to 1$. By \eqref{eq:p-gen} and the dominated convergence theorem,
we also have
$P_t f_n\to P_t 1$
and $\int_0^t P_s\, \pL f_n\,ds\to 0$.
We use {\rm (CoJ4)} for $f_n$ and let $n\to\infty$.
\qed

\begin{remark}\label{rem:t-global}
Since $T>0$ is arbitrary, we actually have that
\eqref{eq:Q_n_via_q_n},
\eqref{def:q_n-gen},
\eqref{eq:Q_via_q},
\eqref{def:q-gen},
\eqref{eq:P_via_p},
\eqref{def:p-gen},
and
\eqref{ineq:p-non-neg}
hold for all $t>0$, $x,y\in\Rd$,
where the sets
$\mathcal{S}_n, \mathcal{S}_q, \mathcal{S}_p, \mathcal{S}_{p_+} \subseteq (0,\infty)\times\Rd\times\Rd$
and all sections
$\mathcal{S}_n^{t,x},\mathcal{S}_q^{t,x}, \mathcal{S}_p^{t,x}, \mathcal{S}_{p_+}^{t,x}\subseteq \Rd$
are of full measure.
\end{remark}

\vspace{\baselineskip}

\subsection{Uniqueness}

In this section we prove Theorem~\ref{thm:uniq}.
A general idea of the proof is that the assumptions of the theorem guarantee that
\begin{align}\label{eq:1=2}
\iint\limits_{(s,\infty)\times \Rd}  \mu_1(s,x,dudz) \psi(u,z) 
= \iint\limits_{(s,\infty)\times \Rd}  \mu_2(s,x,dudz)  \psi(u,z) \,,
\end{align}
where
\begin{align}\label{eq:phi-psi}
\big[ \partial_u + \pL_z \big]\phi(u,z)=\psi(u,z)\,.
\end{align}
The result shall follow if 
the class of $\psi$ for which
\eqref{eq:1=2} holds is rich enough.
Put differently, we intend to show that given $\psi$ the equation \eqref{eq:phi-psi} can be solved, that is, we can find $\phi$.
Since
$P_t$ is a candidate for the fundamental solution of $\partial_t u =\pL u$,
we anticipate that
$$
\phi(u,z)=-\int_{u}^{\infty} P_{r-u}\psi(r,\cdot) (z)\,dr\,.
$$
The problem is that our knowledge about the regularity of $P_t  \psi$ 
is limited, we therefore only use an approximate form of $\phi$ created by substituting $P_{t,\varepsilon}\psi$ for the integrand.
More details can be found in the proof.
We point out that in the proof we use various (arbitrarily large) values of~$T$.

\vspace{\baselineskip}

\noindent
{\bf Proof of Theorem~\ref{thm:uniq}.}
We divide the proof into five steps.

\noindent
{\it Step 1.} 
We claim that
\begin{align}\label{eq:uniq-ext}
\iint\limits_{(s,\infty)\times \Rd}  \mu_j(s,x,dudz)\Big[\partial_u\, \Phi(u,z) + \pL_z \Phi(u,z)  \Big] = - \Phi(s,x)\,,
\qquad j=1,2\,,
\end{align}
holds for every function $\Phi\colon \RR\times \Rd \to \Rd$ which satisfies the following conditions:
\begin{enumerate}
\item[(i)] there are $\mathcal{S}<s<\mathcal{T}$ such that $\Phi(u,z)=0$ if $u<\mathcal{S}$ or $u> \mathcal{T}$,
\item[(ii)] $\partial_u \Phi(u,z)$ exists and is continuous on $\RR\times\Rd$, and $\|\partial_u \Phi\|_{\infty}<\infty$,
\item[(iii)] $\partial_z^{\bbbeta}\Phi(u,z)$ exists and is continuous on $\RR\times\Rd$, and $\|\partial_z^{\bbbeta} \Phi\|_{\infty}<\infty$ for every $|\bbbeta|\leq 2$.
\end{enumerate}
Note that such $\Phi$ is bounded and continuous on $\RR\times\Rd$. We consider
$$
\phi_n(u,z)= (\Phi * \varphi_{n})(u,z)\cdot \varphi(z/n)\,,
$$
where $\varphi_{n}$ is a standard mollifier ('$*$' stands for convolution) in $\RR\times\Rd$, and  $\varphi\in C_c^{\infty}(\Rd)$ 
is an arbitrary fixed function
such that $0\leq \varphi \leq 1$,
$\varphi(z)=1$ if $|z|\leq 1$
and
$\varphi(z)=0$ if $|z|\geq 2$.
It is not hard to verify that
$\phi_n\in C_c^{\infty}(\RR\times\Rd)$
and that
the following pointwise
convergence holds
$$
\phi_n \longrightarrow \Phi\,, \qquad
\partial_u \phi_n \longrightarrow \partial_u \Phi\,, \qquad
\pL_z \phi_n \longrightarrow \pL_z \Phi\,,
$$
and that there is a constant $c>0$ such that
$\|\partial_u \phi_n\|_{\infty}+\|\pL_z \phi_n\|_{\infty} \leq c <\infty$
for all $n\in\NN$.
Now, we use
\eqref{eq:uniq} with $\phi=\phi_n$ and the dominated convergence theorem to pass with $n$ to infinity under the integral, which
yields \eqref{eq:uniq-ext}.

\vspace{\baselineskip}

\noindent
{\it Step 2.}
Let $\theta\in C_c^{\infty}(\RR)$, $\xi\in C_c^{\infty}(\Rd)$.
We claim that for each $\varepsilon\in (0,1]$ the function $\widetilde{\Phi}_{\varepsilon}$ defined below satisfies (i), (ii), (iii).
We consider
$$
\widetilde{\Phi}_{\varepsilon}(u,z):=\vartheta(u) \Phi_{\varepsilon}(u,z)\,,
$$
where $\vartheta\in C^{\infty}(\RR)$ is such that $\vartheta(u)=0$ if $u\leq \mathcal{S}$, and $\vartheta(u)=1$ if $u\geq s$
for some
$\mathcal{S}<s$, and
\begin{align*}
\Phi_{\varepsilon}(u,z):=
-\int_{u}^{\infty} \theta(r) P_{r-u,\varepsilon}\,\xi(z)\,dr\,.
\end{align*}
Let $\mathcal{T} > s$ be such that ${\rm supp}\, \theta \subseteq (-\infty,\mathcal{T})$.
Note that
to have $\widetilde{\Phi}_{\varepsilon}$ properly defined we only need $\Phi_{\varepsilon}(u,z)$ for $u\geq \mathcal{S}$. Therefore,
we use $P_t$ constructed on $(0,T]$, where
$T=\mathcal{T}-\mathcal{S}+1$,
which gives access to $P_{t,\varepsilon}$ for $t\in (0,T-1]$ and $\varepsilon \in (0,1]$,
see Remark~\ref{rem:1toT}.
Clearly, (i) holds for $\widetilde{\Phi}_{\varepsilon}$ with $\mathcal{S}$ and $\mathcal{T}$.
Hence, in what follows, we consider $\widetilde{\Phi}_{\varepsilon}(u,z)$ and $\Phi_{\varepsilon}(u,z)$ only for $\mathcal{S}\leq u \leq \mathcal{T}$.
We focus on~(ii). Note that
$$
\Phi_{\varepsilon}(u,z)=-\int_0^{\mathcal{T}-\mathcal{S}} \theta(v+u) P_{v,\varepsilon} \,\xi(z) \,dv\,, \qquad
\partial_u \Phi_{\varepsilon}(u,z)=-\int_0^{\mathcal{T}-\mathcal{S}} \theta'(v+u) P_{v,\varepsilon} \,\xi(z) \,dv\,.
$$
In the second equality above we use
Corollary~\ref{cor:approx_sol-gen_prop}
part (4) to differentiate under the integral. It also guarantees 
continuity and boundedness of $\Phi_{\varepsilon}$ and $\partial_u \Phi_{\varepsilon}$,
hence also of
$\partial_u \widetilde{\Phi}_{\varepsilon}(u,z)
= 
\vartheta'(u) \Phi_{\varepsilon}(u,z)
+\vartheta(u) \partial_u \Phi_{\varepsilon}(u,z)$.
We shall now consider (iii).
By
\eqref{lem:point-1-aux},
Fact~\ref{fact:Q_well_Bb},
\cite[Lemma~5.6]{MR3996792},
Lemma~\ref{lem:time_conv}
and the dominated convergence theorem, we have for $|\bbbeta|\leq 1$,
\begin{align*}
\partial_z^{\bbbeta} \int_0^{\mathcal{T}-\mathcal{S}} \theta(v+u) P^0_{v+\varepsilon} \,\xi(z) \,dv
=\int_0^{\mathcal{T}-\mathcal{S}} \theta(v+u) \int_{\Rd} \partial_z^{\bbbeta} \,p_0(v+\varepsilon,z,y)\xi(y) \,dy\,dv
\end{align*}
and
\begin{align*}
\partial_z^{\bbbeta} \int_0^{\mathcal{T}-\mathcal{S}} \theta(v+u) & \int_0^v P^0_{v-r+\varepsilon} Q_r \,\xi(z)\,drdv\\
=&
\int_0^{\mathcal{T}-\mathcal{S}} \theta(v+u) \int_0^v \int_{\Rd}\partial_z^{\bbbeta}\, p_0(v-r+\varepsilon,z,y) Q_r\, \xi(y)\,dy\,drdv\,.
\end{align*}
Similar argument using \eqref{ineq:L1_uni_time-1}
and \cite[Corollary~5.10]{MR3996792} proves the above formulas for all $|\bbbeta|\leq 2$.
Clearly, by \eqref{def:approx_sol}, $\partial_z^{\bbbeta} \,\Phi_{\varepsilon}$
is the sum of the above expressions.
Therefore, the same estimates that allowed differentiation under the integrals,
yield
continuity (see \eqref{ineq:remove_x}) and boundedness of
$\partial_z^{\bbbeta} \,\Phi_{\varepsilon}$, see also Lemma~\ref{lem:a-continuity}.
Finally, the same holds for
$\partial_z^{\bbbeta}\, \widetilde{\Phi}_{\varepsilon}(u,z)=\vartheta(u) \partial_z^{\bbbeta} \,\Phi_{\varepsilon}(u,z)$.

\vspace{\baselineskip}

\noindent
{\it Step 3.} 
We claim that
for all $s\leq u \leq \mathcal{T}$, $z\in\Rd$,
\begin{align*}
\big[\partial_u+\pL_z\big] \widetilde{\Phi}_{\varepsilon}(u,z)= \theta(u) P_{0,\varepsilon}\,\xi(z)+ \int_0^{\mathcal{T}-\mathcal{S}} \theta(v+u)\, (\partial_v - \pL_z) P_{v,\varepsilon}\,\xi(z)\,dv\,.
\end{align*}
Note that $\widetilde{\Phi}_{\varepsilon}=\Phi_{\varepsilon}$ if $u\geq s$.
Due to
Corollary~\ref{cor:H6}, see also
Remark~\ref{rem:H5_H6},
and
Corollary~\ref{cor:approx_sol-gen_prop}
part (4)
we can integrate by parts to get
$\partial_u \Phi_{\varepsilon}(u,z)
=\theta(u)P_{0,\varepsilon}\,\xi(z)+\int_0^{\mathcal{T}-\mathcal{S}} \theta(v+u)\partial_v P_{v,\varepsilon}\,\xi(z)\,dv$.
Furthermore, we have
$\pL_z \Phi(u,z)=-\int_0^{\mathcal{T}-\mathcal{S}} \theta(v+u)\, \pL_z P_{v,\varepsilon}\,\xi(z)\,dv$,
by literally the same proof as that of Lemma~\ref{lem:point-4}.
Adding the two expressions gives the claim.

\vspace{\baselineskip}

\noindent
{\it Step 4.}
Let $\theta\in C_c^{\infty}(\RR)$, $\xi\in C_c^{\infty}(\Rd)$.
We claim that \eqref{eq:1=2}
holds with $\psi(u,z)=\theta(u)\xi(z)$.
From~{\it Step 1} and {\it Step 2} we get
\eqref{eq:1=2} with
$\psi(u,z)=[\partial_u+\pL_z] \widetilde{\Phi}_{\varepsilon}(u,z)$.
Next, by
{\it Step 3},
{\rm (H0)} (see Corollary~\ref{cor:H0})
and
Fact~\ref{fact:approx_sol-zero-2}
we get
$\lim_{\varepsilon\to 0^+}[\partial_u+\pL_z] \widetilde{\Phi}_{\varepsilon}(u,z)=\theta(u)\xi(z)$. Therefore, by the dominated convergence theorem we obtain the claim.

\vspace{\baselineskip}

\noindent
{\it Step 5.} We finally claim that the consequent in Theorem~\ref{thm:uniq} is true.
Clearly, by {\it Step 4} and approximation by standard mollification,
we extend \eqref{eq:1=2}
to hold for every $\psi(u,z)=\theta(u)\xi(z)$,
where $\theta\in C_c(\RR)$, $\xi\in C_c(\Rd)$.
Next, by Urysohn's lemma and the dominated convergence theorem, \eqref{eq:1=2} holds for every 
$\psi(u,z)=\ind_{J}(u) \ind_{V}(z)$,
where $J\subset \RR$ is open and bounded, and $V\subseteq \Rd$ is open, 
see \cite[(2) on p. 49]{MR924157}.
Finally, we get \eqref{eq:1=2} for $\psi(u,z)=\ind_{E}(u,z)$, where $E$ is like in the statement,
by the monotone class theorem
\cite[Theorem~2.3]{MR0264757}.
\qed

\subsection{Existence and properties}
In this section we prove
Theorem~\ref{thm:exist} and Theorem~\ref{thm:sem_prop}.

\vspace{\baselineskip}

Due to
Sections~\ref{ssec:a-hc}, \ref{ssec:b-hc},
\ref{ssec:ker}
we can use all the conclusions of results in Sections~\ref{sec:f_a_a},
\ref{sec:coj-gen}, \ref{sec:realization}.

\vspace{\baselineskip}

\noindent
{\bf Proof of Theorem~\ref{thm:sem_prop}.}
Since $T>0$ is arbitrary and the formulas 
\eqref{def:op_P}, \eqref{def:op_Q}
and \eqref{def:op_Qn}
are consistent, the operators $P_t$ are defined for all $t>0$. Furthermore,
\eqref{eq:P_via_p}
and
\eqref{def:p-gen}
hold for all $t>0$, $x,y\in\Rd$, see
Remark~\ref{rem:t-global}.
We also use
Remark~\ref{rem:1toT}.
The operator $P_t$ maps from $C_0(\Rd)$ to $C_0(\Rd)$ by Fact~\ref{fact:PQ_well_C0}.
The semigroup property and
positivity stem from Proposition~\ref{prop:coj123},
where {\rm (CoJ1)} and {\rm (CoJ3)} are validated.
The strong continuity follows from Fact~\ref{fact:P-s-cont}.
The contractivity results from {\rm (CoJ1)} and {\rm (CoJ2)}, therefore again from Proposition~\ref{prop:coj123}.
Furthermore,

part (i) follows from Lemma~\ref{lem:CoJ5},

part (ii) follows from Fact~\ref{fact:strong-F}. see also Lemma~\ref{lem:strong-F}, 

part (iii) follows from {\rm (CoJ4)} validated in Proposition~\ref{prop:coj4}, see \eqref{def:DpL} and 
\eqref{def:pD},

part (iv) and (v)  follow from Proposition~\ref{prop:diff_closure}, proved in the next section.
\qed

\noindent
{\bf Proof of Theorem~\ref{thm:exist}.}
Let $\phi\in C_c^{\infty}(\RR\times\Rd)$.
Then $\xi(t)=\phi(t,\cdot)$
satisfies the assumptions of
\cite[Theorem~4.1]{MR3514392},
see \cite[Proof of (1.16)]{MR3514392}.
Note that here we 
use Theorem~\ref{thm:sem_prop}.
Thus
for every $(s,x)\in\RR\times\Rd$ and $\phi \in C_c^{\infty}(\RR\times\Rd)$ we have
\begin{align}\label{eq:fs}
\int_s^{\infty} P_{u-s} \big[ \partial_u \phi(u,\cdot)+ \pL\, \phi(u,\cdot) \big](x)\, du = -\phi(s,x)\,.
\end{align}
This ends the proof by \eqref{eq:P_via_p} and \eqref{eq:p-gen}.
\qed

\subsection{Time-derivative and generator}
In this section we calculate the derivative of $P_tf(x)$ with respect to $t>0$, and use it to analyse the generator.
The calculations are based on
\begin{align}\label{eq:P-alt}
P_t f=P_t^0 f+\int_0^{t/2} P_{t-s}^0\, Q_s f\, ds+\int_0^{t/2}P_s^0\, Q_{t-s} f\, ds\,,
\end{align}
the series representation $Q_tf=\sum_{n=0}^{\infty}Q_t^n f$, and for $n=1,\ldots$, splitting the integrals:
\begin{align}\label{eq:Qn-alt}
Q_t^{n} f = \int_0^{t/2} Q_{t-s}^0\, Q_{s}^{n-1}f \,ds+\int_0^{t/2} Q_{s}^0\, Q_{t-s}^{n-1}f \,ds\,.
\end{align}
See \eqref{def:op_P}, \eqref{def:op_Q} and \eqref{eq:Qn}.
In what follows the constant $C_3$ is taken from
Lemma~\ref{lem:well_def}.

\begin{lemma}\label{lem:der_q0-t}
There is $c>0$ such that for all $t\in (0,T]$, $x\in\Rd$ and $f\in B_b(\Rd)$,
$$
\partial_t Q_t^0 f(x)=\int_{\Rd} \partial_t q_0(t,x,y) f(y)\,dy\,,
$$
and
$$
\| \partial_t Q_t^0 f \|_{\infty}
\leq c \,C_3
(t^{-1 }r_t^{\efcs-1})t^{-1}r_t^{\varepsilon_0}\|f\|_{\infty}\,.
$$
\end{lemma}
\pf
The equality follows from 
\eqref{eq:delta-alt}, Remark~\ref{rem:gen_not},
\eqref{ineq:L1_uni_time-1}, \cite[Lemma~5.6]{MR3996792} and
the dominated convergence theorem.
The inequality follows from 
Lemmas~\ref{lem:q_0-aux*1}, \ref{lem:q_0-aux}
and~\ref{lem:q_0-aux*2}
together with Corollary~\ref{cor-shifts} and Lemma~\ref{lem:conv}(a), cf. proof of Lemma~\ref{lem:well_def}.
\qed

We recall that the value of $\varepsilon_0$ is fixed and defined by \eqref{def:ve_0-A} and \eqref{def:ve_0-A*}.

\begin{lemma}\label{lem:der_qn-t}
There are  $C,c>0$ such that for all $t\in (0,T]$, $x\in\Rd$,
$n\in\NN$ and $f\in B_b(\Rd)$,
$$
\partial_t Q_t^nf(x)= \int_0^{t/2} \partial_t Q_{t-s}^0\, Q_{s}^{n-1}f (x) \,ds+\int_0^{t/2} Q_{s}^0\, \partial_t Q_{t-s}^{n-1}f(x) \,ds+Q_{t/2}^0\, Q_{t/2}^{n-1}f(x)\,,
$$
and
$$
\| \partial_t Q_t^nf\|_{\infty}
\leq c \,C_3 (t^{-1}r_t^{\efcs-1})\,
(C C_3)^n \, \prod_{k=1}^n B\!\left(\frac{\varepsilon_0}{2},\frac{k\varepsilon_0}{2}\right) t^{-1} r_t^{(n+1)\varepsilon_0}\,
  \|f\|_{\infty}\,.
$$
\end{lemma}
\pf
Let $c$ be the constant from Lemma~\ref{lem:der_q0-t}. Note that
by \eqref{set:h-scaling} there is $c_1=c_1(T,\lah,\lch,h)>0$ such that $r_{t/2} \geq c_1 r_t$ for $t\in(0,T]$, see \cite[Lemma~2.3]{MR4140542}.
We let $\lambda_0=1$, and for $n\in\NN$ we put
$$
\lambda_n=  C_3^n \prod_{k=1}^n B\!\left(\frac{\varepsilon_0}{2},\frac{k\varepsilon_0}{2}\right)
\qquad \mbox{and}\qquad
C=3\max\left\{1, 2 c_1^{\efcs-1},\frac{4 r_T^{1-\efcs}}{c\Gamma(\frac{\varepsilon_0}{2})}2^{\frac{\varepsilon_0}{2}}\right\}.
$$
Our aim is to prove
differentiability
and 
$\|\partial_t Q_t^n f\|_{\infty}\leq c\, C_3 (t^{-1}r_t^{\efcs-1})\, C^n\lambda_n   \, t^{-1} r_t^{(n+1)\varepsilon_0}\,
  \|f\|_{\infty}$ for $n\in \NN_0$.
If $n=0$ the statement holds by Lemma~\ref{lem:der_q0-t}.
We assume the differentiability and the estimate for $n\in\NN_0$.
By \eqref{eq:Qn-alt}
we have
$Q_t^{n+1} f = \int_0^{t/2} Q_{t-s}^0\, Q_{s}^{n}f \,ds+\int_0^{t/2} Q_{s}^0\, Q_{t-s}^{n}f \,ds$. 
In what follows we use properties of $Q_t^nf$ from Fact~\ref{fact:Qn_well_Bb}.
We have
\begin{align}\label{eq:d1}
\begin{aligned}
&\int_0^{t/2} Q_{t-s}^0\, Q_{s}^{n}f (x)\,ds
-\int_0^{t_1/2} Q_{t_1-s}^0\, Q_{s}^{n}f(x) \,ds\\
&\qquad = \int_0^{t_1/2} \int_{t_1}^t \partial_{\theta} Q_{\theta-s}^0\, Q_s^n f (x)\,d\theta ds
+\int_{t_1/2}^{t/2} Q_{t-s}^0\, Q_{s}^{n}f (x)\,ds=: I_1+I_2\,,
\end{aligned}
\end{align}
and
\begin{align}\label{eq:d2}
\begin{aligned}
&\int_0^{t/2} Q_{s}^0\, Q_{t-s}^{n}f(x) \,ds
-\int_0^{t_1/2} Q_{s}^0\, Q_{t_1-s}^{n}f(x) \,ds\\
&\qquad = \int_0^{t_1/2}\int_{\Rd} q_0(s,x,y) \int_{t_1}^t\partial_{\theta} Q_{\theta-s}^nf(y)\,d\theta \,dyds
+\int_{t_1/2}^{t/2} Q_s^0 Q_{t-s}^nf(x)\,ds
=: I_3+I_4\,.
\end{aligned}
\end{align}
Dividing by $t-t_1$, using the induction hypothesis and \eqref{ineq:H3}
justify the usage of the dominated convergence theorem for parts $I_1$ and $I_3$ as $t\to t_1$.
The continuity of $Q_{u}^0Q_v^n f$ for $u,v>0$
gives that the parts corresponding to $I_2$ and $I_4$ converge to $\frac12 Q_{t_1/2}^0\, Q_{t_1/2}^{n}f(x)$ each.
Thus we get the differentiability and
\begin{align*}
\partial_t Q_t^{n+1}f(x)= \int_0^{t/2} \partial_t Q_{t-s}^0\, Q_{s}^{n}f (x) \,ds+\int_0^{t/2} Q_{s}^0\, \partial_t Q_{t-s}^{n}f(x) \,ds+Q_{t/2}^0\, Q_{t/2}^{n}f(x)\,.
\end{align*}
We then prove the estimates. We frequently use
Lemma~\ref{lem:time_conv}
and the monotonicity of $r_t$, see Remark~\ref{rem:r_t}.
By Lemma~\ref{lem:der_q0-t} and
Fact~\ref{fact:Qn_well_Bb},
\begin{align*}
\int_0^{t/2} \|\partial_t Q_{t-s}^0\, Q_{s}^{n}f \|_{\infty}\,ds
&\leq    c\, C_3  (t/2)^{-1} r_{t/2}^{\efcs-1}\, 
 C_3 \lambda_n \int_0^t (t-s)^{-1} r_{t-s}^{\varepsilon_0}\, s^{-1} r_s^{(n+1)\varepsilon_0}ds\\
&\leq c\, C_3 (t^{-1}r_{t}^{\efcs-1})\,  2 c_1^{\efcs-1} \lambda_{n+1} \, t^{-1} r_t^{(n+2)\varepsilon_0}\,.
\end{align*}
By
Lemma~\ref{lem:well_def} and inductive hypothesis
\begin{align*}
\int_0^{t/2}\| Q_{s}^0\, \partial_t Q_{t-s}^{n}f\|_{\infty} \,ds
&\leq  C_3\, c\, C_3 (t/2)^{-1} r_{t/2}^{\efcs-1}\, C^n \lambda_n \int_0^t s^{-1} r_s^{\varepsilon_0} (t-s)^{-1} r_{t-s}^{(n+1)\varepsilon_0}ds\\
&\leq c\,C_3 (t^{-1}r_{t}^{\efcs-1})\, 2 c_1^{\efcs-1} C^n \lambda_{n+1} \, t^{-1} r_t^{(n+2)\varepsilon_0}\,.
\end{align*}
Following \cite[p. 129]{MR3765882} we have $\frac1{B(\varepsilon_0/2,(n+1)\varepsilon_0/2)}\leq 2^{(n+1)\frac{\varepsilon_0}{2}}/\Gamma(\frac{\varepsilon_0}{2})$, $n=0,1,\ldots$.
Now, by
Lemma~\ref{lem:well_def}
and Fact~\ref{fact:Qn_well_Bb} we get
\begin{align*}
\|Q_{t/2}^0\, Q_{t/2}^{n}f\|_{\infty}
& \leq C_3 (t/2)^{-1} r_{t/2}^{\varepsilon_0}\,C_3 \lambda_n \,(t/2)^{-1}\, r_{t/2}^{(n+1)\varepsilon_0}\\
&\leq C_3 (t^{-1}r_{t}^{\efcs-1})\, 4 r_T^{1-\efcs}\,  C_3 \lambda_n \, t^{-1}r_t^{(n+2)\varepsilon_0}\\
&\leq C_3 (t^{-1} r_{t}^{\efcs-1}) \, 
 \frac{4 r_T^{1-\efcs}}{\Gamma(\frac{\varepsilon_0}{2})} 2^{(n+1)\frac{\varepsilon_0}{2}}
\lambda_{n+1}\, t^{-1} r_{t}^{(n+2)\varepsilon_0}\\
&\leq c \,C_3 
(t^{-1} r_{t}^{\efcs-1}) 
\left(\frac{4 r_T^{1-\efcs}}{c\Gamma(\frac{\varepsilon_0}{2})} 2^{\frac{\varepsilon_0}{2}}\right)  C^n 
\lambda_{n+1}
\, t^{-1} r_{t}^{(n+2)\varepsilon_0}\,.
\end{align*}
Summing up the estimates gives the desired bound for $n+1$.
\qed

\begin{lemma}\label{lem:der_q-t}
There is $c>0$ such that for all $t\in (0,T]$, $x\in\Rd$ and $f\in B_b(\Rd)$,
$$
\partial_t Q_tf(x)=\sum_{n=0}^{\infty} \partial_t Q_t^n f(x)\,,
$$ 
and
$$
\| \partial_t Q_t f\|_{\infty}
\leq c (t^{-1}r_t^{\efcs-1})\, t^{-1} r_t^{\varepsilon_0}\, \|f\|_{\infty}\,.
$$
\end{lemma}
\pf
Fix $f$ and $x$.
For every $t\in (0,T]$, by Fact~\ref{fact:Q_well_Bb} the series
$\sum_{n=0}^{\infty} Q_t^n f(x)$ converges. By
Lemma~\ref{lem:der_qn-t},
the series $\sum_{n=0}^{\infty} \partial_t Q_t^n f(x)$
converges locally uniformly in $t\in (0,T]$.
Thus,
we can differentiate term by term.
The estimate follows from Lemma~\ref{lem:der_qn-t}.
\qed

\begin{lemma}
There is $c>0$ such that for all $t\in (0,T]$ and $f\in C_0(\Rd)$,
\begin{align}\label{eq:der_p0-t-b}
\|\partial_t P_t^0f \|_{\infty}\leq c t^{-1} r_t^{\efcs-1} \|f\|_{\infty}\,,
\end{align}
and
\begin{align}\label{eq:der_p0-t-c}
\lim_{s\to t}\|\partial_s P_{s}^0f -  \partial_t P_{t}^0f\|_{\infty}=0\,.
\end{align}
\end{lemma}
\pf
The estimate \eqref{eq:der_p0-t-b} follows from Lemma~\ref{lem:point-1},
Proposition~\ref{prop:gen_est_time}, 
Corollary~\ref{cor-shifts} and Lemma~\ref{lem:conv}(a).
Now, by
\eqref{eq:par_t_is_L} and \eqref{ineq:L1_uni_time-1}
for $t,s\geq \tau$,
$$
|\partial_s p_0(s,x,y)-\partial_t p_0(t,x,y)|
\leq c|t-s| \rr_{\tau} (y-x)\,.
$$
Thus, using Lemma~\ref{lem:point-1}
and \cite[Lemma~5.6]{MR3996792}
the equality \eqref{eq:der_p0-t-c} holds.
\qed

\begin{theorem}\label{thm:der_p-t}
For all $t\in (0,T]$, $x\in\Rd$ and $f\in C_0(\Rd)$,
\begin{align}
\partial_t P_t f(x)
&=\partial_t P_t^0 f(x)
+\int_0^{t/2} \partial_t P_{t-s}^0 Q_sf(x)\,ds
+\int_0^{t/2} P_s^0 \partial_t Q_{t-s}f(x)\,ds
+ P_{t/2}^0\, Q_{t/2}f(x)\label{eq:der1} \\
&=\partial_t P_t^0 f(x)+Q_tf(x)+\lim_{\varepsilon\to 0^+}\int_0^{t-\varepsilon}\partial_t P_{t-s}^0 Q_sf(x)\,ds \label{eq:der2}\\
&= \pL_x P_t^0 f(x) + \lim_{\varepsilon \to 0^+} \pL_x \int_0^{t-\varepsilon}  P_{t-s}^0\, Q_sf(x)\,ds\,. \label{eq:der3}
\end{align}
\end{theorem}
\pf
First recall that $Q_tf\in C_0(\Rd)$ and $\partial_t Q_tf \in B_b(\Rd)$, see Fact~\ref{fact:PQ_well_C0} and Lemma~\ref{lem:der_q-t}.
We get \eqref{eq:der1} by differentiating \eqref{eq:P-alt}.
The exact calculations are like in \eqref{eq:d1} and \eqref{eq:d2}
with $Q_t^0$ replaced by $P_t^0$, and $Q_t^n$ changed to $Q_t$.
The usage of the dominated convergence theorem is justified by
\eqref{eq:der_p0-t-b},
Fact~\ref{fact:Q_well_Bb}  part (2),
\eqref{ineq:H1} and
Lemma~\ref{lem:der_q-t}.
See also Lemma~\ref{lem:continuity} part~(5).
Note that we have
$
\partial_s P_{t-s}^0 Q_sf(x)
= -\partial_t P_{t-s}^0 Q_sf(x)
+ P_{t-s}^0 \partial_s Q_sf(x)
$.
Thus
\begin{align*}
\int_0^{t/2} P_s^0 \partial_t Q_{t-s}f(x)\,ds
&=\lim_{\varepsilon \to 0^+}\int_{t/2}^{t-\varepsilon} P_{t-s}^0 \partial_s Q_sf(x)\,ds\\
&=\lim_{\varepsilon \to 0^+} \left( \int_{t/2}^{t-\varepsilon} \partial_t P_{t-s}^0 Q_sf(x)\,ds
+P_{\varepsilon}^0 Q_{t-\varepsilon}f(x)-P_{t/2}^0Q_{t/2}f(x)\right).
\end{align*}
By part (5) of Lemma~\ref{lem:continuity},
this already proves \eqref{eq:der2}.
Recall that $(\partial_t-\pL_x) P_t^0f(x)=- Q_t^0 f (x)$ by  Lemma~\ref{lem:point-1}.
Then
$$
\partial_t P_t^0 f(x)= -Q_t^0f(x)+\pL_x P_t^0f(x)\,,
$$
and
$$
\lim_{\varepsilon\to 0^+}\int_0^{t-\varepsilon}\partial_t P_{t-s}^0 Q_sf(x)\,ds=
-\int_0^{t} Q_{t-s}^0 Q_sf(x)\,ds
+\lim_{\varepsilon\to 0^+}\int_0^{t-\varepsilon} \pL_x P_{t-s}^0 Q_sf(x)\,ds\,.
$$
It suffices to add the above equalities and use \eqref{eq:Q} to obtain \eqref{eq:der3}.
\qed

\begin{lemma}\label{lem:der-t-approx}
For all 
$t\in (0,T)$ and $f\in C_0(\Rd)$
$$
\lim_{\varepsilon\to 0^+}\|\partial_t P_{t,\varepsilon}f-\partial_t P_t f\|_{\infty}=0\,.
$$
\end{lemma}
\pf
Let $t\in(0,T)$ and $\varsigma>0$ be fixed such that
$t\in (0,T-\varsigma]$.
We use \eqref{eq:der1} for $\partial_t P_tf(x)$.
Similarly to
Theorem~\ref{thm:der_p-t}, we get
\begin{align*}
\partial_t P_{t,\varepsilon} f(x)
=\partial_t P_{t+\varepsilon}^0 f(x)
&+\int_0^{t/2} \partial_t P_{t-s+\varepsilon}^0 Q_sf(x)\,ds\\
&\quad+\int_0^{t/2} P_{s+\varepsilon}^0 \partial_t Q_{t-s}f(x)\,ds
+ P_{t/2+\varepsilon}^0\, Q_{t/2}f(x)\,.
\end{align*}
We have $\partial_t P_{t+\varepsilon}^0 f \to \partial_t P_{t}^0 f$ and 
$P_{t/2+\varepsilon}^0\, Q_{t/2}f \to P_{t/2}^0\, Q_{t/2}f$ in the norm by \eqref{eq:der_p0-t-c} and {\rm (H2)}, respectively.
Next, 
by \eqref{eq:der_p0-t-b}, Fact~\ref{fact:Q_well_Bb} part (2),
{\rm (H1)} and Lemma~\ref{lem:der_q-t}
for all $s\in (0,t/2)$ and 
$\varepsilon\in [0,\varsigma]$
we get
$
\|\partial_t P_{t-s+\varepsilon}^0 Q_sf\|_{\infty}
\leq c s^{-1} r_s^{\varepsilon_0} \|f\|_{\infty}$,
$\|P_{s+\varepsilon}^0 \partial_t Q_{t-s}f\|_{\infty}\leq c\|f\|_{\infty}$.
Thus the difference of the integral parts of $\partial_tP_{t,\varepsilon}f$ and $\partial_t P_t$ converge to zero in the norm by 
\eqref{eq:der_p0-t-c} and {\rm (H2)}, and
the dominated convergence theorem.
\qed

We note that \eqref{eq:der2}
and \eqref{eq:der3}
resemble formulae displayed in Corollary~\ref{cor:point-1-3}.
Due to Lemma~\ref{lem:der-t-approx}
the relation is even stronger, see also Fact~\ref{fact:approx_sol-zero-1}.

\begin{corollary}\label{cor:der_p-tb}
For all $t\in (0,T)$ and $f\in C_0(\Rd)$,
\begin{align}
\partial_t P_t f(x)&=
\partial_t P_t^0 f(x)+Q_tf(x)+\lim_{\varepsilon\to 0^+}\int_0^t\partial_t P_{t-s+\varepsilon}^0 Q_sf(x)\,ds \label{eq:der2b} \\
&=\pL_x P_t^0 f(x) + \lim_{\varepsilon \to 0^+}   \pL_x \int_0^t  P_{t-s+\varepsilon}^0\, Q_sf(x)\,ds\,.
\label{eq:der2b}
\end{align}
\end{corollary}

We next prove parts (iv) and (v) of Theorem~\ref{thm:sem_prop}. In the proof we use  various (arbitrarily large) values of~$T$.

\begin{proposition}\label{prop:diff_closure}
The semigroup $(P_t)_{t>0}$ is differentiable and its generator
$(\mathcal{A},D(\mathcal{A}))$
 is the closure of $(\pL,C_c^{\infty}(\Rd))$. In particular, for all $t>0$, $x\in\Rd$,
\begin{align*}
\mathcal{A} P_tf(x)=\partial_t P_tf(x)\,.
\end{align*}
\end{proposition}
\pf 
Let $(\mathcal{A}_c,D(\mathcal{A}_c)):=(\pL,C_c^{\infty}(\Rd))$ and denote its closure by $(\bar{\mathcal{A}}_c,D(\bar{\mathcal{A}}_c))$.
Clearly, we have $\bar{\mathcal{A}}_c\subseteq \mathcal{A}$. To obtain that $\mathcal{A}\subseteq\bar{\mathcal{A}}_c$ it suffices to show $D(\mathcal{A})\subseteq D(\bar{\mathcal{A}}_c)$.
To this end in each of the steps below we exploit that $\bar{\mathcal{A}}_c$ is a closed operator.
Note also that the differentiability follows from {\it Step 2}.

\noindent
{\it Step 1}. We claim that $f\in D(\bar{\mathcal{A}}_c)$ for every $f\in C_0^2(\Rd)$.  Indeed, by convolution with a standard 
mollifier and multiplication by a cut-off function
we get
$f_n(x)=(f*\varphi_n)(x)\cdot \varphi(x/n)$,
which is such that $f_n\in C_c^{\infty}(\Rd)$, $f_n\to f$ and
$\bar{\mathcal{A}}_c f_n = \pL f_n\to \pL f$
in the norm
(the latter follows from
$\|\pL g\|_{\infty} \leq c \sum_{|\bbbeta|\leq 2} \|\partial^{\bbbeta} g\|_{\infty}$
and
a fact that $f_n$ converges to $f$
in the norm together with all its derivatives up to order $2$).

\noindent
{\it Step 2}. We claim that $P_t f\in D(\bar{\mathcal{A}}_c)$ for every $f\in C_0(\Rd)$.
Indeed, we let $f_n=P_{t,1/n} f$, see~\eqref{def:approx_sol}. Then
$f_n\in C_0^2(\Rd)$
by \eqref{ineq:L1_uni_time-1}.
Furthermore,
$f_n \to P_t f$  by
Corollary~\ref{cor:approx_sol-gen_prop}
part (1).
Finally, by
{\it Step~1}.,
Theorem~\ref{thm:sem_prop}(iii),
Lemma~\ref{fact:approx_sol-zero-1}
and~\ref{lem:der-t-approx} we get
$\bar{\mathcal{A}}_c f_n= \mathcal{A} f_n = \pL f_n = (\pL -\partial_t)P_{t,1/n} f+\partial_t P_{t,1/n}f \to \partial_t P_tf$.

\noindent
{\it Step 3}. We claim that $f\in D(\bar{\mathcal{A}}_c)$ for every $f\in D(\mathcal{A})$.
Indeed, if $f_n=P_{1/n}f$, then $f_n\in D(\bar{\mathcal{A}}_c)$, $f_n\to f$ and, since the generator commutes on its domain with the semigroup,
$\bar{\mathcal{A}}_c f_n = \mathcal{A} f_n= \mathcal{A}P_{1/n}f = P_{1/n}\mathcal{A}f\to \mathcal{A}f$, see Fact~\ref{fact:P-s-cont}.
\qed

\subsection{Pointwise estimates}
Note that  
Lemma~\ref{lem:q_0-aux}
provides a pointwise upper bound for $q_0(t,x,y)$
for all \mbox{$t\in (0,\indTr]$} and \mbox{$x,y\in\Rd$}.
We shall now
propagate it on $q_n$ defined in \eqref{def:q_n-gen}, and then further on $q$ and $p$ given in
\eqref{def:q-gen} and
\eqref{def:p-gen}.

\begin{lemma}\label{lem:q-bound}
Assume $\Sa$ and \eqref{cond:pointwise}.
For every $\varepsilon\in (0,\boost)$
there exists a constant  $c>0$
such that for all 
$t\in (0,\indTr]$ and $x,y\in\Rd$,
\begin{align*}
|q(t,x,y)|\leq c\left(\erry{\varepsilon}{0}+\sum_{j=0}^{\efnj}\erry{\indcsi{j}-1}{\indhei{j}}\right)(t,x,y)\,.
\end{align*}
The constant $c$ can be chosen to depend only on
$d,\lmcc,\ka,\lah,\lch,h,\indnj,\boost, \varepsilon$, $\min\limits_{j=0,\ldots,\efnj} (\indcsi{j})$.
\end{lemma}
\pf
We set $T=\indTr$ for the whole proof.
Fix $\varepsilon\in (0,\boost)$ and
choose $0<\dM< \varepsilon $ such that 
$\varepsilon+\dM<\boost$.
Now, by \eqref{cond:pointwise} for each $j=0,\ldots,\efnj$ there exists $0<\ell_j<(\lah/2)\land (\efcs \indhei{j})$ such that
\begin{align}\label{eq:def-Lj}
\varepsilon+\dM=2(\ell_j+\indcsi{j}-1)\,.
\end{align}
Define $\ell:=\max\{\ell_j\colon\,\, j=0,\ldots,\efnj\}$. Note that
$$
0<\ell<\lah/2\,,\qquad \dM < \min_{j=0,\ldots,\efnj}( \ell_j+\indcsi{j}-1) \,.
$$
We will show by induction that for all $n\in\NN$ and
$t\in (0,T]$, $x,y\in\Rd$, 
\begin{align}\label{ineq:q_n-bound}
|q_n(t,x,y)|\leq \lambda_n\left(\erry{\varepsilon+n \dM}{0}+ \sum_{j=0}^{\efnj} \erry{\indcsi{j}-1+n\dM}{\indhei{j}}\right)(t,x,y)\,,
\end{align}
where
$$
\lambda_n= c_0 \left(\frac{c}{\lah-2\ell} \right)^n \, \prod_{k=1}^{n} B(\dM/2,k\dM/2)\,,
$$
a constant $c_0$ is taken from Lemma~\ref{lem:q_0-aux}
and $c$ is another constant that we will now specify.
To this end we write 
$D_1$ for a constant from Corollary~\ref{cor-shifts}, similarly
$D_2$ is from Lemma~\ref{lem:conv}(b),
$D_3=\lch (r_T\vee 1)^{2}$, $D_4=r_T\vee 1$, $D_5=5(\efnj+1)^2c_0$. Finally we write 
$$c:=D_1 D_2 D_3^{2/\lah} D_4 D_5\,.$$

Before we proceed with the proof of \eqref{ineq:q_n-bound} we establish two auxiliary inequalities which will facilitate calculations.
For the first auxiliary inequality we let
$n\in\NN_0$
and apply Lemma~\ref{lem:conv}(c) with
\begin{align*}
\beta_0=2\ell,\ \beta_1= \indhei{j},\ \beta_2 = \indhei{i},\quad 
&n_1=n_2=\varepsilon+\dM+2-\indcsi{j}-\indcsi{i}\,,\\
&m_1=\dM+1-\indcsi{j}, \ m_2=\dM+1-\indcsi{i}\,,\\
&\gamma_1 =\indcsi{j}-1, \ \gamma_2=\indcsi{i}-1+n\dM\,.
\end{align*}
Note that by the definition of $\ell_j$ in \eqref{eq:def-Lj} and the subsequent properties of $\dM$, the assumptions of Lemma~\ref{lem:conv}(c) are satisfied. Indeed,  for all $i,j=0,\ldots, \efnj$ we have 
\begin{align*}
&0<\varepsilon+\dM+2-\indcsi{j}-\indcsi{i}=\ell_j+\ell_i\leq (2\ell)\wedge\efcs(\indhei{j}+\indhei{i})\,,\\
&0<\dM+1-\indcsi{j}<\ell_j \leq (2\ell)\wedge\efcs\indhei{j}\,.
\end{align*}
Therefore we get
\begin{align}
&\int_0^t\int_{\Rd}
\errx{\indcsi{j}-1}{\indhei{j}}(t-s,x,z) \erry{\indcsi{i}-1+n\dM}{\indhei{i}}(s,z,y)\,dzds  \nonumber \\
&\leq  
\frac{D_2 D_3^{2/\lah}}{\lah-2\ell}\,
B\!\left(\frac{\dM}{2},\frac{(n+1)\dM}{2}\right) \left(2\erry{\varepsilon+(n+1)\dM}{0}+ \erry{\indcsi{i}-1+(n+1)\dM}{\indhei{i}}
+\erry{\indcsi{j}-1+(n+1)\dM}{\indhei{j}} \right)(t,x,y)\,. \label{ineq:aux1}
\end{align}
For the second auxiliary inequality
 we use Lemma~\ref{lem:conv}(c) with
\begin{align*}
\beta_0=2\ell,\ \beta_1=\indhei{j}, \ \beta_2=0,  \quad 
&n_1= n_2=\dM+1-\indcsi{j}\,,\\
&m_1=\dM+1-\indcsi{j}, \  m_2=0\,,\\
&\gamma_1 =\indcsi{j}-1, \  \gamma_2=\varepsilon+n\dM\,.
\end{align*}
The assumptions are satisfied again and we get
\begin{align}
&\int_0^t\int_{\Rd}
\errx{\indcsi{j}-1}{\indhei{j}}(t-s,x,z) \erry{\varepsilon+n\dM}{0}(s,z,y)\,dzds \nonumber \\
&\leq \frac{D_2 D_3^{1/\lah}}{\lah-2\ell}\,
B\!\left(\frac{\dM}{2},\frac{\varepsilon+n\dM}{2}\right)\left(3\erry{\varepsilon+(n+1)\dM}{0}+\erry{\varepsilon+\indcsi{j}-1+n\dM}{\indhei{j}}\right)(t,x,y) \nonumber \\
&\leq \frac{D_2 D_3^{1/\lah}D_4}{\lah-2\ell}\,B\!\left(\frac{\dM}{2},\frac{(n+1)\dM}{2}\right)\left(3\erry{\varepsilon+(n+1)\dM}{0}+\erry{\indcsi{j}-1+(n+1)\dM}{\indhei{j}}\right)(t,x,y)\,.\label{ineq:aux2}
\end{align}
In the second inequality we used the monotonicity of the Beta function, we wrote that
$\varepsilon+n\dM=\varepsilon-\dM+(n+1)\dM $, applied Remark~\ref{rem:r_t} and $D_4^{\varepsilon-\dM}\leq D_4$.

We will now show \eqref{ineq:q_n-bound}
for $n=1$.
By Lemma~\ref{lem:q_0-aux},
Corollary~\ref{cor-shifts} and \eqref{ineq:aux1} we obtain
\begin{align*}
|q_1(t,x,y)|
&\leq \int_0^t \int_{\Rd}|q_0(t-s,x,z)| |q_0(s,z,y)|\, dzds\\
&\leq c_0^2 D_1 \sum_{i,j=0}^{\efnj}
\int_0^t\int_{\Rd}
\errx{\indcsi{j}-1}{\indhei{j}}(t-s,x,z) \erry{\indcsi{i}-1}{\indhei{i}}(s,z,y)\,dzds\\
&\leq c_0^2\frac{D_1 D_2 D_3^{2/\lah}}{\lah-2\ell}B(\dM/2,\dM/2)\sum_{i,j=0}^{\efnj}\left(2\erry{\varepsilon+\dM}{0}+ \erry{\indcsi{j}-1+\dM}{\indhei{j}}
+\erry{\indcsi{i}-1+\dM}{\indhei{i}} \right)(t,x,y)\\
&\leq 2c_0^2\frac{D_1 D_2 D_3^{2/\lah}}{\lah-2\ell}(\efnj+1)^2B(\dM/2,\dM/2)\left(\erry{\varepsilon+\dM}{0}+\sum_{j=0}^{\efnj}\erry{\indcsi{j}-1+\dM}{\indhei{j}}\right)(t,x,y)\\
&\leq \lambda_1 \left(\erry{\varepsilon+\dM}{0}+\sum_{j=0}^{\efnj}\erry{\indcsi{j}-1+\dM}{\indhei{j}}\right)(t,x,y)\,.
\end{align*}

We assume that \eqref{ineq:q_n-bound} holds for $n\in\NN$.
By Lemma~\ref{lem:q_0-aux}, \eqref{ineq:q_n-bound}
and Corollary~\ref{cor-shifts} we have
\begin{align*}
|q_{n+1}(t,x,y)|
&\leq \int_0^t \int_{\Rd}|q_0(t-s,x,z)| |q_{n}(s,z,y)|\, dzds \\
&\leq c_0\lambda_n
\int_0^t \int_{\Rd}
\sum_{j=0}^{\efnj} \erry{\indcsi{j}-1}{\indhei{j}}(t-s,x,z)
\left(\erry{\varepsilon+n \dM}{0}+ \sum_{i=0}^{\efnj} \erry{\indcsi{i}-1+n\dM}{\indhei{i}}\right)(s,z,y)\,dzds\\
&\leq c_0\lambda_n D_1\sum_{j=0}^{\efnj}\int_{0}^{t}\int_{\Rd}\errx{\indcsi{j}-1}{\indhei{j}}(t-s,x,z)\erry{\varepsilon+n \dM}{0}(s,z,y)\,dzds\\
&\quad +c_0\lambda_n D_1\sum_{i,j=0}^{\efnj}\int_{0}^{t}\int_{\Rd}\errx{\indcsi{j}-1}{\indhei{j}}(t-s,x,z)\erry{\indcsi{i}-1+n\dM}{\indhei{i}}(s,z,y)\,dzds\eqqcolon I_1+I_2\,.
\end{align*}
We estimate $I_1$ by applying \eqref{ineq:aux2} and using
$D_3^{1/\lah}\leq D_3^{2/\lah}$, $(\efnj+1)\leq (\efnj+1)^2$ to obtain
\begin{align*}
I_1 &\leq c_0\lambda_n\frac{D_1 D_2 D_3^{1/\lah} D_4}{\lah-2\ell}\,
B\!\left(\frac{\dM}{2},\frac{(n+1)\dM}{2}\right)\!(\efnj+1)\!
\left(3\erry{\varepsilon+(n+1)\dM}{0}+\sum_{j=0}^{\efnj}\erry{\indcsi{j}-1+(n+1)\dM}{\indhei{j}}\right)\!(t,x,y)\\
&\leq \lambda_n\frac{\frac{3}{5}D_1 D_2 D_3^{2/\lah}D_4 D_5}{\lah-2\ell}\,
B\!\left(\frac{\dM}{2},\frac{(n+1)\dM}{2}\right)\left(\erry{\varepsilon+(n+1)\dM}{0}+\sum_{j=0}^{\efnj}\erry{\indcsi{j}-1+(n+1)\dM}{\indhei{j}}\right)(t,x,y)\\
	&=\frac{3}{5}\lambda_{n+1}\left(\erry{\varepsilon+(n+1)\dM}{0}+\sum_{j=0}^{\efnj}\erry{\indcsi{j}-1+(n+1)\dM}{\indhei{j}}\right)(t,x,y)\,.
\end{align*}
For $I_2$ we apply \eqref{ineq:aux1} to get
\begin{align*}
I_2
&\leq c_0\lambda_n \frac{D_1 D_2 D_3^{2/\lah}}{\lah-2\ell}\,
B\!\left(\frac{\dM}{2},\frac{(n+1)\dM}{2}\right)\!(N+1)^2
\!\left(2\erry{\varepsilon+(n+1)\dM}{0}+2\sum_{j=0}^{\efnj}\erry{\indcsi{j}-1+(n+1)\dM}{\indhei{j}}\right)\!(t,x,y)\\
&\leq \lambda_n\frac{\frac{2}{5}D_1 D_2 D_3^{2/\lah}D_4 D_5}{\lah-2\ell}\,
B\!\left(\frac{\dM}{2},\frac{(n+1)\dM}{2}\right)\left(\erry{\varepsilon+(n+1)\dM}{0}+\sum_{j=0}^{\efnj}\erry{\indcsi{j}-1+(n+1)\dM}{\indhei{j}}\right)(t,x,y)\\
&=\frac{2}{5}\lambda_{n+1}\left(\erry{\varepsilon+(n+1)\dM}{0}+\sum_{j=0}^{\efnj}\erry{\indcsi{j}-1+(n+1)\dM}{\indhei{j}}\right)(t,x,y)\,.
\end{align*}
Ultimately, after adding the estimates for $I_1$ and $I_2$ we get
\begin{equation*}
|q_{n+1}(t,x,y)|\leq \lambda_{n+1}\left(\erry{\varepsilon+(n+1)\dM}{0}+\sum_{j=0}^{\efnj}\erry{\indcsi{j}-1+(n+1)\dM}{\indhei{j}}\right)(t,x,y)\,.
\end{equation*}
This ends the proof of \eqref{ineq:q_n-bound}.
Using Remark~\ref{rem:r_t} 
we derive from
\eqref{ineq:q_n-bound} that
$$
| q_n(t,x,y) |\leq \lambda_n \,r_T^{n\dM}
\left(\erry{\varepsilon}{0}+ \sum_{j=0}^{\efnj} \erry{\indcsi{j}-1}{\indhei{j}}\right)(t,x,y)\,.
$$
It remains to recognize that
$\sum_{n=1}^{\infty} \lambda_n r_T^{n\dM}<\infty$, for instance by noticing
$$
\lambda_n \, r_T^{n\dM}
=c_0\Gamma(\dM/2)\left(\frac{c \, \Gamma(\dM/2) r_T^\dM}{\lah-2\ell}\right)^n\frac{1}{\Gamma((n+1)\dM/2)}\,.
$$
Actually, $r_T= 1$, $\dM$ depends only on $\varepsilon$ and $\boost$, while $\ell$ depends only on $\varepsilon+\dM$ and 
$\min\limits_{j=0,\ldots,\efnj} (\indcsi{j})$.

\qed

\begin{proposition}\label{prop:remainder-bound}
Assume $\Sa$ and \eqref{cond:pointwise}.
For every $\varepsilon\in (0,\boost)$
there exists a constant $c>0$ 
such that for all 
$t\in (0,\indTr]$ and $x,y\in\Rd$,
$$
|p(t,x,y)-p_0(t,x,y)|\leq ct \left(\erry{\varepsilon \land \varepsilon_0}{0}+\sum_{j=0}^{\efnj}\erry{\indcsi{j}-1}{\indhei{j}}\right)(t,x,y)\,.
$$
The constant $c$ can be chosen to depend only on
$d,\lmcc,\ka,\lah,\lch,h,\indnj,\boost, \varepsilon, \varepsilon_0$, $\min\limits_{j=0,\ldots,\efnj} (\indcsi{j})$.
\end{proposition}
\pf
We set $T=\indTr$.
By Proposition~\ref{prop:gen_est}, Lemma~\ref{lem:q-bound} and Corollary~\ref{cor-shifts} we 
have
\begin{align*}
|p(t,x,y)-p_0(t,x,y)| &\leq  c \int_0^t \int_{\Rd}(t-s) \errx{0}{0}(t-s,x,z)\erry{\varepsilon}{0}(s,z,y)\,dzds\\
&\quad+c\sum_{j=0}^{\efnj}\int_{0}^{t}\int_{\Rd}(t-s)\errx{0}{0}(t-s,x,z)\erry{\indcsi{j}-1}{\indhei{j}}(s,z,y)\,dzds\eqqcolon I_1+I_2 \,.
\end{align*}
Using Lemma~\ref{lem:conv}(c) with
\begin{align*}
\beta_0=0, \ \beta_1=0, \ \beta_2=0, \quad &n_1=n_2= m_1=m_2=0\,,\\
&\gamma_1=0, \ \gamma_2=\varepsilon\,,
\end{align*}
we get $I_1\leq (c/\lah) t\erry{\varepsilon}{0}(t,x,y)$.
Due to the definition of $\varepsilon_0$, see \eqref{def:ve_0-A} and \eqref{def:ve_0-A*}, for each $j=0,\ldots,\efnj$ there exists
$0<\ell_j < \lah \land (\efcs \indhei{j})$
such that $\varepsilon_0=\ell_j+\indcsi{j}-1$.
Let $\ell :=\max\{\ell_j\colon \, j=0,\ldots,\efnj\}$ and note that $0<\ell<\lah$.
We use Lemma~\ref{lem:conv}(c) with
\begin{align*}
\beta_0=\ell, \ \beta_1=0, \ \beta_2=\indhei{j}, \quad &n_1=n_2=\varepsilon_0+1-s_j\,,\\
&m_1=0, \ m_2=\varepsilon_0+1-s_j\,,\\
&\gamma_1=0, \ \gamma_2=\indcsi{j}-1\,,
\end{align*}
to get that
\begin{align*}
I_2\leq \frac{c}{\lah-\ell}\, t\left(\erry{\varepsilon_0}{0}+\sum_{j=0}^{\efnj}\erry{\indcsi{j}-1}{\indhei{j}}\right)(t,x,y)\,.
\end{align*}
This ends the proof.
\qed

\noindent
{\bf Proof of Theorem~\ref{thm:pointwise}.}
It suffices to use Proposition~\ref{prop:remainder-bound} and estimates of $p_0(t,x,y)$ provided by Proposition~\ref{prop:gen_est}, see Section~\ref{ssec:a-zoet}.
\qed

\appendix

\section{Proofs of Section~\ref{sec:f_a_a}}\label{sec:a-proofs}

\noindent
\textbf{Proof of Lemma~\ref{lem:time_conv}.}
We have by ($\aR$) that
$r_{t(1-u)}\leq \sqrt{1-u} \,r_t$ and
$r_{tu}\leq \sqrt{u} \,r_t$
for $t,u\in (0,1]$.
Then left hand side of the first inequality is
$$
t^{-1} \int_0^1 (1-u)^{-1} r_{t(1-u)}^{\varepsilon}\, u^{-1} r_{tu}^{k \varepsilon}\, \,du 
\leq t^{-1}r_t^{(k+1)\varepsilon} \int_0^1 (1-u)^{-1+\varepsilon/2} u^{-1+k\varepsilon/2}\,du\,,
$$
while of the second one is
$$
\int_0^1 u^{-1} r_{tu}^{\varepsilon}\,du \leq r_t^{\varepsilon} \int_0^1 u^{-1+\varepsilon/2}\,du\,.
$$
\qed

\noindent
\textbf{Proof of Fact~\ref{fact:Qn_well_Bb}.}
We will tacitly use the monotonicity of $r_t$.
Clearly,
the mapping
$(t,s_n,\ldots,s_1) \mapsto Q_{t-s_n}^0 \ldots Q_{s_1}^0 f$
is continuous 
in the norm
on $0<s_1<\ldots <s_n<t$.
Furthermore, by 
Lemma~\ref{lem:time_conv} we have
\begin{align}\label{ineq:Q_n_norm}
\idotsint\limits_{0<s_1<\ldots<s_n<t} \| Q_{t-s_n}^0 \ldots Q_{s_1}^0 f\|_{\infty}\,\, ds_1 \ldots ds_n 
&\leq  C_3^{n+1} \prod_{k=1}^n B\!\left(\frac{\varepsilon_0}{2},\frac{k\varepsilon_0}{2}\right) t^{-1} r_t^{(n+1)\varepsilon_0}\, \|f\|_{\infty}\,,
\end{align}
where $C_3$ and $\varepsilon_0$ are
taken from {\rm (H3)}.
Thus the Bocher integral is well defined in $(B_b(\Rd),\|\cdot\|_{\infty})$, see \cite[Proposition~1.2.2]{MR3617205}.
This proves (1) and (2).
Using Fubini's theorem
\cite[Proposition~1.2.7]{MR3617205} we get for $n=2,3\ldots$,
\begin{align*}
\idotsint\limits_{0<s_1<\ldots<s_n<t}\!\!\! Q_{t-s_n}^0 \ldots Q_{s_1}^0 f \,\,ds_1 \ldots ds_n 
=
\int_0^t \!\! Q_{t-s}^0 \! \left(\,\,\, \idotsint\limits_{0<s_1<\ldots<s_{n-1}<s}\!
Q_{s-s_{n-1}}^0 \ldots Q_{s_1}^0 f 
\,\,ds_1 \ldots d_{s_{n-1}} \!\! \right)\! ds,
\end{align*}
which gives \eqref{eq:Qn}.
Finally, we prove the continuity by induction.
For $n=0$ it holds by {\rm (H4)}.
Since
\begin{align*}
Q_t^n f - Q_{t_1}^n f
=  Q_t^0f- Q_{t_1}^0f+
\int_0^t Q_{t-s}^0\, Q_{s}^{n-1}f \,ds
 - \int_0^{t_1} Q_{t_1-s}^0\, Q_{s}^{n-1}f \,ds\,,
\end{align*}
and due to {\rm (H4)}, it suffices to consider
the second difference, which we also split 
into two parts given $\varepsilon_1 \in (0,t_1)$ and $|t-t_1|<\varepsilon_1/2$,
\begin{align*}
\int_0^t Q_{t-s}^0\, Q_{s}^{n-1}f \,ds
 - \int_0^{t_1} Q_{t_1-s}^0\, Q_{s}^{n-1}f \,ds
&=  \int_0^{t_1-\varepsilon_1} Q_{t-s}^0\, Q_{s}^{n-1}f \,ds
 - \int_0^{t_1-\varepsilon_1} Q_{t_1-s}^0\, Q_{s}^{n-1}f \,ds\\
& + 
\int_{t_1-\varepsilon_1}^t Q_{t-s}^0\, Q_{s}^{n-1}f \,ds
-\int_{t_1-\varepsilon_1}^{t_1} Q_{t_1-s}^0\, Q_{s}^{n-1}f \,ds\,.
\end{align*}
The first part tends to zero as $t\to t_1$
by the dominated convergence theorem 
\cite[Proposition~1.2.5]{MR3617205}, because
$\|Q_{t-s}^0 Q_s^{n-1} f\|_{\infty}\leq
c (t-s)^{-1} r_{t-s}^{\varepsilon_0}
\,s^{-1}r_s^{\varepsilon_0} \|f\|_{\infty}
\leq  c
(\varepsilon_1/2)^{-1} r_1^{\varepsilon_0}
s^{-1}r_s^{\varepsilon_0} \|f\|_{\infty}$
and the latter is integrable over $(0,t_1-\varepsilon_1)$, see Lemma~\ref{lem:time_conv}.
The continuity will follow if we show that
given $\varepsilon>0$ there is $\varepsilon_1\in (0,t_1/2)$ such that
for $|\tau-t_1|<\varepsilon_1/2$,
\begin{align*} 
\| \int_{t_1-\varepsilon_1}^{\tau} Q_{\tau-s}^0\, Q_{s}^{n-1}f \,ds \|_{\infty} <\varepsilon\,.
\end{align*}
Note that
$
\|Q_{\tau-s}^0\, Q_{s}^{n-1}f\|_{\infty}
\leq c (\tau-s)^{-1} r_{\tau-s}^{\varepsilon_0}\, (t_1/2)^{-1} r_1^{\varepsilon_0}\|f\|_{\infty}$,
therefore 
$\|\int_{t_1-\varepsilon_1}^{\tau}Q_{\tau-s}^0\, Q_{s}f \,ds\|_{\infty}\leq c \int_0^{(3/2)\varepsilon_1} u^{-1} r_u^{\varepsilon_0} du \,(t_1/2)^{-1} r_1^{\varepsilon_0}\|f\|_{\infty}$,
which can be made arbitrarily small by the choice of $\varepsilon_1$.
This  eventually ends the proof.
\qed

\noindent
\textbf{Proof of Fact~\ref{fact:Q_well_Bb}.}
Since
\begin{align*}
C_3^{n+1} \prod_{k=1}^n B\!\left(\frac{\varepsilon_0}{2},\frac{k\varepsilon_0}{2}\right) r_t^{(n+1)\varepsilon_0}
=\frac{\left[C_3\Gamma\!\left(\dfrac{\varepsilon_0}{2}\right)\right]^{n+1}}{\Gamma\!\left((n+1)\dfrac{\varepsilon_0}{2}\right)} \,\, r_t^{(n+1)\varepsilon_0}
\leq  
\frac{\left[C_3 \Gamma\!\left(\dfrac{\varepsilon_0}{2}\right)\right]^{n+1} r_1^{n \varepsilon_0}}{\Gamma\!\left((n+1)\dfrac{\varepsilon_0}{2}\right)}\,\, r_t^{\varepsilon_0}\,,
\end{align*}
where $\Gamma(a)$ is the Gamma function, and the latter is 
sumable over $n\in \NN$, 
by Fact~\ref{fact:Qn_well_Bb}
the series in \eqref{def:op_Q}
converges absolutely in the norm. 
Therefore, 
$Q_t\colon B_b(\Rd)\to B_b(\Rd)$
as well as 
\begin{align}\label{ineq:Q_norm}
\|Q_t f\|_{\infty}\leq \tilde{c}\, t^{-1} r_t^{\varepsilon_0} \|f\|_{\infty}\,,
\end{align}
 where $\tilde{c}=\tilde{c}(C_3,\varepsilon_0,r_1)$. 
Actually, the series
\eqref{def:op_Q} converges locally uniformly in $t\in (0,1]$, thus the continuity follows from that of $Q_t^nf$. 
Now by \eqref{eq:Qn},
\begin{align*}
\sum_{n=1}^N Q_t^nf=\int_0^t Q_{t-s}^0 \left(Q_s^0f  + \sum_{n=2}^N Q_s^{n-1} f \right)
= \int_0^t Q_{t-s}^0 \left(Q_s^0f  + \sum_{n=1}^{N-1} Q_s^{n} f \right).
\end{align*}
By the dominated convergence theorem \cite[Proposition~1.2.5]{MR3617205} we can pass with $N\to \infty$
under the integral,
which together with \eqref{def:op_Q} gives \eqref{eq:Q}.
\qed

\noindent
\textbf{Proof of Fact~\ref{fact:P_well_Bb}.}
Due to Fact~\ref{fact:Q_well_Bb}
the mapping $(t,s)\mapsto P_{t-s}^0 Q_s f$ is continuous in the norm on $0<s<t$. 
Lemma~\ref{lem:time_conv} provides
\begin{align}\label{ineq:R_bound}
\int_0^t \|P_{t-s}^0 Q_s f \|_{\infty}\, ds
\leq  C_1 c B(1,\varepsilon_0/2)\, r_t^{\varepsilon_0} \|f \|_{\infty}\,,
\end{align}
where $c=c(C_3,\varepsilon_0,r_1)$ is taken from Fact~\ref{fact:Q_well_Bb}, see \eqref{ineq:Q_norm}.
In particular, the Bocher integral is well defined in $(B_b(\Rd),\|\cdot\|_{\infty})$,
$P_t\colon B_b(\Rd)\to B_b(\Rd)$ and
\begin{align}\label{ineq:P_norm}
\|P_t f\|_{\infty} \leq \tilde{c} \|f\|_{\infty}\,,
\end{align}
where $\tilde{c}=\tilde{c}(C_1,C_3,\varepsilon_0,r_1)$.
To prove the continuity one writes the difference
\begin{align*}
P_t f - P_{t_1} f = P_t^0 f - P_{t_1}^0 f
+\int_0^t P_{t-s}^0\, Q_s f\, ds
-\int_0^{t_1} P_{t_1-s}^0\, Q_s f\, ds\,,
\end{align*}
then uses {\rm (H2)} for the first term,
and the dominated convergence theorem for the second (like in the proof of part (2) of Lemma~\ref{lem:continuity}). We skip further details. 
\qed

\noindent
\textbf{Proof of Fact~\ref{fact:strong-F}.}
For every $f\in B_b(\Rd)$ we treat $P_t^0 f$ and $Q_t^0f$ as elements of the Banach space $(C_b(\Rd),\|\cdot\|_{\infty})$. All  Bochner integrals in \eqref{def:op_P} and \eqref{def:op_Qn} exist, and the series
\eqref{def:op_Q} converges in that space, see \eqref{ineq:Q_n_norm},
\eqref{ineq:Q_norm},
\eqref{ineq:R_bound}. This ends the argumentation.
\qed

\noindent
\textbf{Proof of Fact~\ref{fact:PQ_well_C0}.}
The proof is exactly the same as that for Fact~\ref{fact:Q_well_Bb} and Fact~\ref{fact:P_well_Bb}. The only difference is that now all the objects are the elements of $C_0(\Rd)$, but the norm $\|\cdot\|_{\infty}$ remains the same.
\qed

\noindent
\textbf{Proof of Fact~\ref{fact:P-s-cont}.}
The results is a consequence of \eqref{def:op_P}, {\rm{(H0)}} and \eqref{ineq:R_bound} that give
$$
\lim_{t\to 0^+}\|\int_0^t P_{t-s}^0 Q_s f \,ds \|_{\infty} = 0\,.
$$
\qed

\noindent
\textbf{Proof of Lemma~\ref{lem:continuity}.}
Since all the domains of the given mappings are compact subsets of $\RR^2$, it suffices to show pointwise continuity.

The first statement follows from {\rm (H0)} and {\rm (H2)}. The third one is a consequence of {\rm (H4)}. For the fifth, note that if $(t,\varepsilon), (t_1,\varepsilon_1)\in [\tau,1]\times [0,1]$ the difference
$$
P_{\varepsilon}^0\,Q_t f - P_{\varepsilon_1}^0\,Q_{t_1} f =
P_{\varepsilon}^0 [Q_t f-Q_{t_1} f]  + P_{\varepsilon}^0\,Q_{t_1} f - P_{\varepsilon_1}^0\,Q_{t_1} f\,,
$$
tends to zero as $(t,\varepsilon)\to (t_1,\varepsilon_1)$, by
{\rm (H1)} and Lemma~\ref{fact:Q_well_Bb}
applied to the first term,
and {\rm (H0)} or {\rm (H2)}
applied to the second term.

Now, we prove (2).
For $(t,\varepsilon)\in \Omega_{[0,1]}$ we have
$$
\int_0^t P_{t-s+\varepsilon}^0\,Q_s f \,ds
=
\int_0^1 \ind_{0<s<t}\, P_{t-s+\varepsilon}^0\,Q_s f \,ds\,.
$$
Using {\rm (H2)} we get that
the integrand converges to
$\ind_{0<s<t_1}P_{t_1-s+\varepsilon_1}^0\,Q_s f$ in the norm,
whenever $(t,\varepsilon) \to (t_1,\varepsilon_1) \in \Omega_{[0,1]}$ and $s\neq t_1$.
This also includes the case $t_1=0$.
Further, by {\rm (H1)} and
Lemma~\ref{fact:Q_well_Bb},
$$
\| \ind_{0<s<t}\, P_{t-s+\varepsilon}^0\,Q_s f \|_{\infty}
\leq C_1 c\,s^{-1} r_s^{\varepsilon_0}\|f\|_{\infty}\,,
$$
and the latter is integrable by
Lemma~\ref{lem:time_conv}.
Therefore, the continuity follows
from the dominated convergence theorem 
\cite[Proposition~1.2.5]{MR3617205}.

Finally, we prove (4).
For $(t,\varepsilon), (t_1,\varepsilon_1)\in \Omega_{[\tau,1]}$,
given $\varepsilon_2 \in (0,t_1)$ and $|t-t_1|<\varepsilon_2/2$,
we have
\begin{align*}
\int_0^t Q_{t-s+\varepsilon}^0\, Q_{s}f \,ds
 - \int_0^{t_1} Q_{t_1-s+\varepsilon_1}^0\, Q_{s}f \,ds
&=  \int_0^{t_1-\varepsilon_2} Q_{t-s+\varepsilon}^0\, Q_{s}f \,ds
 - \int_0^{t_1-\varepsilon_2} Q_{t_1-s+\varepsilon_1}^0\, Q_{s}f \,ds\\
& + 
\int_{t_1-\varepsilon_2}^t Q_{t-s+\varepsilon}^0\, Q_{s}f \,ds
-\int_{t_1-\varepsilon_2}^{t_1} Q_{t_1-s+\varepsilon_1}^0\, Q_{s}f \,ds\,.
\end{align*}
The first part tends to zero as $(t,\varepsilon)\to (t_1,\varepsilon_1)$
by {\rm (H4)}
and
the dominated convergence theorem 
\cite[Proposition~1.2.5]{MR3617205}, which is justified since
by {\rm (H3)} and Fact~\ref{fact:Q_well_Bb}
we have
$\|Q_{t-s+\varepsilon}^0 Q_s f\|_{\infty}\leq
 C_3 c
(\varepsilon_2/2)^{-1} r_1^{\varepsilon_0}
s^{-1}r_s^{\varepsilon_0} \|f\|_{\infty}$
and the latter is integrable over $(0,t_1-\varepsilon_2)$, see Lemma~\ref{lem:time_conv}.
It remains to show that
given $\varepsilon_3>0$ there is $\varepsilon_2\in (0,t_1/2)$ such that
for 
all $(t,\varepsilon)\in \Omega_{[\tau,1]}$ satisfying
$|t-t_1|<\varepsilon_2/2$,
\begin{align*} 
\| \int_{t_1-\varepsilon_2}^t Q_{t-s+\varepsilon}^0\, Q_{s}f \,ds \|_{\infty} <\varepsilon_3\,.
\end{align*}
Note that
$
\|Q_{t-s+\varepsilon}^0\, Q_{s}f\|_{\infty}
\leq C_3 c (t-s+\varepsilon)^{-1} r_{t-s+\varepsilon}^{\varepsilon_0}\, (t_1/2)^{-1} r_1^{\varepsilon_0}\|f\|_{\infty}$,
therefore we get the bound
$\|\int_{t_1-\varepsilon_2}^tQ_{t-s+\varepsilon}^0\, Q_{s}f \,ds\|_{\infty}\leq C_3 c \int_{\varepsilon}^{\varepsilon+t-t_1+\varepsilon_2} u^{-1} r_u^{\varepsilon_0} du \,(t_1/2)^{-1} r_1^{\varepsilon_0}\|f\|_{\infty}$,
which can be made arbitrarily small by the choice of $\varepsilon_2$,
because $0< t-t_1+\varepsilon_2<(3/2)\varepsilon_2$,
see Lemma~\ref{lem:time_conv}.
\qed

\section{Shifts and convolutions}\label{sec:sac}

In the first observation we show that 
as soon as the integrals defining $\efdrf_r^x$ make sense, then \eqref{set:indrf-holder-s-B} extends in the spatial variable.

\begin{fact}\label{lem:drf-sp-ext}
For every $R>0$ there exists a constant
$c=c(R,\ka)$ such that
$$
t |\efdrf_{r_t}^{x}-\efdrf_{r_t}^{y}|\leq c \sum_{j=1}^{\indnj} |x-y|^{\indhei{j}}\, r_t^{\indcsi{j}}\,,
$$
for all $x,y\in\Rd$ satisfying $|x-y|\leq R$ and $t\in (0,\indTr]$.
\end{fact}
\pf
Since we assume $t\in (0,\indTr]$ we have $r_t\in (0,1]$.
Due to  \eqref{set:indrf-holder-s-B} we only need to discuss $1<|x-y|\leq R$. We take a sequence $x=w_0, w_1,\ldots,w_m=y$ such that $|w_i-w_{i-1}|\leq 1$.
Then, by the triangle inequality 
and \eqref{set:indrf-holder-s-B},
we have 
$t |\efdrf_{r_t}^{x}-\efdrf_{r_t}^{y}|
 \leq \sum_{j=1}^{\indnj} m \ka |x-y|^{\indhei{j}}\, r_t^{\indcsi{j}}$.
We can always take $m\leq R+1$.
\qed

We collect several properties of $r_t$ that are often used throughout the paper.
We only use  the definition of $h$ and
a fact that $\nu$ is a L{\'e}vy measure.

\begin{remark}\label{rem:r_t}
Note that $r_t$  is continuous and increasing in $t>0$,
and
for all $t>0$, $\lambda\in (0,1]$ we have
$r_{\lambda t} \leq \sqrt{\lambda}\, r_t$,
see \cite[Lemma~2.1]{MR4140542}.
For all $0<t\leq T$, $\beta \in [0,2]$
we have $r_t^{\beta} \leq r_T^2\vee 1$.
If~\eqref{set:h-scaling} holds, then for $r_t\in (0,1]$
we get
$r_t^{\lah} t^{-1}= r_t^{\lah} h(r_t)\geq \lch^{-1} h(1)$.
\end{remark}

In the next auxiliary observation we 
additionally use 
\eqref{set:k-bound} and
\eqref{set:indrf-cancellation-scale}, 
see also Remark~\ref{rem:exdrf}.

\begin{fact}\label{fact:efdrf-bound}
For every $T>0$ there exists 
a constant $c=c(T,\lmcc,\ka,\lch,h)$ such that for all 
$w\in\Rd$ and
$0<u \leq t\leq T$,
$$
u |\efdrf_{r_u}^w |\leq c r_t^{\efcs} \,.
$$
\end{fact}
\pf
If $r_u\in (0,1]$ it follows from \eqref{set:indrf-cancellation-scale} and Remark~\ref{rem:r_t}.
If $r_u \in (1,r_T]$, then we use
\mbox{$|\exdrf(x)|\leq \ka h(1)$},
we bound the integral part of
$u| \efdrf_{r_u}^w| $ by
$T r_u \ka \lmcc h(1)$,
and we apply
$r_u\leq r_T^{1-\efcs} r_t^{\efcs}$.

\qed

Now we deal with the effective drift $\efdrf_{r_t}^x$ as an argument of the bound function $\rr_t(x)$.
The following results are true if $\Sb$ or $\Sa$  holds, but in fact we only use
\eqref{set:h-scaling}--\eqref{set:k-bound}, \eqref{set:indrf-cancellation-scale}, 
\eqref{set:indrf-holder-s-B} or \eqref{set:indrf-holder-s-A},
and
$\efcs \indhei{j}+\indcsi{j}-1\geq 0$ for each $j=1,\ldots,\efnj$.

\begin{lemma}\label{lem:shifts-1} 
For every $T>0$ there exists a constant 
$c=c(T,\lmcc,\ka,\lah,\lch,h,\indnj)$
 such that for all $0<s<t\leq T$ and $x,y\in\Rd$,
\begin{align*}
\rr_t(y-x-(t-s)\efdrf_{r_{t-s}}^{x}-s\efdrf_{r_s}^{y})\leq c \rr_t(y-x-t\efdrf_{r_t}^y) \,.
\end{align*}
\end{lemma}
\pf
By Fact~\ref{fact:efdrf-bound} there is $M=M(T,\lmcc,\ka,\lch,h)$
such that $u|\efdrf_{r_u}^w| \leq M r_t^{\efcs}$ 
for all $w\in\Rd$ and $0<u\leq t\leq T$.
We consider three cases.

\noindent
{\it Case 1.} $r_t\in (1, r_T]$. 
Let $z=t\efdrf_{r_t}^y -(t-s)\efdrf_{r_{t-s}}^{x}-s\efdrf_{r_s}^{y}$.
Then 
$|z|\leq |t\efdrf_{r_t}^y|+|(t-s)\efdrf_{r_{t-s}}^{x}-s\efdrf_{r_s}^{y}|
\leq 3M  r_t^{\efcs}\leq 3M  r_t$, and by \cite[Corollary~5.10]{MR3996792}
we have
\begin{align*}
\rr_t(y-x-(t-s)\efdrf_{r_{t-s}}^{x}-s\efdrf_{r_s}^{y}-z+z) &\leq \rr_t(y-x-(t-s)\efdrf_{r_{t-s}}^{x}-s\efdrf_{r_s}^{y}-z)
\\ 
&= c \rr_t(y-x-t\efdrf_{r_t}^y) \,.
\end{align*}
\noindent
{\it Case 2.} $|y-x|\geq 5 M r_t^{\efcs}$. By the triangle inequality
\begin{align*}
|y-x-(t-s)\efdrf^x_{r_{t-s}}-s\efdrf^y_{r_s}|
&\geq
|y-x|
-|(t-s)\efdrf^x_{r_{t-s}}+s\efdrf^y_{r_s}|\\
&\geq 
\frac12 |y-x- t\efdrf^y_{r_t}|
+\frac12 |y-x|\\
&\quad-\frac12 t|\efdrf^y_{r_t}| 
-|(t-s)\efdrf^x_{r_{t-s}}+s\efdrf^y_{r_s}|\\
& \geq \frac12 |y-x- t\efdrf^y_{r_t}|
+\frac12 |y-x| - \frac52 M r_t^{\efcs}\\
&\geq \frac12 |y-x- t\efdrf^y_{r_t}|\,.
\end{align*}
Thus, by the monotonicity of $r^{-d}K(r)$ and scaling of $K(r)$, see  \cite[Lemma~5.1 (3) and (4)]{MR3996792},
\begin{align*}
\rr_t(y-x-(t-s)\efdrf^x_{r_{t-s}}-s\efdrf^y_{r_s})&\leq 
 \rr_t\left( \frac12(y-x-t\efdrf^y_{r_t})\right)
\leq
2^{d+2}\, \rr_t(y-x-t\efdrf^y_{r_t})\,.
\end{align*}

\noindent
{\it Case 3.} $r_t\in(0, 1]$ and $|y-x|< 5 M r_t^{\efcs}$. Due to \cite[Corollary~5.10]{MR3996792}
it suffices to show that there is 
$a=a(\lmcc,\ka,\lch,h,\indnj)$
 such that for all $0<s<t\leq \indTr $ and $x,y\in\Rd$,
\begin{align}\label{ineq:case3}
|t\efdrf^y_{r_t}
-(t-s)\efdrf^x_{r_{t-s}}-s\efdrf^y_{r_s}|\leq  a \,r_t\,.
\end{align}
We have
\begin{align*}
t\efdrf^y_{r_t}
-(t-s)\efdrf^x_{r_{t-s}}-s\efdrf^y_{r_s}
&=
(t-s)(\efdrf^y_{r_t}-\efdrf^x_{r_t})
+(t-s)(\efdrf^x_{r_t}-\efdrf^x_{r_{t-s}})
+ s(\efdrf^y_{r_t}-\efdrf^y_{r_s})\\
&=A_1+A_2+A_3\,.
\end{align*}
We estimate $A_2$ and $A_3$ using
\cite[(8)]{MR4140542}
 as follows
\begin{align*}
|A_2|+|A_3|\leq (t-s) r_t \, h(r_{t-s})
+s\, r_t\, h(r_s)
=2r_t\,.
\end{align*}
Now, by
 \eqref{set:indrf-holder-s-B}
(actually by Fact~\ref{lem:drf-sp-ext} with $R=5M$) or \eqref{set:indrf-holder-s-A}  we have
\begin{align*}
|A_1|
\leq t |\efdrf_{r_t}^y-\efdrf_{r_t}^x|
\leq c \sum_{j=1}^{\efnj} 
|y-x|^{\indhei{j}}\, r_t^{\indcsi{j}}
\leq c \sum_{j=1}^{\efnj} 
r_t^{\efcs \indhei{j}+\indcsi{j}}
= c \left(\sum_{j=1}^{\efnj} 
r_t^{\efcs \indhei{j}+\indcsi{j}-1}\right) r_t\,.
\end{align*}
Therefore, 
the expression in the bracket is bounded if
$\efcs \indhei{j}+\indcsi{j}-1\geq 0$ for each $j=1,\ldots,\efnj$.
It is guaranteed by
our assumptions.
This finally gives \eqref{ineq:case3},
and ends the proof.
\qed

By nearly the same proof as that of Lemma~\ref{lem:shifts-1} 
we can obtain a similar result with $t\efdrf^y_{r_t}$
replaced by $t\efdrf^x_{r_t}$ on the right hand side of the inequality. 
The only noticeable difference is in {\it Case 3.}, where we would have to consider
$$t\efdrf^x_{r_t}
-(t-s)\efdrf^x_{r_{t-s}}-s\efdrf^y_{r_s}
=
-s(\efdrf^y_{r_t}-\efdrf^x_{r_t})
+(t-s)(\efdrf^x_{r_t}-\efdrf^x_{r_{t-s}})
+ s(\efdrf^y_{r_t}-\efdrf^y_{r_s})\,.$$
We omit further details.
\begin{lemma}\label{lem:shifts-2}
For every $T>0$ there exists a constant 
$c=c(T,\lmcc,\ka,\lah,\lch,h,\indnj)$
such that for all  $0<s<t\leq T$ and $x,y\in \Rd$,
\begin{align*}
\rr_t(y-x-(t-s)\efdrf^x_{r_{t-s}}-s\efdrf^y_{r_s})\leq c \rr_t(y-x-t\efdrf^x_{r_t})\,.
\end{align*}
\end{lemma}

Recall that $r_t$ is continuous and $\lim_{t\to 0^+} r_t =0$, see
Remark~\ref{rem:r_t}. 
Thus, for every $w\in\Rd$,
$t \,\efdrf^w_{r_t}$ is continuous in $t>0$ and
we  have
$\lim_{t\to 0^+} t \,\efdrf^w_{r_t}=0$ by 
\eqref{set:indrf-cancellation-scale}.
Recall from \cite[Lemma~5.1]{MR3996792}
that $K(r)$ is continuous.
All the above allows us to pass with $s \to t^-$ in
Lemma~\ref{lem:shifts-2}. The outcome is formulated in the following corollary.
One could also give a direct proof without referring to continuity, again similar to that of Lemma~\ref{lem:shifts-1}.

\begin{corollary}\label{cor-shifts}
For every $T>0$ there exists a constant $c=c(T,\lmcc,\ka,\lah,\lch,h,\indnj)$ such that for all $t\in(0,T]$
and $x,y\in \Rd$,
\begin{align*}
\rr_t(y-x-t\efdrf^y_{r_t})\leq c \rr_t(y-x-t\efdrf^x_{r_t})\,.
\end{align*}
\end{corollary}

\begin{lemma}\label{lem:conv}
Let $\beta_0\in  [0,\lah)$.
\begin{itemize}
\item[(a)]
For every $T>0$  there exists a constant $c_1=c_1(d,T,\lmcc,\ka,\lch,h)$
such that for all $t\in(0,T]$, $x,y\in\Rd$  and $\beta\in [0,2]$,
\begin{align*}
\int_{\Rd} \errx{0}{\beta}(t,x,z)dz \leq \frac{c_1}{\lah-\beta_0}\, t^{-1} r_t^{\beta_0 \land \efcs \beta}\,,
\end{align*}
and
\begin{align*}
\int_{\Rd} \erry{0}{\beta}(t,z,y)dz \leq \frac{c_1}{\lah-\beta_0} \,t^{-1}r_t^{\beta_0 \land \efcs \beta}\,.
\end{align*}
\item[(b)]
For every $T>0$ there exists a constant 
$c_2=c_2(d,T,\lmcc,\ka,\lah,\lch,h,\indnj)$
such that for all
$\beta_1,\beta_2\in [0,1]$ and $n_1,n_2,m_1,m_2\in [0,\beta_0]$ with $n_1,n_2\leq \efcs(\beta_1+\beta_2)$, $m_1\leq \efcs \beta_1$, $m_2\leq \efcs\beta_2$
and all $0<s<t\leq T$, $x,y\in\Rd$,
\begin{align*}
\int_{\Rd} 
\errx{0}{\beta_1}(t-s,x,z)
\erry{0}{\beta_2}(s,z,y)dz
\leq \frac{c_2}{\lah-\beta_0} &\left[ 
\left((t-s)^{-1}r_{t-s}^{n_1}+s^{-1}r_s^{n_2}\right) \erry{0}{0}(t,x,y)\right. \\
&\quad + (t-s)^{-1} r_{t-s}^{m_1} \erry{0}{\beta_2}(t,x,y) \\
&\quad  +  \left.  s^{-1} r_{s}^{m_2} \erry{0}{\beta_1}(t,x,y) \right].
\end{align*}
\item[(c)]
Let $T>0$. For all $\gamma_1,\gamma_2\in\RR$,
$\beta_1,\beta_2, \theta_1,\theta_2\in [0,1]$ and $n_1,n_2,m_1,m_2\in [0,\beta_0]$ with $n_1,n_2\leq \efcs(\beta_1+\beta_2)$, $m_1\leq \efcs\beta_1$, $m_2\leq \efcs\beta_2$
and 
satisfying
\begin{align*}
\pi_i := 
\left( \frac{\gamma_i+n_i\land m_i}{2}\right) \land \left(\frac{\gamma_i+n_i\land m_i}{\lah}\right) +1-\theta_i >0\,,
\qquad i=1,2\,,
\end{align*}
and all $x,y\in\Rd$, $t\in (0,T]$ 
we have
\begin{align*}
&\int_0^t \int_{\Rd}
(t-s)^{1-\theta_1} \errx{\gamma_1}{\beta_1}(t-s,x,z)\,
 s^{1-\theta_2}\erry{\gamma_2}{\beta_2}(s,z,y)dzds\\
&\leq \frac{c_3}{\lah-\beta_0}\, t^{2-\theta_1-\theta_2}
\left( \erry{\gamma_1+\gamma_2+ n_1}{0}
+\erry{\gamma_1+\gamma_2+ n_2}{0}
+ \erry{\gamma_1+\gamma_2+ m_1}{\beta_2}
+\erry{\gamma_1+\gamma_2+ m_2}{\beta_1} \right)(t,x,y)\,,
\end{align*}
where 
\begin{align*}
c_3:= c_4\,
B(\pi_1 ,\pi_2)\,,
\qquad c_4:=c_2 [\lch (1 \vee r_T^2 )]^{-\frac{\gamma_1\land 0+\gamma_2\land 0}{\lah}}\,.
\end{align*}
\end{itemize}
\end{lemma}
\pf
We only show the first inequality in part (a), because the other one has exactly the same proof.
Note that 
$\left(|z-x|^{\beta}\land 1\right)\leq 2 \left(|z-x-t\efdrf_{r_t}^x|^{\beta}\land 1\right)+
2 |t\exdrf_{r_t}^x|^{\beta}
$. Thus, by
Fact~\ref{fact:efdrf-bound} we get
\begin{align*}
\errx{0}{\beta}(t,x,z)
&= \left(|z-x|^{\beta}\land 1\right) \err{0}{0}(t,z-x-t\efdrf_{r_t}^x)\\
&\leq 2 \err{0}{\beta_0 \land \efcs\beta}(t,z-x-t\efdrf_{r_t}^x)
+  2 c^{\beta} \,\err{\efcs \beta}{0}(t,z-x-t\efdrf_{r_t}^x)\,,
\end{align*}
where
\begin{align*}	
\err{\gamma}{\beta}(t,x):= r_t^{\gamma} \left(|x|^{\beta}\land 1\right) t^{-1} \rr_t(x)\,.
\end{align*}
Now, part (a) follows from \cite[Lemma~5.12]{MR3996792} and \cite[Remark~2.12]{MR4140542}.
Next,
applying \cite[Proposition~5.13]{MR3996792}
and Lemma~\ref{lem:shifts-1} we get
\begin{align*}
\errx{0}{0}(t-s,x,z)\land
\erry{0}{0}(s,z,y)
&\leq c t^{-1}\rr_t(y-x-(t-s)\efdrf_{r_{t-s}}^{x}-s\efdrf_{r_s}^{y})\\
&\leq c t^{-1}\rr_t(y-x-t\efdrf_{r_t}^y)\\
&= c \erry{0}{0}(t,x,y)\,.
\end{align*}
Hence
\begin{align*}
\frac{\errx{0}{0}(t-s,x,z)
\,\erry{0}{0}(s,z,y)}{c \erry{0}{0}(t,x,y)}
\leq \errx{0}{0}(t-s,x,z)
+\erry{0}{0}(s,z,y)\,,
\end{align*}
where 
$c=c(d,T,\lmcc,\ka,\lah,\lch,h,\indnj)$.
Recall that (see \cite[page 6059]{MR3996792})
\begin{align*}
&\left(|x-z|^{\beta_1} \land 1\right)
\left(|z-y|^{\beta_2} \land 1\right)
\leq \left(|x-z|^{\beta_1+\beta_2} \land 1\right)
+ \left(|x-z|^{\beta_1} \land 1\right) \left(|x-y|^{\beta_2} \land 1\right),\\
&\left(|x-z|^{\beta_1} \land 1\right)
\left(|z-y|^{\beta_2} \land 1\right)
\leq
\left(|z-y|^{\beta_1+\beta_2} \land 1\right)
+ \left(|x-y|^{\beta_1} \land 1\right) \left(|z-y|^{\beta_2} \land 1\right).
\end{align*}
Thus
\begin{align*}
&\frac{\errx{0}{\beta_1}(t-s,x,z)
\,\erry{0}{\beta_2}(s,z,y)}{c \erry{0}{0}(t,x,y)}\\
&\qquad\leq
\left(\left(|x-z|^{\beta_1+\beta_2} \land 1\right)
 + \left(|x-z|^{\beta_1} \land 1\right) \left(|x-y|^{\beta_2} \land 1\right)\right)
\errx{0}{0}(t-s,x,z)\\
&\qquad\quad+
\left(
\left(|z-y|^{\beta_1+\beta_2} \land 1\right)
+ \left(|x-y|^{\beta_1} \land 1\right) \left(|z-y|^{\beta_2} \land 1\right)\right)
\erry{0}{0}(s,z,y)\\
&\qquad\leq 
\errx{0}{\beta_1+\beta_2}(t-s,x,z)
+\left(|x-y|^{\beta_2} \land 1\right)\errx{0}{\beta_1}(t-s,x,z)\\
&\qquad\quad +\erry{0}{\beta_1+\beta_2}(s,z,y)
+\left(|x-y|^{\beta_1} \land 1\right)\erry{0}{\beta_2}(s,z,y)\,.
\end{align*}
Part (b) follows by integrating in $z$, using part (a) and adjusting parameters by Remark~\ref{rem:r_t}.
Now, multiply the inequality in (b) by
$$
(t-s)^{1-\theta_1} r_{t-s}^{\gamma_1} \,s^{1-\theta_2} r_s^{\gamma_2}\,,
$$
integrate in $s$ and use \cite[Lemma~5.15]{MR3996792} to obtain (c) with a constant
\begin{align*}
c_4\cdot
\max\left\{
 B\left(k_1+1-\theta_1,\frac{\gamma_2}{2}\land \frac{\gamma_2}{\lah}+2-\theta_2\right)\!;
B\left( \frac{\gamma_1}{2}\land \frac{\gamma_1}{\lah}+2-\theta_1, l_1+1-\theta_2 \right)\!; \right.\\
 \left.
 B\left(k_2+1-\theta_1,
\frac{\gamma_2}{2}\land \frac{\gamma_2}{\lah}+2-\theta_2\right)\!;
B\left( \frac{\gamma_1}{2}\land \frac{\gamma_1}{\lah}+2-\theta_1 , l_2+1-\theta_2\right) \right\},
\end{align*}
where
$k_1=(\frac{\gamma_1+ n_1}{2}\land \frac{\gamma_1+ n_1}{\lah})$,
$k_2=(\frac{\gamma_1+ m_1}{2}\land \frac{\gamma_1+ m_1}{\lah})$, 
$l_1=(\frac{\gamma_2+ n_2}{2}\land \frac{\gamma_1+ n_2}{\lah})$,
$l_2=(\frac{\gamma_2+ m_2}{2}\land \frac{\gamma_1+ m_2}{\lah})$.
Finally, use the monotonicity of the beta function.
\qed

In the next lemma we only use 
the definition of $h$,
a fact that $\nu$ is a L{\'e}vy measure.

\begin{lemma}\label{lem:drf}
Assume that $\lah\in(0,1)$ in \eqref{set:h-scaling}.
For all $r\in (0,1]$,
\begin{align*}
\int_{r \leq |z|<1} |z| \nu(|z|)dz \leq 
\frac{2 \lch}{1-\lah} r^{\lah}h(r)\,.
\end{align*}
\end{lemma}
\pf
We use \cite[(5)]{MR4140542} with $f(s)=s$ and
$N(dz)=\ind_{|z|<1} \nu(|z|)dz$, \cite[Lemma~2.1 6.]{MR4140542}, 
\eqref{set:h-scaling} and the assumption $\lah\in(0,1)$,
\begin{align*}
\int_{r \leq |z|<1} |z| \nu(|z|)dz 
&= \int_0^{\infty} \int_{|z|\geq r\vee s} \ind_{|z|<1}  \nu(|z|)dz ds
\leq \int_0^1 \int_{|z|\geq r\vee s}  \nu(|z|)dz ds\\
&\leq \int_0^1 h(r\vee s)ds
= r h(r) +\int_r^1 h(s)ds\\
&\leq r h(r) + \lch \int_r^1 (r/s)^{\lah} h(r)ds
\leq r h(r)+ \frac{\lch}{1-\lah} r^{\lah} h(r)\,.
\end{align*}
\qed

\section{Heat kernels with frozen coefficients}\label{sec:a-frozen}

We analyse the function
$p^{\mathfrak{K}_w}(t,x,y)$ as introduced in
Section~\ref{ssec:a-zoet}. 
The first result provides the 
upper estimates for that function and its derivatives with respect to spatial variable $x\in\Rd$.
We embark by assuming
\eqref{set:h-scaling}--\eqref{set:k-bound}.
To shorten the notation we let 
$
\param
$\label{param}
to represent $\lmcc,\ka,\lah,\lch,h$.

\begin{proposition}\label{prop:gen_est}
For every $T>0$ 
and $\bbbeta\in \mathbb{N}_0^d$
there exists a constant $c(d,T,\bbbeta,\param)$ such that for all $t\in (0,T]$ and $x,y,w\in\Rd$,
$$
|\partial_x^{\bbbeta} p^{\mathfrak{K}_w}(t,x,y)|\leq c\,
r_t^{-|\bbbeta|} \rr_t(y-x-t\efdrf^w_{r_t})\,.
$$
\end{proposition}
\pf 
The inequality follows from
\cite[Theorem~5.2]{GS-2020}
and \cite[Remark~2.12]{MR4140542}.
\qed

Here are a few more rather general properties.
\begin{remark}\label{rem:gen_not}
From
\cite[Lemma~6.1]{MR4140542}
we have
\begin{align}\label{eq:par_t_is_L}
\partial_t p^{\mathfrak{K}_w}(t,x,y)
= \pL^{\mathfrak{K}_w}_x p^{\mathfrak{K}_w}(t,x,y)\,.
\end{align}
Furthermore, by \cite[(96), (97)]{MR4140542}
the operators $\partial_t$, $\partial_x^{\bbbeta}$
and $\pL_x^{\mathfrak{K}_v}$
commute on $p^{\mathfrak{K}_w}(t,x,y)$.
\end{remark}

In Proposition~\ref{prop:gen_est_time} below, after two auxiliary results, we estimate  $\partial_t\, p^{\mathfrak{K}_w}(t,x,y)$. 

\begin{lemma}\label{lem:delta_est}
For every $T>0$ 
and $\bbbeta\in \mathbb{N}_0^d$
there exists a constant $c=c(d,T,\bbbeta,\param)$ such that
for all
$t\in (0,T]$, $x,y,z,w\in\Rd$,
\begin{align*}
|\partial_x^{\bbbeta} \delta_{r_t}^{\mathfrak{K}_w} (t,x,y;z)|
\leq c \left(\frac{|z|}{r_t}\right)^2 r_t ^{-|\bbbeta|} \rr_t(y-x-t\efdrf^w_{r_t})
\qquad \mbox{if}\quad |z|<r_t\,,
\end{align*}
\begin{align*}
| \partial_x^{\bbbeta} \delta_{r_t}^{\mathfrak{K}_w} (t,x,y;z)|
\leq c\, r_t ^{-|\bbbeta|} \left( \rr_t(y-x-z-t \efdrf^w_{r_t})
+ \rr_t(y-x-t\efdrf^w_{r_t})\right)
\qquad \mbox{if}\quad |z|\geq r_t\,.
\end{align*}
\end{lemma}
\pf
The second inequality follows directly from
Proposition~\ref{prop:gen_est}.
For the first one we
use
\begin{align}\label{ineq:increm}
|f(x+z)-f(x) - \left<z,\nabla f(x)\right>|
\leq |z|^2 \sum_{|\bbgamma|=2}
\int_0^1 \int_0^1 |\partial_x^{\bbgamma}f (x+\theta'\theta z)  | \, d\theta' d\theta\,.
\end{align}
It gives
$$
|\partial_x^{\bbbeta}\, \delta_{r_t}^{\mathfrak{K}_w} (t,x,y;z)|
\leq |z|^2 \sum_{|\bbgamma|=2} \int_0^1\int_0^1 | \partial_x^{\bbbeta+\bbgamma}\, p^{\mathfrak{K}_w} (t,x+\theta'\theta z,y)|\,.
$$
Then we apply 
Proposition~\ref{prop:gen_est}
with \cite[Corollary~5.10]{MR3996792}.
\qed

\begin{corollary}\label{cor:delta-est}
For every $T>0$ 
and $\bbbeta\in \mathbb{N}_0^d$
there exists a constant $c=c(d,T,\bbbeta,\param)$ such that
for all
$t\in (0,T]$, $x,y,w\in\Rd$,
\begin{align*}
\int_{\Rd} | \partial_x^{\bbbeta}\, \delta_{r_t}^{\mathfrak{K}_w} (t,x,y;z)|\, J(z)dz
\leq c\, t^{-1} r_t^{-|\bbbeta|} \rr_t(y-x-t\efdrf^w_{r_t}) \,.
\end{align*}
\end{corollary}
\pf
We split the integral into two regions: $|z|<r_t$ and $|z|\geq r_t$, and we apply respective estimates from Lemma~\ref{lem:delta_est}. On $|z|<r_t$ we further use \eqref{set:J} and \cite[Lemma~2.1 6.]{MR4140542}, on $|z|\geq r_t$ we use \cite[Lemma~2.1 6.]{MR4140542} and \cite[Lemma~5.9]{MR3996792}.
\qed

Now, we additionally use \eqref{set:indrf-cancellation-scale}.
\begin{proposition}\label{prop:gen_est_time}
For every $T>0$ 
and $\bbbeta\in \mathbb{N}_0^d$
there exists a constant $c(d,T,\bbbeta, \param)$ such that for all $t\in (0,T]$ and $x,y,w,v\in\Rd$,
$$
|\pL_x^{\mathfrak{K}_{v}} \partial_x^{\bbbeta}\, p^{\mathfrak{K}_w}(t,x,y)|\leq c\,
t^{-1} r_t^{-|\bbbeta|+ \efcs-1} \rr_t(y-x-t\efdrf^w_{r_t})\,,
$$
and
$$
|\partial_t\, \partial_x^{\bbbeta}\, p^{\mathfrak{K}_w}(t,x,y)|\leq c\,
t^{-1} r_t^{-|\bbbeta|+ \efcs-1} \rr_t(y-x-t\efdrf^w_{r_t})\,.
$$
\end{proposition}
\pf
First we apply
\eqref{eq:pL} with $r=r_t$ to get
\begin{align*}
\pL_x^{\mathfrak{K}_{v}} \partial_x^{\bbbeta}\, p^{\mathfrak{K}_{w}}(t,x,y)
=\efdrf_{r_t}^{v}\cdot \nabla_x \partial_x^{\bbbeta}\, p^{\mathfrak{K}_{w}}(t,x,y)
+
\int_{\Rd} \partial_x^{\bbbeta}\,\delta_{r_t}^{\mathfrak{K}_{w}}(t,x,y;z)\kappa(v,z)J(z)dz\,.
\end{align*}
Using
Corollary~\ref{cor:delta-est} we bound the integral part by $c\, t^{-1} r_t^{-|\bbbeta|} \rr_t(y-x-t\efdrf^w_{r_t})$.
Next, by
Fact~\ref{fact:efdrf-bound} and
Proposition~\ref{prop:gen_est},
\begin{align*}
|\efdrf_{r_t}^v\cdot \nabla_x \partial_x^{\bbbeta}\, p^{\mathfrak{K}_w}(t,x,y)|
\leq c\, t^{-1} r_t^{\efcs- 1-|\bbbeta|}
\rr_t(y-x-t\efdrf^w_{r_t})\,.
\end{align*}
We adjust powers using Remark~\ref{rem:r_t}.
The second inequality is now a special case, see \eqref{eq:par_t_is_L}.
\qed

The next results provide estimates for $\partial_t^2\, p^{\mathfrak{K}_w}(t,x,y)$.
\begin{lemma}\label{lem:delta_est-t}
For every $T>0$ 
there exists a constant $c=c(d,T,\param)$ such that for all
$t\in (0,T]$, $x,y,z,w\in\Rd$,
\begin{align*}
|\partial_t \, \delta_r^{\mathfrak{K}_w} (t,x,y;z)|_{_{r=r_t}}
\leq c \left(\frac{|z|}{r_t}\right)^2 t^{-1 }r_t^{\efcs-1} \rr_t(y-x-t\efdrf^w_{r_t})
\quad \mbox{if}\quad |z|<r_t\,,
\end{align*}
\begin{align*}
| \partial_t\, \delta_{r}^{\mathfrak{K}_w} (t,x,y;z)|_{_{r=r_t}}
\leq c\, t^{-1 }r_t^{\efcs-1} \left( \rr_t(y-x-z-t \efdrf^w_{r_t})
+ \rr_t(y-x-t\efdrf^w_{r_t})\right)
\quad \mbox{if}\quad |z|\geq r_t\,.
\end{align*}
\end{lemma}
\pf
The second inequality follows immediately from 
Proposition~\ref{prop:gen_est_time}. For the first one we use
\eqref{ineq:increm}, see Remark~\ref{rem:gen_not}, to get
\begin{align*}
|\partial_t\, \delta_r^{\mathfrak{K}_w} (t,x,y;z)|_{_{r=r_t}}
\leq |z|^2 \sum_{|\bbgamma|=2} \int_0^1\int_0^1 |
\partial_t\, \partial_x^{\bbgamma}\, p^{\mathfrak{K}_w} (t,x+\theta'\theta z,y)|\,.
\end{align*}
Now, we apply Proposition~\ref{prop:gen_est_time} with \cite[Corollary~5.10]{MR3996792}.
\qed

Similarly to how
Lemma~\ref{lem:delta_est}
leads to
Corollary~\ref{cor:delta-est}, now
Lemma~\ref{lem:delta_est-t} gives the following.
\begin{corollary}\label{cor:delta-est-t}
For every $T>0$ 
there exists a constant $c=c(d,T,\param)$ such that
for all
$t\in (0,T]$, $x,y,w\in\Rd$,
\begin{align*}
\int_{\Rd} |\partial_t\, \delta_{r}^{\mathfrak{K}_w} (t,x,y;z)|_{_{r=r_t}}\, J(z)dz
\leq c\, t^{-2} r_t^{\efcs-1} \rr_t(y-x-t\efdrf^w_{r_t}) \,.
\end{align*}
\end{corollary}

\begin{proposition}\label{prop:gen_est_time2}
For every $T>0$ there exists a constant $c=c(d,T,\param)$ such that
for all
$t\in (0,T]$, $x,y,w,v\in\Rd$,
$$
|\pL_x^{\mathfrak{K}_{v}} \partial_t\, p^{\mathfrak{K}_w}(t,x,y)|\leq c (t^{-1} r_t^{\efcs-1})^2\,  \rr_t(y-x-t\efdrf^w_{r_t}) \,. 
$$
\end{proposition}
\pf The prove is similar to the proof of Proposition~\ref{prop:gen_est_time}.
\qed

We collect further estimates that are handy when verifying integrability.

\begin{corollary}\label{cor:L1_uni_time}
For every $T>\tau>0$
there exists a constant $c>0$ such that
for all
$t\in [\tau,T]$, $x,y,z,w,v\in\Rd$, $k\in\{0,1\}$ and
$|\bbbeta|\leq 2$,
\begin{align}\label{ineq:L1_uni_time-1}
|\partial_t^k \,\partial_x^{\bbbeta}\,p^{\mathfrak{K}_w}(t,x,y)|
+|\pL_x^{\mathfrak{K}_{v}}\partial_t^k \,p^{\mathfrak{K}_w}(t,x,y)|
\leq c \rr_{\tau}(y-x)\,,
\end{align}
\begin{align}\label{ineq:L1_uni_time-2}
|\partial_t^k \,\delta_1^{\mathfrak{K}_w}(t,x,y;z)|
\leq c |z|^2 \,\rr_{\tau}(y-x)
\qquad \mbox{if}\quad |z|< 1\,,
\end{align}
\begin{align}\label{ineq:L1_uni_time-3}
|\partial_t^k \,\delta_1^{\mathfrak{K}_w}(t,x,y;z)|
\leq c \rr_{\tau}(y-x-z)+c \rr_{\tau}(y-x)
\qquad \mbox{if}\quad |z|\geq 1\,.
\end{align}
\end{corollary}
\pf
We repeatedly use that
$r_t$ is increasing and 
 $t^{-1}\rr_t(x)$ is non-increasing  in $t>0$.
By Fact~\ref{fact:efdrf-bound} we have
$|t\efdrf^w_{r_t}|\leq c r_t^{\efcs} \leq c (r_T^{\efcs} r_{\tau}^{-1})\, r_{\tau} $,
therefore using \cite[Corollary~5.10]{MR3996792} we 
get
$$
t^{-1}\rr_t(y-x-t\efdrf^w_{r_t}) \leq c  
 \tau^{-1} \rr_{\tau}(y-x)\,.
$$
This together with Propositions~\ref{prop:gen_est}, \ref{prop:gen_est_time} and~\ref{prop:gen_est_time2} already proves \eqref{ineq:L1_uni_time-1} and \eqref{ineq:L1_uni_time-3}.
For $|z|<1= (r_{\tau}^{-1})\, r_{\tau}$
we have by \eqref{ineq:increm} and \eqref{ineq:L1_uni_time-1}, see Remark~\ref{rem:gen_not},
\begin{align*}
|\partial_t^k \, \delta_1^{\mathfrak{K}_w} (t,x,y;z)|
&\leq |z|^2 \sum_{|\bbgamma|=2} \int_0^1\int_0^1 | \partial_t^k \, \partial_x^{\bbgamma}\,p^{\mathfrak{K}_w}(t,x+\theta'\theta z,y)|\,d\theta' d\theta\\
&\leq c |z|^2 \sum_{|\bbgamma|=2} \int_0^1\int_0^1
\rr_{\tau}(y-x-\theta'\theta z) \,d\theta' d\theta\,.
\end{align*}
Finally, it suffices to use \cite[Corollary~5.10]{MR3996792}.
\qed

In the next result we show a concentration of mass around the starting point $x\in\Rd$.
\begin{lemma}\label{lem:strong_at_zero-aux}
For every $\varepsilon>0$,
$$
\lim_{t\to 0^+} \sup_{x\in\Rd} \int_{|x-y|>\varepsilon} p^{\mathfrak{K}_y}(t,x,y)\,dy=0\,,
$$
and
$$
\lim_{t\to 0^+} \sup_{x\in\Rd} \int_{|x-y|>\varepsilon} p^{\mathfrak{K}_x}(t,x,y)\,dy=0\,.
$$
\end{lemma}
\pf
By
\eqref{set:indrf-cancellation-scale}
we have $\lim_{t\to 0^+}\sup_{w\in\Rd} |t \efdrf_{r_t}^w|=0$, hence $|t \efdrf_{r_t}^w|\leq \varepsilon/2$ and consequently
$|y-x-t \efdrf_{r_t}^w|\geq |y-x|/2$
for sufficiently small $t$ if $|y-x|>\varepsilon$. 
Then
using Proposition~\ref{prop:gen_est} 
and the monotonicity of $r^{-d}K(r)$
\cite[Lemma~5.1]{MR3996792}
we get
$$
p^{\mathfrak{K}_w}(t,x,y)
\leq c \frac{t K(|y-x-t \efdrf_{r_t}^w|)}{|y-x-t \efdrf_{r_t}^w|^d}
\leq c  \frac{t K(|y-x|/2)}{(|y-x|/2)^d}\,,\qquad w \in\Rd\,.
$$
The result follows
because $\int_{|z|>\varepsilon} K(|z|/2)(|z|/2)^{-d} \,dz<\infty$, see \cite[Lemma~2.2]{MR4140542}.
\qed

In the next few results we investigate the regularity of $p^{\mathfrak{K}_{w}}(t,x,y)$ with respect to $w\in\Rd$.
\begin{lemma}\label{lem:difference}
For all $t>0$, $x,y,w_1,w_2,v\in\Rd$ and $s\in (0,t)$,
\begin{align*}
\frac{d}{d s} \int_{\Rd} &p^{\mathfrak{K}_{w_1}}(s,x,z)\, p^{\mathfrak{K}_{w_2}}(t-s,z,y)\,dz\\
&= \int_{\Rd} \pL_x^{\mathfrak{K}_{w_1}}p^{\mathfrak{K}_{w_1}}(s,x,z) \, p^{\mathfrak{K}_{w_2}}(t-s,z,y)\,dz
- \int_{\Rd} p^{\mathfrak{K}_{w_1}}(s,x,z)\, \pL_z^{\mathfrak{K}_{w_2}} p^{\mathfrak{K}_{w_2}}(t-s,z,y) \,dz\,,
\end{align*}
and
\begin{align*}
\int_{\Rd} \pL^{\mathfrak{K}_v}_x p^{\mathfrak{K}_{w_1}}(s,x,z)
p^{\mathfrak{K}_{w_2}}(t-s,z,y)\,dz
= &\int_{\Rd} p^{\mathfrak{K}_{w_1}}(s,x,z) \, \pL_z^{\mathfrak{K}_v} p^{\mathfrak{K}_{w_2}}(t-s,z,y)\,dz\,.
\end{align*}
\end{lemma}
\pf
The first equality follows by differentiating under the integral sign, i.e.,
by passing 
to the limit in the difference quotient under the integral using the dominated convergence theorem, which is justified by 
\eqref{ineq:L1_uni_time-1}, and then by applying \eqref{eq:par_t_is_L}.
To prove the second equality 
we use \eqref{eq:pL} with $r=1$ to get
\begin{align*}
\int_{\Rd} \pL_x^{\mathfrak{K}_{v}}p^{\mathfrak{K}_{w_1}}(s,x,z) \, p^{\mathfrak{K}_{w_2}}(t-s,z,y)\,dz
=\int_{\Rd}
\exdrf(v)\cdot \nabla_x \, p^{\mathfrak{K}_{w_1}}(s,x,z) \, p^{\mathfrak{K}_{w_2}}(t-s,z,y)\,dz\\ 
+ \int_{\Rd}\int_{\Rd} \delta_1^{\mathfrak{K}_{w_1}}(t,x,z;\vartheta)\, p^{\mathfrak{K}_{w_2}}(t-s,z,y) \,\kappa(v,\vartheta)J(\vartheta)d\vartheta 
\,dz\,.
\end{align*}
Then we apply Fubini's theorem, what is justified by Corollary~\ref{cor:L1_uni_time},
\cite[Lemma~2.1 6.]{MR4140542} and \cite[Lemma~5.9]{MR3996792}, and finally we use
$p^{\mathfrak{K}_{w_1}}(t,x+\vartheta,z)=p^{\mathfrak{K}_{w_1}}(t,x,z-\vartheta)$
and
integration by parts combined with
$\nabla_x \, p^{\mathfrak{K}_{w_1}}(t,x,z)= - \nabla_z \, p^{\mathfrak{K}_{w_1}}(t,x,z)$,
cf. the proof of \cite[Lemma~2.10]{MR3996792}.
\qed

From now on we already assume
\eqref{set:h-scaling}--\eqref{set:k-holder}, 
\eqref{set:indrf-cancellation-scale},
\eqref{set:indrf-holder-s-B}. 
\begin{lemma}\label{lem:small_difference}
There exists a constant $c>0$ such that for all $t\in (0,\indTr]$, $x,y\in\Rd$ satisfying
$|x-y|\leq 1$ we have
$$
|p^{\mathfrak{K}_x}(t,x,y)- p^{\mathfrak{K}_y}(t,x,y)|
\leq c t \sum_{j=0}^{\efnj} 
\erry{\indcsi{j}-1}{\indhei{j}}(t,x,y)\,.
$$
\end{lemma}
\pf
By Lemma~\ref{lem:difference} and \eqref{eq:par_t_is_L} we have
\begin{align*}
p^{\mathfrak{K}_x}(t,x,y)-p^{\mathfrak{K}_y}(t,x,y)
&= 
\lim_{\varepsilon \to 0^+} \int_{\varepsilon}^{t/2}
\int_{\Rd} p^{\mathfrak{K}_x}(s,x,z) \left( \pL_z^{\mathfrak{K}_x} 
-  \pL_z^{\mathfrak{K}_y}\right) p^{\mathfrak{K}_y}(t-s,z,y)\,dzds\\
&+
\lim_{\varepsilon \to 0^+ } \int_{t/2}^{t-\varepsilon}
\int_{\Rd} \left( \pL_x^{\mathfrak{K}_x} 
-  \pL_x^{\mathfrak{K}_y}\right) p^{\mathfrak{K}_x}(s,x,z)\,  p^{\mathfrak{K}_y}(t-s,z,y)\,dzds
\,.
\end{align*}
We only deal with the first component, the proof for the other one is the same. 
Using \eqref{eq:pL} with $r=r_{t-s}$ in the first equality, then applying 
\eqref{set:indrf-holder-s-B},
Proposition~\ref{prop:gen_est},
\eqref{set:k-holder},
Corollary~\ref{cor:delta-est}
and Remark~\ref{rem:r_t},
we get
\begin{align*}
&\int_{\varepsilon}^{t/2}
\int_{\Rd}  p^{\mathfrak{K}_x}(s,x,z)\,  | \!\left( \pL_z^{\mathfrak{K}_x} 
-  \pL_z^{\mathfrak{K}_y}\right) p^{\mathfrak{K}_y}(t-s,z,y)|\,dzds\\
&= \int_{\varepsilon}^{t/2}
\int_{\Rd}  p^{\mathfrak{K}_x}(s,x,z)
\bigg| (\efdrf_{r_{t-s}}^{x}-\efdrf_{r_{t-s}}^{y}) \cdot \nabla_z p^{\mathfrak{K}_y}(t-s,z,y)\\
&\hspace{0.25\linewidth} + \int_{\Rd} \delta_{r_{t-s}}^{\mathfrak{K}_y} (t-s,z,y;w) (\kappa(x,w)-\kappa(y,w))J(w)dw\bigg|\, dzds\\
&\leq 
c \int_{\varepsilon}^{t/2}
\int_{\Rd} \rr_s(z-x-s\efdrf^x_{r_s})
 \left(\sum_{j=0}^{\efnj} |x-y|^{\indhei{j}}\, r_{t-s}^{\indcsi{j}-1}\right) (t-s)^{-1} \rr_{t-s}(y-z-(t-s)\efdrf^y_{r_{t-s}})\,dzds\\
&= c   \int_{\varepsilon}^{t/2}  \left(\sum_{j=0}^{\efnj} |x-y|^{\indhei{j}}\, r_{t-s}^{\indcsi{j}-1}\right) \int_{\Rd} s \errx{0}{0}(s,x,z) \erry{0}{0}(t-s,z,y)\,dzds  \\
&\leq 
c \left(\sum_{j=0}^{\efnj} |x-y|^{\indhei{j}}\, r_{t/2}^{\indcsi{j}-1}\right) 
\int_{\varepsilon}^{t/2}
\int_{\Rd} s \errx{0}{0}(s,x,z) \erry{0}{0}(t-s,z,y)\,dzds\,.
\end{align*}
It suffices to use $r_{t/2}\geq c r_t$, which follows from
\cite[Lemma~2.3]{MR4140542},
and 
\begin{align*}
\int_{\varepsilon}^{t/2}
\int_{\Rd} s \errx{0}{0}(s,x,z) \erry{0}{0}(t-s,z,y)\,dzds
\leq c\int_{\varepsilon}^{t/2}
(1+s(t-s)^{-1})ds \,\erry{0}{0}(t,x,y)\,,
\end{align*}
which results from
Lemma~\ref{lem:conv}(b).
\qed

We shall investigate $q_0(t,x,y)$, see \eqref{eq:delta-alt}.

\begin{lemma}\label{lem:q_0-aux*1}
For every $R>0$ there exists a constant $c=c(d,R,\param)>0$ such that for all $t\in (0,\indTr]$, $k\in\{0,1\}$ and $x,y\in\Rd$ satisfying $|x-y|\leq R$ we have
$$
|(\partial_t-\pL_x^{\mathfrak{K}_x}) \partial_t^k p^{\mathfrak{K}_y}(t,x,y)|
\leq c (t^{-1 }r_t^{\efcs-1})^k
\sum_{j=0}^{\efnj} \erry{\indcsi{j}-1}{\indhei{j}}(t,x,y)\,.
$$
\end{lemma}
\pf
We use Remark~\ref{rem:gen_not} and \eqref{eq:delta-alt} to get
\begin{align*}
-(\partial_t-\pL_x^{\mathfrak{K}_x})
\partial_t^k &p^{\mathfrak{K}_y}(t,x,y)
=(\pL_x^{\mathfrak{K}_x}-\pL_x^{\mathfrak{K}_y}) \partial_t^k p^{\mathfrak{K}_y}(t,x,y)\\
&=(\efdrf_{r}^{x}-\efdrf_{r}^{y}) \cdot \nabla_x \,\partial_t^k p^{\mathfrak{K}_y}(t,x,y)
+
\int_{\Rd}  \partial_t^k \delta_{r}^{\mathfrak{K}_y} (t,x,y;z) (\kappa(x,z)-\kappa(y,z))J(z)dz\,.
\end{align*}
Putting $r=r_t$, using  Corollaries~\ref{cor:delta-est} and~\ref{cor:delta-est-t},
\eqref{set:k-bound} and \eqref{set:k-holder}, we bound the integral part by
$$
c \ka(|x-y|^{\khe}\land 1)t^{-1} (t^{-1} r_t^{\efcs-1})^k \,\rr_{t}(y-x-t\efdrf^y_{r_t})
=c \ka (t^{-1} r_t^{\efcs-1})^k \erry{0}{\khe}(t,x,y)\,.
$$
Further, by Propositions~\ref{prop:gen_est} and~\ref{prop:gen_est_time} we get
\begin{align}
|(\efdrf_{r_t}^{x}-\efdrf_{r_t}^{y}) \cdot \nabla_x \,\partial_t^k \,p^{\mathfrak{K}_y}(t,x,y)|
&\leq 
c |\efdrf_{r_t}^{x}-\efdrf_{r_t}^{y}|\,
r_t^{-1} (t^{-1} r_t^{\efcs-1})^k \, \rr_{t}(y-x-t\efdrf^y_{r_t}) \nonumber \\
&=
c t|\efdrf_{r_t}^{x}-\efdrf_{r_t}^{y}|\,
(t^{-1} r_t^{\efcs-1})^k \,\erry{-1}{0}(t,x,y)\,. \label{ineq:aux3}
\end{align}
The result follows from
\eqref{set:indrf-holder-s-B}
(actually from Fact~\ref{lem:drf-sp-ext}).
\qed

\begin{lemma}\label{lem:q_0-aux}
If \eqref{set:indrf-holder-s-A} holds, then
there exists a constant $c=c(d,\param)>0$ such that for all $t\in (0,\indTr]$, $k\in\{0,1\}$ and $x,y\in\Rd$ we have
$$
|(\partial_t-\pL_x^{\mathfrak{K}_x}) \partial_t^k p^{\mathfrak{K}_y}(t,x,y)|
\leq c (t^{-1 }r_t^{\efcs-1})^k
\sum_{j=0}^{\efnj} \erry{\indcsi{j}-1}{\indhei{j}}(t,x,y)\,.
$$
\end{lemma}
\pf
The proof is basically the same as that of Lemma~\ref{lem:q_0-aux*1}, with the advantage here that we can bound \eqref{ineq:aux3} for all $x,y\in\Rd$.
\qed

The next result will be needed when considering \eqref{set:indrf-holder-s-B} instead of the stronger condition \eqref{set:indrf-holder-s-A}.

\begin{lemma}\label{lem:q_0-aux*2}
There exist $R=R(\lmcc,\ka,\lch,h)>0$ and $c=c(d,\param)>0$ such that
for all $t\in (0,\indTr]$, $k\in\{0,1\}$ we have
$$
\sup_{x\in\Rd} \int_{|x-y|\geq R} |(\partial_t-\pL_x^{\mathfrak{K}_x}) \partial_t^k p^{\mathfrak{K}_y}(t,x,y)|\,dy \leq c t^{-1}
(t^{-1 }r_t^{\efcs-1})^k\,
r_t^{\lah+\efcs-1}\,.
$$
\end{lemma}
\pf
The proof resembles that of Lemma~\ref{lem:strong_at_zero-aux}.
First, 
according to Propositions~\ref{prop:gen_est_time}
and~\ref{prop:gen_est_time2} we have
\begin{align*}
|(\partial_t-\pL_x^{\mathfrak{K}_x})\partial_t^k p^{\mathfrak{K}_y}(t,x,y)|
&\leq |\pL_x^{\mathfrak{K}_y}\partial_t^k
p^{\mathfrak{K}_y}(t,x,y)|
+|\pL_x^{\mathfrak{K}_x} \partial_t^k p^{\mathfrak{K}_y}(t,x,y)|\\
&\leq c
(t^{-1} r_t^{\efcs-1})^{1+k}\, \rr_t(y-x-t\efdrf^y_{r_t})\,.
\end{align*}
By 
\eqref{set:indrf-cancellation-scale}
we get $|t \efdrf_{r_t}^y|\leq c_0 r_t^{\efcs} \leq c_0$ and we define $R:=2c_0$.
Then $|y-x-t \efdrf_{r_t}^y|\geq |y-x|/2$
 if $|y-x|\geq R$.
Using the monotonicity of $r^{-d}K(r)$
\cite[Lemma~5.1]{MR3996792}
we get
$$
\rr_t(y-x-t\efdrf^y_{r_t})
\leq  \frac{t K(|y-x-t \efdrf_{r_t}^y|)}{|y-x-t \efdrf_{r_t}^y|^d}
\leq \frac{t K(|y-x|/2)}{(|y-x|/2)^d}\,.
$$
We bound $t\leq c r_t^{\lah}$ as in Remark~\ref{rem:r_t}.
The result follows
since $\int_{|z|\geq R} K(|z|/2)(|z|/2)^{-d} \,dz=2^{d-1}\omega_d h(R/2)$, see \cite[Lemma~2.2]{MR4140542}.
Here $\omega_d=2\pi^{d/2}/\Gamma(d/2)$ is the surface measure of the unit sphere in $\Rd$.
\qed

Finally, in the next result we require that
 $\Sb$ or $\Sa$ holds.

\begin{proposition}\label{prop:strong_at_zero}
We have
\begin{align*}
\lim_{t \to 0^+} \sup_{x\in\Rd} \left| \int_{\Rd} p^{\mathfrak{K}_y}(t,x,y)\,dy -1 \right|=0\,.
\end{align*}
\end{proposition}
\pf
We use Lemma~\ref{lem:small_difference},
\begin{align*}
&\left| \int_{\Rd} p^{\mathfrak{K}_y}(t,x,y)\,dy -1 \right|
\leq
 \int_{\Rd} \left| p^{\mathfrak{K}_y}(t,x,y)-p^{\mathfrak{K}_x}(t,x,y)  \right| dy\\
&\quad \leq c
\int_{|x-y|\leq 1} t \sum_{j=0}^{\efnj} 
\erry{\indcsi{j}-1}{\indhei{j}}(t,x,y)\,dy
+\int_{|x-y|> 1}
\left( p^{\mathfrak{K}_y}(t,x,y)+p^{\mathfrak{K}_x}(t,x,y)\right)dy\,.
\end{align*}
By Lemma~\ref{lem:strong_at_zero-aux} the second term converges to zero as desired.
We bound the first term using
Corollary~\ref{cor-shifts}
and Lemma~\ref{lem:conv}(a),
\begin{align*}
\int_{|x-y|\leq 1} t \sum_{j=0}^{\efnj} 
\erry{\indcsi{j}-1}{\indhei{j}}(t,x,y)\,dy
\leq c
\int_{\Rd} t \sum_{j=0}^{\efnj} 
\errx{\indcsi{j}-1}{\indhei{j}}(t,x,y)\,dy
\leq c \sum_{j=0}^{\efnj} r_t^{\beta_0 \land (\efcs \indhei{j})+ \indcsi{j}-1}\,,
\end{align*}
where due to our assumptions
$\beta_0\in [0,\lah)$ exists such that
$\beta_0 \land (\efcs \indhei{j})+ \indcsi{j}-1>0$ for all 
$j=0,\ldots, \efnj$.
It remains to apply $\sup_{x\in\Rd}$ and pass with $t$ to zero.
\qed

\begin{lemma}\label{lem:a-continuity}
The functions $\partial_x^{\bbbeta} p^{\mathfrak{K}_w}(t,x,y)$ and $\pL_x^{\mathfrak{K}_v}p^{\mathfrak{K}_w}(t,x,y)$ are jointly continuous in $t>0$, $x,y,w,v\in\Rd$.  
\end{lemma}
\pf
Due to \eqref{eq:pL} the continuity of $\pL_x^{\mathfrak{K}_v}p^{\mathfrak{K}_w}(t,x,y)$ will follow from  
that of 
$\exdrf(v)$, $\kappa(v,\cdot)$,
$\partial_x^{\bbbeta} p^{\mathfrak{K}_w}(t,x,y)$ and
the dominated convergence theorem,
which is justified by 
\eqref{set:J},
\eqref{set:k-bound}, 
\eqref{ineq:L1_uni_time-2},
\eqref{ineq:L1_uni_time-3}
and
\cite[Lemma~5.1 (8), Lemma~5.9]{MR3996792}.
Now, by \cite[(96) and (97)]{MR3996792},
$$
\partial_x^{\bbbeta} p^{\mathfrak{K}_w}(t,x,y)=(2\pi)^{-d}\int_{\Rd}i^{|\bbbeta|} z^{\bbbeta} e^{-i\left<y-x,z\right>} e^{-t \Psi_w(z)}\,dz\,,
$$
where $$\Psi_w(x)=-i\left<x,b(w)\right>-\int_{\Rd}\left(e^{i\left<x,z\right>}-1-i\left<x,z\right>\ind_{|z|<1}\right)\kappa(w,z)J(z)dz\,,$$
is continuous in $w\in\Rd$.
Note that 
${\rm Re}[\Psi_w(x)]\geq \ka^{-1} \lmcc^{-1}  \int_{\Rd}\left(1-\cos\left<x,z\right>\right)\nu(|z|)dz \geq c h(1/|x|)$,
where the last inequality follows from \cite[(85) and (86)]{MR3996792}.
Combining that with \eqref{set:h-scaling}
we have
$|e^{-t \Psi_w(z)}|\leq e^{- c |z|^{\lah}}$,
$|z|\geq 1$.
Therefore
we can use the dominated convergence theorem in the integral representation of $\partial_x^{\bbbeta} p^{\mathfrak{K}_w}(t,x,y)$ to obtain the continuity.
\qed

\section*{Acknowledgements}
The authors thank Krzysztof Bogdan and Alexei Kulik for helpful discussions and remarks.
We also thank the referee for a careful reading and comments that helped to improve the presentation of the paper's content.
The research was partially supported through the DFG-NCN Beethoven Classic 3 programme, contract no. 2018/31/G/ST1/02252 (National Science Center, Poland) and SCHI-419/11–1 (DFG, Germany).

\small
\bibliographystyle{abbrv}

\end{document}